\documentclass[a4,10pt]{article}
\usepackage{amsmath,amscd,amssymb}

\newcommand{\nbigd}{\mathcal{D}}
\newcommand{\nbige}{\mathcal{E}}
\newcommand{\nbigf}{\mathcal{F}}
\newcommand{\nbigg}{\mathcal{G}}
\newcommand{\nbigh}{\mathcal{H}}
\newcommand{\nbigi}{\mathcal{I}}
\newcommand{\nbigj}{\mathcal{J}}
\newcommand{\nbigk}{\mathcal{K}}
\newcommand{\nbigl}{\mathcal{L}}

\newcommand{\nbign}{\mathcal{N}}
\newcommand{\nbigo}{\mathcal{O}}
\newcommand{\nbigp}{\mathcal{P}}
\newcommand{\nbigq}{\mathcal{Q}}
\newcommand{\nbigr}{\mathcal{R}}
\newcommand{\nbigs}{\mathcal{S}}

\newcommand{\nbigu}{\mathcal{U}}
\newcommand{\nbigv}{\mathcal{V}}

\newcommand{\nbigx}{\mathcal{X}}

\newcommand{\proj}{\mathbb{P}}
\newcommand{\seisuu}{{\mathbb Z}}
\newcommand{\rnum}{{\mathbb Q}}

\newcommand{\cnum}{{\mathbb C}}
\newcommand{\real}{{\mathbb R}}

\newcommand{\DD}{\mathbb{D}}

\newcommand{\gbigd}{\mathfrak D}

\newcommand{\gbigh}{\mathfrak H}

\newcommand{\gbigp}{\mathfrak P}
\newcommand{\gbigq}{\mathfrak Q}
\newcommand{\gbigr}{\mathfrak R}
\newcommand{\gbigs}{\mathfrak S}

\newcommand{\gbigx}{\mathfrak X}

\newcommand{\gminia}{\mathfrak a}
\newcommand{\gminib}{\mathfrak b}

\newcommand{\gminiu}{\mathfrak u}
\newcommand{\gminiv}{\mathfrak v}

\newcommand{\vece}{{\boldsymbol e}}

\newcommand{\vecv}{{\boldsymbol v}}

\newcommand{\vecw}{{\boldsymbol w}}

\newcommand{\vecalpha}{{\boldsymbol \alpha}}
\newcommand{\veca}{{\boldsymbol a}}
\newcommand{\vecb}{{\boldsymbol b}}

\newcommand{\vecdelta}{{\boldsymbol \delta}}

\newcommand{\vecN}{{\boldsymbol N}}

\newcommand{\vecx}{{\boldsymbol x}}

\newcommand{\vecp}{{\boldsymbol p}}

\newcommand{\vecS}{{\boldsymbol S}}

\newcommand{\lrarr}{\longrightarrow}

\newcommand{\pf}{{\bf Proof}\hspace{.1in}}
\newcommand{\qed}{\mbox{\rule{1.2mm}{3mm}}}

\def\Hom{\mathop{\rm Hom}\nolimits}

\def\End{\mathop{\rm End}\nolimits}

\def\Re{\mathop{\rm Re}\nolimits}

\def\Gr{\mathop{\rm Gr}\nolimits}

\def\rank{\mathop{\rm rank}\nolimits}
\def\Spec{\mathop{\rm Spec}\nolimits}

\def\modulo{\mathop{\rm modulo}\nolimits}

\def\Gr{\mathop{\rm Gr}\nolimits}

\def\Res{\mathop{\rm Res}\nolimits}

\def\degpar{\mathop{\rm par\textrm{-}deg}\nolimits}

\def\can{\mathop{\rm can}\nolimits}

\def\id{\mathop{\rm id}\nolimits}

\def\gcd{\mathop{\rm g.c.d.}\nolimits}

\def\Irr{\mathop{\rm Irr}\nolimits}

\def\diag{\mathop{\rm diag}\nolimits}

\newcommand{\del}{\partial}
\newcommand{\delbar}{\overline{\del}}

\newcommand{\pardeg}{\degpar}

\newcommand{\nhom}{{\mathcal Hom}}

\newcommand{\sankaku}{\triangle}

\newcommand{\harmonicbundle}{(E,\delbar_E,\theta,h)}

\newcommand{\barz}{\overline{z}}
\newcommand{\zbar}{\barz}

\newcommand{\baralpha}{\overline{\alpha}}
\newcommand{\alphabar}{\baralpha}
\newcommand{\barlambda}{\overline{\lambda}}
\newcommand{\lambdabar}{\barlambda}
\newcommand{\fbar}{\overline{f}}

\newcommand{\nbigxlambda}{\nbigx^{\lambda}}

\newcommand{\Par}{{\mathcal Par}}

\newcommand{\lefttop}[1]{{}^{#1}\!}

\def\reg{\mathop{\rm reg}\nolimits}

\newcommand{\openclosed}[2]{]#1,#2]}

\newcommand{\vecvhat}{\widehat{\vecv}}

\newcommand{\Vhat}{\widehat{V}}
\newcommand{\nablahat}{\widehat{\nabla}}

\newcommand{\nablatilde}{\widetilde{\nabla}}
\newcommand{\vecvtilde}{\widetilde{\vecv}}
\newcommand{\vecwtilde}{\widetilde{\vecw}}

\newcommand{\gbar}{\overline{g}}

\newcommand{\kappatilde}{\widetilde{\kappa}}

\newcommand{\nbigetilde}{\widetilde{\nbige}}

\newcommand{\ptilde}{\widetilde{p}}

\newcommand{\phat}{\widehat{p}}

\newcommand{\ftilde}{\widetilde{f}}

\newcommand{\Phibar}{\overline{\Phi}}

\newcommand{\stilde}{\widetilde{s}}
\newcommand{\mubar}{\overline{\mu}}

\newcommand{\DDtilde}{\widetilde{\DD}}

\newcommand{\nbigehat}{\widehat{\nbige}}

\newcommand{\ellsitabar}{\underline{\ell}}

\newcommand{\vecy}{\boldsymbol y}

\def\Gal{\mathop{\rm Gal}\nolimits}

\def\op{\mathop{\rm op}\nolimits}

\def\hol{\mathop{\rm hol}\nolimits}

\def\St{\mathop{\rm St}\nolimits}

\newcommand{\vbar}{\overline{v}}
\newcommand{\ubar}{\overline{u}}

\newcommand{\Ztilde}{\widetilde{Z}}

\newcommand{\wtilde}{\widetilde{w}}

\newcommand{\Ttilde}{\widetilde{T}}

\newcommand{\betabar}{\overline{\beta}}

\newcommand{\Ktilde}{\widetilde{K}}

\newcommand{\TTtilde}{{T}\widetilde{T}}

\newcommand{\gminiatilde}{\widetilde{\gminia}}

\newcommand{\Vbar}{\overline{V}}

\newcommand{\nbigqtilde}{\widetilde{\nbigq}}

\newcommand{\kappabar}{\overline{\kappa}}

\newcommand{\inftyhat}{\widehat{\infty}}

\newcommand{\qbar}{\overline{q}}
\newcommand{\vecphat}{\widehat{\vecp}}

\newcommand{\vecptilde}{\widetilde{\vecp}}

\newcommand{\ttilde}{\widetilde{t}}
\newcommand{\nablacirc}{\nabla^{\circ}}
\newcommand{\nbigvcirc}{\nbigv^{\circ}}
\newcommand{\gminivtilde}{\widetilde{\gminiv}}

\newtheorem{thm}{Theorem}[section]
\newtheorem{cor}[thm]{Corollary}

\newtheorem{rem}[thm]{Remark}
\newtheorem{lem}[thm]{Lemma}
\newtheorem{prop}[thm]{Proposition}
\newtheorem{df}[thm]{Definition}

\begin{document}

\title{Harmonic bundles
and Toda lattices with opposite sign}

\author{Takuro Mochizuki}
\date{}
\maketitle

\begin{abstract}

We study a certain type of wild harmonic bundles
in relation with a Toda equation.
We explain how to obtain
a classification of the real valued solutions of the Toda equation
in terms of their parabolic weights,
from the viewpoint of the Kobayashi-Hitchin correspondence.
Then, we study the associated integrable variation of twistor structure.
In particular, we give a criterion for the existence of
an integral structure.
It follows from two results.
One is the explicit computation of the Stokes factors
of a certain meromorphic flat bundle.
The other is an explicit description of
the associated meromorphic flat bundle.
We use the opposite filtration of the limit mixed twistor structure
with an induced torus action.

\vspace{.1in}

\end{abstract}

\section{Introduction}
According to M. Guest and C.-S. Lin \cite{Guest-Lin},
the following equation is called 
the Toda lattice with opposite sign:
\begin{equation}
 \label{eq;12.10.7.1}
 2\delbar_z\del_zw_i
-e^{2(w_i-w_{i-1})}
+e^{2(w_{i+1}-w_i)}
=0
\quad
(i=1,\ldots,r)
\end{equation}
Here, 
we use the convention $w_{r+i}=w_i$ $(\forall i)$,
and $w_i$ are $\real$-valued functions
on $\cnum^{\ast}$.

S. Cecotti and C. Vafa found that the equation
(\ref{eq;12.10.7.1}) appeared from 
the $tt^{\ast}$-equation for several significant models
in physics
\cite{Cecotti-Vafa91}, \cite{Cecotti-Vafa93}.
Moreover, they observed that,
if $n=2$,
the equation (\ref{eq;12.10.7.1})
can be reduced to a Painlev\'{e} III equation
under some natural requirements.
Based on the analysis in
\cite{Its-Novokshenov} and \cite{McCoy-Tracy-Wu},
they mathematically verified the existence and 
uniqueness of the solutions adapted to the conditions
required by the models.
They found that the solutions have abundant information
of the models, and called them magical solutions.
Among their rich studies,
they proposed a problem to classify the solutions
of the equation whose associated $\frac{\infty}{2}$-variations
of Hodge structure have $\seisuu$-structure,
which are expected to be significant in physics.

As far as the author knows,
Guest and Lin are the first mathematicians
who systematically studied such issues.
They classified
the $\real$-valued global solutions of the equation (\ref{eq;12.10.7.1}),
by imposing an additional symmetry,
if it is reduced to the equation
for two unknown functions.
Moreover, collaborating with A. Its,
they studied the Stokes structure
of the associated $\frac{\infty}{2}$-VHS,
and they determined when it has 
a $\seisuu$-structure.

In this paper,
we study the general case.
We shall give a classification of
the $\real$-valued solutions
in terms of the parabolic weights,
and give a purely algebraic
criterion for the existence of $\seisuu$-structure.

As is well known, at least implicitly,
the equation (\ref{eq;12.10.7.1})
is closely related to the Hitchin equation 
for Higgs bundles on curves.
Indeed, the solutions of (\ref{eq;12.10.7.1})
are equivalent to
Toda-like harmonic bundles.
(See \S\ref{subsection;13.1.8.1}.)
Hence, it is quite natural 
to study the problem from the viewpoint of 
the Kobayashi-Hitchin correspondence
and the asymptotic behaviour 
of harmonic bundles.
We prove that the solutions (\ref{eq;12.10.7.1})
are classified by the parabolic weights.
(See Theorem \ref{thm;12.10.6.100} and 
the remark right after it.)
We also describe their asymptotic behaviour
and additional symmetry.

\begin{rem}
The harmonic bundles considered in this paper are
close to cyclotomic harmonic bundles of C. Simpson
in {\rm\cite{Simpson-middle-convolution}}.
The author found a more related work 
due to D. Baraglia {\rm\cite{Baraglia}},
who clarified the relation of his cyclic harmonic bundles
and affine Toda equations without poles
on smooth compact Riemann surfaces.
(See also {\rm\cite{Aldrovandi-Falqui}}.)
He has already used efficiently
the idea to deduce the symmetry of a harmonic metric
from the symmetry of the underlying Higgs bundle.
\end{rem}

\begin{rem}
There are many researches on the Toda equation.
However, at this moment,
the author does not know whether we can deduce
the classification of $\real$-valued solutions
{\rm(\ref{eq;12.10.7.1})}
from the previous results
(for example {\rm\cite{Tracy-Widom-1998}}, {\rm\cite{Takasaki}}
and {\rm\cite{Widom-1997}}).
Because our method gives only the existence and uniqueness,
it would be significant to obtain explicit descriptions
of the solutions by another method.
\end{rem}

As for $\seisuu$-structure,
our argument proceeds as follows.
First, we closely study the Stokes structure of
a certain meromorphic flat bundle (\ref{eq;12.10.2.1}).
In general, it is rather difficult to understand
the Stokes structure of a given meromorphic flat bundle.
But, in our case, it is not difficult.
The Stokes factors are related to the monodromy 
in a simple way,
and we obtain that the non-zero entries
of the Stokes factors are described
by the coefficients of the characteristic polynomial
of the monodromy.
In particular, we obtain that,
the meromorphic flat bundle has a $\seisuu$-structure,
if and only if the coefficients of the characteristic polynomial
are integers.

Second, to apply the result in the previous paragraph,
we need to know the precise form
of the meromorphic flat bundle
associated to a Toda-like harmonic bundle.
One of the important ingredients is
the limit mixed twistor structure
associated to the harmonic bundle with homogeneity.
It is naturally equipped with a torus action,
and hence it comes from a mixed Hodge structure.
(The mixed Hodge structure seems to be related with
that in \cite{Hertling-Sevenheck}.
See Remark \ref{rem;13.1.2.1}.)
Once it is obtained,
we can apply an argument of M. Saito
which is familiar in the construction of
Frobenius manifold.
(\cite{Saito-Brieskorn}.
 See \cite{Douai-Sabbah},
 \cite{Reichelt-logarithmic},
 \cite{Reichelt-Sevenheck}.)
Namely, by using the canonical decomposition
of the mixed Hodge structure due to P. Deligne,
we can find an opposite filtration to the Hodge filtration.
Then, we can extend the associated meromorphic flat bundle
to tr-TLE-structure in the sense of C. Hertling \cite{Hertling}.
Then, it is not difficult to deduce the more precise form
of the connection in the case of Toda-like harmonic bundle
by using its high symmetry.
Getting together the above two results,
we obtain the desired criterion.

\vspace{.1in}
It is reasonable to divide this paper into two parts,
Part I: the classification of $\real$-valued solutions
(\S\ref{section;12.12.16.2}--\ref{section;12.12.16.3}),
and 
Part II: the study on the Stokes structure
of the associated $\frac{\infty}{2}$-VHS
(\S\ref{section;12.10.8.10}--\ref{section;12.12.14.200}).
They are rather independent, although related,
and we use rather different techniques.

In \S\ref{section;12.12.16.2},
we give a review of 
a general theory of
filtered Higgs bundles and wild harmonic bundles
on curves for the convenience of the readers.
In \S\ref{section;12.12.16.3},
we apply it to a class of filtered Higgs bundles
and harmonic bundles on $(\proj^1,\{0,\infty\})$,
and we explain how to deduce a classification of
the $\real$-valued solutions of the Toda lattice with opposite sign
by parabolic weights.

In \S\ref{section;12.10.8.10},
after a review of $\seisuu$-structure of a meromorphic flat bundle
possibly with irregular singularity of pure slope,
we will study the Stokes structure
of a certain meromorphic flat bundle (\ref{eq;12.10.2.1}).
In particular, we give a criterion
for the existence of a $\seisuu$-structure
on the meromorphic flat bundle or its pull back
by a ramified covering.
In \S\ref{section;12.12.14.100},
we describe a rather general theory of
harmonic bundle with homogeneity
and the associated integrable variation of twistor structure,
which can be regarded
as a variant of 
\cite{Hertling},
\cite{mochi2},
\cite{sabbah2}
and \cite{s3}.
We would like to give details on the role of
Euler field in our situation.
We explain how a torus action is induced 
on the limit mixed twistor structure in some situation,
and how we obtain a tr-TLE-structure.
In \S\ref{section;12.12.14.200},
we apply the results in \S\ref{section;12.10.8.10}
and \S\ref{section;12.12.14.100}
to Toda-like harmonic bundles.
In \S\ref{subsection;12.12.16.11},
we give a complement on the tr-TLE structure
associated to a Toda-like harmonic bundle,
in particular, around 
the irregular singular point $\infty$.
We use it to obtain a precise description of
the associated meromorphic flat bundle.
Then, we obtain in \S\ref{subsection;12.12.16.12}
the criterion when they have $\seisuu$-structure.

\paragraph{Acknowledgement}

This study is inspired by the intriguing works
\cite{Cecotti-Vafa91},
\cite{Cecotti-Vafa93},
\cite{Dorfmeister-Guest-Rossman},
\cite{Guest-Lin},
\cite{Guest-Lin2}
and \cite{Guest-Lin-Its}.
I am heartily grateful to Martin Guest 
for his detailed and useful comments 
on the earlier version of this manuscript.
I also thank Hiroshi Iritani who attracted my attention
to the real structure of a quantum $D$-module.
I am grateful to Akira Ishii and Yoshifumi Tsuchimoto
for their constant encouragements.
I thank Claude Sabbah for numerous discussion and his kindness
for years.
A part of this paper is based on my talk in the conference
``Various Aspects on the Painlev\'{e} Equations''
held in November 2012.
It is my great pleasure to
express my gratitude to the organizers.
This study is supported by 
Grant-in-Aid for Scientific Research (C) 10315971.

\part{}

\section{Preliminary}
\label{section;12.12.16.2}
\subsection{Filtered Higgs bundle on a curve}

\subsubsection{Filtered Higgs bundle on a disc with a marked point}

We give a review of filtered Higgs bundle on curves.
We shall explain it in the case that
$X$ is a disc 
$\bigl\{z\in\cnum\,\big|\,|z|<1
 \bigr\}$
with a point $D=\{0\}$.
\paragraph{Filtered bundle}
A filtered bundle on $(X,D)$
is a locally free $\nbigo_{X}(\ast D)$-module
$\nbigv$ of finite rank
with an increasing filtration
by locally free $\nbigo_X$-submodules
$\nbigp_a\nbigv\subset\nbigv$ $(a\in\real)$
satisfying the following conditions,
(i) $\nbigp_{a}\nbigv_{|X\setminus D}
=\nbigv_{|X\setminus D}$,
(ii)
$\nbigp_a\nbigv=\bigcap_{b>a}\nbigp_b\nbigv$,
(iii)
$\nbigp_{a-1}\nbigv=
 \nbigp_a\nbigv\otimes\nbigo(-D)
=z\nbigp_a\nbigv$.
Let $\nbigp_{\ast}\nbigv$ denote
the $\nbigo_X(\ast D)$-module
$\nbigv$ with the filtration
$\bigl\{\nbigp_a\nbigv\,\big|\,
 a\in\real\bigr\}$,
and we call it a filtered bundle on $(X,D)$.
We say that $\nbigv$ is the underlying
$\nbigo_X(\ast D)$-module of $\nbigp_{\ast}\nbigv$,
and that $\nbigp_{\ast}\nbigv$ is a filtered bundle
over the $\nbigo_X(\ast D)$-module $\nbigv$.
The filtration is also denoted by $\nbigp_{\ast}\nbigv$.

We set $\nbigp_{<a}\nbigv:=\bigcup_{c<a}\nbigp_c\nbigv$
and $\Gr^{\nbigp}_a(\nbigv):=
 \nbigp_a\nbigv/\nbigp_{<a}\nbigv$.
It is naturally regarded as 
a finite dimensional vector space over $\cnum$.
We set
$\Par\bigl(
 \nbigp_{\ast}\nbigv
 \bigr)
:=\bigl\{
 a\in\real\,\big|\,
 \Gr^{\nbigp}_a(\nbigv)\neq 0
 \bigr\}$
which is discrete in $\real$.
If $a\in\Par(\nbigp_{\ast}\nbigv)$,
then $a+\seisuu\subset\Par(\nbigp_{\ast}\nbigv)$.

The sheaf $\nbigo_{X}(\ast D)$ is naturally a filtered bundle,
i.e.,
$\nbigp_a\nbigo_{X}(\ast D)
=\nbigo([a]D)$
for $a\in\real$,
where $[a]$ denote the integer such that
$a-1<[a]\leq a$.
The filtered bundle is denoted just by $\nbigo_{X}(\ast D)$.

A morphism of filtered bundles
$\psi:\nbigp_{\ast}\nbigv_1\lrarr\nbigp_{\ast}\nbigv_2$
is defined to be a morphism of $\nbigo_X(\ast D)$-modules
$\psi:\nbigv_1\lrarr\nbigv_2$
such that $\psi\bigl(\nbigp_a\nbigv_1\bigr)
\subset\nbigp_a\nbigv_2$
for any $a\in\real$.

Let $\nbigp_{\ast}\nbigv_i$ $(i=1,2)$
be filtered bundles on $(X,D)$.
The direct sum 
$\nbigp_{\ast}\nbigv_1\oplus\nbigp_{\ast}\nbigv_2$
is defined as the $\nbigo_{X}(\ast D)$-module
$\nbigv_1\oplus\nbigv_2$
with the filtration 
$\nbigp_a(\nbigv_1\oplus\nbigv_2)
:=\nbigp_{a}\nbigv_1\oplus\nbigp_{a}\nbigv_2$
$(a\in\real)$.
We naturally have
$\Gr^{\nbigp}_{a}\bigl(\nbigv_1\oplus\nbigv_2\bigr)
\simeq\Gr^{\nbigp}_a(\nbigv_1)
\oplus
 \Gr^{\nbigp}_a(\nbigv_2)$.
The tensor product
$\nbigp_{\ast}\nbigv_1\otimes\nbigp_{\ast}\nbigv_2$
is defined as the $\nbigo_{X}(\ast D)$-module
$\nbigv_1\otimes\nbigv_2$
with the filtration
$\nbigp_a(\nbigv_1\otimes\nbigv_2)$
$(a\in\real)$:
\[
 \nbigp_a(\nbigv_1\otimes \nbigv_2)
=\sum_{b+c\leq a}
 \nbigp_b\nbigv_1\otimes\nbigp_c\nbigv_2
=\sum_{\substack{b+c\leq a\\ -1<c\leq 0}}
 \nbigp_b\nbigv_1\otimes\nbigp_c\nbigv_2
\quad
(a\in\real)
\]
We define
$\nhom(\nbigp_{\ast}\nbigv_1,\nbigp_{\ast}\nbigv_2)$
as the $\nbigo_{X}(\ast D)$-module
$\nhom(\nbigv_1,\nbigv_2)$
with the filtration
$\nbigp_a\nhom(\nbigv_1,\nbigv_2)$
$(a\in\real)$:
\[
 \nbigp_a\nhom(\nbigv_1,\nbigv_2):=
 \Bigl\{
 f\in\nhom(\nbigv_1,\nbigv_2)\,\Big|\,
 f\bigl(
 \nbigp_b\nbigv_1
 \bigr)
 \subset
 \nbigp_{b+a}\nbigv_2
 \Bigr\}
\]
As a special case,
we have the dual $\nbigp_{\ast}\nbigv^{\lor}$
of $\nbigp_{\ast}\nbigv$
as $\nhom\bigl(\nbigp_{\ast}\nbigv,\nbigo_{X}(\ast D)\bigr)$.

\paragraph{Pull back, push-forward, descent}

We set $X_1:=\bigl\{w\in\cnum\,\big|\,|w|<1\bigr\}$
and $D_1:=\{w=0\}$,
and let $\varphi:X_1\lrarr X$
be given by $\varphi(w)=w^m$.
Let $\nbigp_{\ast}\nbigv$ be a filtered bundle on $(X,D)$.
We have an induced filtered bundle
$\varphi^{\ast}\nbigp_{\ast}\nbigv$,
called the pull back of
$\nbigp_{\ast}\nbigv$ by $\varphi$.
It is the $\nbigo_{X_1}(\ast D_1)$-module
$\varphi^{\ast}\nbigv$
with the filtrations
$\bigl(
 \nbigp_{a}\varphi^{\ast}\nbigv\,\big|\,
 a\in\real
 \bigr)$:
\[
 \nbigp_a\varphi^{\ast}\nbigv
=\sum_{mb+n\leq a}
 \varphi^{\ast}\nbigp_b\nbigv
\otimes\nbigo(n D_1)
\]

Let $\nbigp_{\ast}\nbigv_1$ 
be a filtered bundle on $(X_1,D_1)$.
We have an induced filtered bundle
$\varphi_{\ast}\nbigp_{\ast}\nbigv_1$
on $(X,D)$,
called the push-forward of $\nbigp_{\ast}\nbigv_1$
by $\varphi$.
It is the $\nbigo_X(\ast D)$-module
$\varphi_{\ast}\nbigv$
with the filtration
$\bigl(
 \nbigp_a(\varphi_{\ast}\nbigv)\,\big|\,
 a\in\real
 \bigr)$,
where 
$\nbigp_a(\varphi_{\ast}\nbigv_1)
=\varphi_{\ast}\nbigp_{am}\nbigv_1$.

Let $\Gal(\varphi)$ be the Galois group of the ramified covering
$\varphi$.
In this case, we have
$\Gal(\varphi)=\bigl\{t\in\cnum\,\big|\,t^m=1\bigr\}$
with the natural action $\rho$ on $X_1$
given by $\rho(t,w)=tw$.
If $\nbigp_{\ast}\nbigv_1$ is 
$\Gal(\varphi)$-equivariant,
then $\varphi_{\ast}\nbigp_{\ast}\nbigv_1$
is equipped with an induced $\Gal(\varphi)$-action.
The invariant part is called the descent of
$\nbigp_{\ast}\nbigv_1$.

If $\nbigp_{\ast}\nbigv$ is a filtered bundle on $(X,D)$,
the pull back $\varphi^{\ast}\nbigp_{\ast}\nbigv$
is naturally $\Gal(\varphi)$-equivariant,
and the descent of $\varphi^{\ast}\nbigp_{\ast}\nbigv$
is naturally isomorphic to $\nbigp_{\ast}\nbigv$.

\paragraph{Subobject}

Let $\nbigp_{\ast}\nbigv$ be a filtered bundle on $(X,D)$.
A subobject of $\nbigp_{\ast}\nbigv$
is an $\nbigo_X(\ast D)$-submodule $\nbigv'\subset\nbigv$
with a filtration $\{\nbigp_{a}\nbigv'\,|\,a\in\real\}$
such that $\nbigp_a\nbigv'\subset \nbigp_a\nbigv$.
It is called strict,
if $\nbigp_a\nbigv'=\nbigp_a\nbigv\cap\nbigv'$ for any $a$.
Note that if we are given an $\nbigo_X(\ast D)$-submodule
$\nbigv'\subset \nbigv$,
it has an induced filtration $\{\nbigp_a\nbigv'\}$
given by
$\nbigp_a\nbigv':=\nbigp_a\nbigv\cap\nbigv'$,
and such $\nbigp_{\ast}\nbigv'$ is strict.
In practice, we have only to consider such strict subobjects.

\begin{rem}
Let $\nbigp_{\ast}\nbigv$ be a filtered bundle on $(X,D)$.
Let $\nbigv=\bigoplus\nbigv_i$
be a direct sum of $\nbigo_X(\ast D)$-modules.
Then, we have the induced subobjects
$\nbigp_{\ast}\nbigv_i$ as above.
In general, we may have
$\nbigp_{\ast}\nbigv\not\simeq
 \bigoplus \nbigp_{\ast}\nbigv_i$.
If $\nbigp_{\ast}\nbigv\simeq
 \bigoplus \nbigp_{\ast}\nbigv_i$,
we say that the filtration $\nbigp_{\ast}\nbigv$
is compatible with the decomposition
$\nbigv=\bigoplus\nbigv_i$.
\end{rem}

\paragraph{Filtered Higgs bundle}

Let $\nbigp_{\ast}\nbigv$ be a filtered bundle on $(X,D)$.
A Higgs field of a filtered bundle $\nbigp_{\ast}\nbigv$ is 
an $\nbigo_X$-morphism
$\theta:\nbigv\lrarr \nbigv\otimes\Omega_{X}^1$.
The pair $(\nbigp_{\ast}\nbigv,\theta)$ is called
a filtered Higgs bundle on $(X,D)$.
We say that
$(\nbigv,\theta)$ is the underlying
meromorphic Higgs bundle of
$(\nbigp_{\ast}\nbigv,\theta)$,
and that
$(\nbigp_{\ast}\nbigv,\theta)$ is a 
filtered Higgs bundle over 
$(\nbigv,\theta)$.

A morphism of filtered Higgs bundles
$\psi:(\nbigp_{\ast}\nbigv_1,\theta_1)
 \lrarr(\nbigp_{\ast}\nbigv_2,\theta_2)$
is defined to be
a morphism of filtered bundles
$\psi:\nbigp_{\ast}\nbigv_1\lrarr\nbigp_{\ast}\nbigv_2$
such that $\psi\circ\theta_1=\theta_2\circ\psi$.
A subobject of $(\nbigp_{\ast}\nbigv,\theta)$
is a subobject $\nbigp_{\ast}\nbigv'\subset\nbigp_{\ast}\nbigv$
such that $\theta(\nbigv')\subset\nbigv'\otimes\Omega_X^1$.

Let $(\nbigp_{\ast}\nbigv_i,\theta_i)$ $(i=1,2)$ be 
filtered Higgs bundles on $(X,D)$.
We have the naturally induced Higgs fields
of $\theta^{(1)}$, $\theta^{(2)}$ and 
$\theta^{(3)}$
on 
$\nbigp_{\ast}\nbigv_1\oplus\nbigp_{\ast}\nbigv_2$,
$\nbigp_{\ast}\nbigv_1\otimes\nbigp_{\ast}\nbigv_2$
and 
$\nhom(\nbigp_{\ast}\nbigv_1,\nbigp_{\ast}\nbigv_2)$,
respectively:
$\theta^{(1)}(v_1,v_2)=\bigl(\theta_1(v_1),\theta_2(v_2)\bigr)$,
$\theta^{(2)}(v_1\otimes v_2)=
 \theta_1(v_1)\otimes v_2+v_1\otimes\theta_2(v_2)$,
$\theta^{(3)}(f)=\theta_2\circ f-f\circ\theta_1$.
In particular, we have the dual
$(\nbigp_{\ast}\nbigv,\theta)^{\lor}
=(\nbigp_{\ast}\nbigv,-\theta^{\lor})$
of a filtered Higgs bundle $(\nbigp_{\ast}\nbigv,\theta)$.

Let $\varphi:X_1\lrarr X$ be a morphism given by 
$\varphi(w)=w^m$.
For a filtered Higgs bundle
$(\nbigp_{\ast}\nbigv,\theta)$ on $(X,D)$,
we have its pull back
$\varphi^{\ast}(\nbigp_{\ast}\nbigv,\theta)
=(\varphi^{\ast}\nbigp_{\ast}\nbigv,\varphi^{\ast}\theta)$ 
on $(X_1,D_1))$.
For a filtered Higgs bundle
$(\nbigp_{\ast}\nbigv_1,\theta_1)$ on $(X_1,D_1)$,
we have its push-forward
$\varphi_{\ast}(\nbigp_{\ast}\nbigv_1,\theta_1)
=(\varphi_{\ast}\nbigp_{\ast}\nbigv_1,\varphi_{\ast}\theta_1)$
on $(X,D)$.
If $(\nbigp_{\ast}\nbigv_1,\theta_1)$ is 
$\Gal(\varphi)$-equivariant,
we have its descent
as the $\Gal(\varphi)$-invariant part of 
$\varphi_{\ast}(\nbigp_{\ast}\nbigv_1,\theta_1)$.

\paragraph{Good filtered Higgs bundle and regular filtered Higgs bundle}

Let $(\nbigp_{\ast}\nbigv,\theta)$ be a filtered Higgs bundle.
It is called regular, if 
$\theta\bigl(
 \nbigp_a\nbigv
 \bigr)
 \subset
 \nbigp_{a}\nbigv\otimes\Omega_{X}^1(\log D)
=\nbigp_{a+1}\nbigv\otimes\Omega_{X}^1$.
In that case, $\theta$ is called logarithmic.
A filtered Higgs bundle
$(\nbigp_{\ast}\nbigv,\theta)$
is called unramifiedly good,
if there exist a subset
$\nbigi\subset z^{-1}\cnum[z^{-1}]$
and a decomposition of filtered bundles
$\nbigp_{\ast}\nbigv
=\bigoplus_{\gminia\in\nbigi}
 \nbigp_{\ast}\nbigv_{\gminia}$
on a neighbourhood $\nbigu$ of $D$,
such that the following holds
for each $\gminia$:
\begin{itemize}
\item
 $\theta(\nbigv_{\gminia})
 \subset\nbigv_{\gminia}$.
 The restriction of $\theta$ to
 $\nbigv_{\gminia}$ is denoted by
 $\theta_{\gminia}$.
\item
 $\theta_{\gminia}-d\gminia\,\id_{\nbigv_{\gminia}}$
 is logarithmic
with respect to 
$\nbigp_{\ast}\nbigv_{\gminia}$.
\end{itemize}
A filtered Higgs bundle
$(\nbigp_{\ast}\nbigv,\theta)$
is called good,
if there exists a ramified covering
$\varphi:(X_1,D_1)\lrarr (X,D)$ given by
$\varphi(w)=w^m$
such that $\varphi^{\ast}(\nbigp_{\ast}\nbigv,\theta)$
is unramifiedly good.
The following lemma is easy to see.
\begin{lem}
\label{lem;12.9.24.40}
Let $(\nbigp_{\ast}\nbigv,\theta)$
be a good filtered Higgs bundle on $(X,D)$.
Assume that $\theta$ is tame, i.e.,
for the expression $\theta=f\,dz/z$,
the coefficients of the characteristic polynomial of $f$
are holomorphic at $0$.
Then, $(\nbigp_{\ast}\nbigv,\theta)$ is regular.
\end{lem}

\subsubsection{Filtered Higgs bundles on a curve}

Let $X$ be a complex curve
with a finite discrete subset $D\subset X$.
The notions of filtered bundle and filtered Higgs bundle
are naturally generalized in this global situation.
We take small neighbourhoods $U_P$ of any $P\in D$.
A filtered bundle $\nbigp_{\ast}\nbigv$ on $(X,D)$
is a locally free $\nbigo_X(\ast D)$-module $\nbigv$
of finite rank with the filtrations
$\nbigp^{(P)}_{\ast}(\nbigv_{|U_P})=
 \bigl(
 \nbigp^{(P)}_a(\nbigv_{|U_P})\,\big|\,
 a\in\real
 \bigr)$ 
by $\nbigo_{U_P}(\ast P)$-modules
for each $P\in U_P$,
such that
$\nbigp^{(P)}_{\ast}(\nbigv_{|U_P})$
is a filtered bundle on $(U_P,P)$.
For each $P\in D$,
we obtain the set
$\Par(\nbigp_{\ast}\nbigv,P):=
 \Par\bigl(
 \nbigp^{(P)}_{\ast}(\nbigv_{|U_P})
 \bigr)$.
For each $\veca=(a_P\,|\,P\in D)\in\real^D$,
we have an $\nbigo_X$-submodule
$\nbigp_{\veca}\nbigv\subset\nbigv$
determined by
$\nbigp_{\veca}\nbigv_{|U_P}
=\nbigp^{(P)}_{a_P}(\nbigv_{|U_P})$
for any $P\in D$.
Obviously,
the tuple
$\{\nbigp_{\veca}\nbigv\,|\,
 \veca\in\real^D\}$
determines $\nbigp_{\ast}\nbigv$.
A filtered Higgs bundle 
$(\nbigp_{\ast}\nbigv,\theta)$
is a filtered bundle $\nbigp_{\ast}\nbigv$
with a Higgs field $\theta$ of $\nbigv$.

Direct sum, tensor product
and inner homomorphism,
subobjects and morphisms are defined
as in the case of filtered bundles
on the disc.
A filtered Higgs bundle
$(\nbigp_{\ast}\nbigv,\theta)$ is called 
logarithmic
(resp. good)
if 
$(\nbigp_{\ast}\nbigv,\theta)_{|U_P}$
is logarithmic
(resp. good) for each $P\in D$.
The other operations are also naturally generalized
in the global case.

\subsubsection{Parabolic degree and stability condition}

Let $X$ be a projective curve with a finite subset $D\subset X$.
Let $\nbigp_{\ast}\nbigv$ be a filtered bundle.
Recall that the parabolic degree of
$\nbigp_{\ast}\nbigv$ is given as follows:
\begin{equation}
 \label{eq;12.10.3.3}
 \pardeg(\nbigp_{\ast}\nbigv):=
 \deg(\nbigp_{\veca}\nbigv)
-\sum_{P\in D}
 \sum_{\substack{
 c\in\Par(\nbigp_{\ast}\nbigv,P)\\
 a_P-1<c\leq a_P
 }}
 c\,\dim \Gr^{\nbigp^{(P)}}_{c}(\nbigv_{|U_P})
\end{equation}
It is independent of the choice of 
$\veca=(a_P\,|\,P\in D)\in \real^D$
in (\ref{eq;12.10.3.3}).
We define
$\mu(\nbigp_{\ast}\nbigv):=
\pardeg(\nbigp_{\ast}\nbigv)/\rank\nbigv$.

Let $(\nbigp_{\ast}\nbigv,\theta)$ be a filtered
Higgs bundle on $(X,D)$.
It is called stable,
if $\mu(\nbigp_{\ast}\nbigv_1)<\mu(\nbigp_{\ast}\nbigv)$
holds for any subobject
$(\nbigp_{\ast}\nbigv_1,\theta_1)
 \subset(\nbigp_{\ast}\nbigv,\theta)$
such that $0<\rank\nbigv_1<\rank\nbigv$.
It is called semistable,
if the condition 
``$\mu(\nbigp_{\ast}\nbigv_1)<\mu(\nbigp_{\ast}\nbigv)$''
is replaced with
``$\mu(\nbigp_{\ast}\nbigv_1)\leq\mu(\nbigp_{\ast}\nbigv)$''.
It is called poly-stable,
if we have a decomposition
 $(\nbigp_{\ast}\nbigv,\theta)=
 \bigoplus(\nbigp_{\ast}\nbigv_i,\theta_i)$
such that
each $(\nbigp_{\ast}\nbigv_i,\theta_i)$
is stable with 
$\mu(\nbigp_{\ast}\nbigv_i)=\mu(\nbigp_{\ast}\nbigv)$.
The following standard fact will be used implicitly.
\begin{lem}
\label{lem;12.12.21.2}
Suppose that $(\nbigp_{\ast}\nbigv_i,\theta_i)$ $(i=1,2)$
are stable.
The space of morphisms 
$(\nbigp_{\ast}\nbigv_1,\theta_1)
\lrarr
 (\nbigp_{\ast}\nbigv_2,\theta_2)$ is at most one dimensional,
and any non-zero morphism is an isomorphism.
\end{lem}

\subsection{Harmonic bundles on curves}

Let $Y$ be a complex curve.
Let $(E,\delbar_E)$ be a holomorphic vector bundle
on $Y$.
Let $\theta$ be a Higgs field of $(E,\delbar_E)$.
Let $h$ be a hermitian metric of $E$.
We have the Chern connection,
whose $(1,0)$-part is denoted by $\del_E$.
Let $\theta^{\dagger}$ denote the adjoint of $\theta$.
The metric $h$ is called  harmonic,
if the connection
$\nabla^1:=\delbar_E+\theta^{\dagger}+\del_E+\theta$
is flat.
In that case,
$(E,\delbar_E,\theta,h)$ is called a harmonic bundle.
Note that 
$(E,\delbar_E,\alpha\theta,h)$ is also 
a harmonic bundle if $\alpha$ is a complex number
with $|\alpha|=1$.
If the base space is higher dimensional,
such a metric is called pluri-harmonic.

\subsubsection{Prolongation of a harmonic bundle on a punctured disc}
\label{subsection;12.9.23.1}

Let $X:=\bigl\{z\in\cnum\,\big|\,|z|<1\bigr\}$
and $D:=\{0\}$.
Let us consider a harmonic bundle
$\harmonicbundle$ on $X\setminus D$.
We have the expression $\theta=f\,dz/z$
where $f$ is a holomorphic section of 
$\End(E)$.
We have the characteristic polynomial
$P=\det(t\id_E-f)=\sum a_j(z)t^j$.
The harmonic bundle is called tame
if $a_j(z)$ are holomorphic on $X$,
and it is called wild if $a_j(z)$ are meromorphic on $X$.
It is called unramifiedly good wild,
if there exist a subset $\nbigi\subset z^{-1}\cnum[z^{-1}]$
and a decomposition
$(E,\theta)=
 \bigoplus_{\gminia\in\nbigi}
 (E_{\gminia},\theta_{\gminia})$
such that $\theta_{\gminia}-d\gminia\,\id_{E_{\gminia}}$
are tame for each $\gminia\in z^{-1}\cnum[z^{-1}]$.
Note that,
if $\harmonicbundle$ is wild,
there exists a ramified covering
$\varphi:(X_1,D_1)\lrarr (X,D)$,
such that
$\varphi^{-1}\harmonicbundle$
is unramifiedly good wild.

\vspace{.1in}

Let $\harmonicbundle$ be a wild harmonic bundle 
on $(X,D)$.
For any $a\in\real$,
we have the $\nbigo_X$-module $\nbigp_aE$
given as follows.
Let $U$ be an open subset of $X$.
If $0\not\in U$,
let $\nbigp_aE(U)$ denote the space of
holomorphic sections of $E$ on $U$.
If $0\in U$,
let $\nbigp_aE(U)$ denote the space of
holomorphic sections $f$ of $E$ on $U\setminus D$
such that
$|f|_h=O(|z|^{-a-\epsilon})$ for any $\epsilon>0$.
We define
$\nbigp E:=\bigcup_{a\in\real}\nbigp_{a}E$,
which is an $\nbigo_X(\ast D)$-module.
The following proposition was proved by Simpson
\cite{Simpson90} in the tame case,
and by the author in the general case
\cite{Mochizuki-wild}.

\begin{prop}\mbox{{}}
$(\nbigp_{\ast}E,\theta)$ is a good filtered Higgs bundle on
$(X,D)$.
If the harmonic bundle is tame,
$(\nbigp_{\ast}E,\theta)$ is regular.
If the harmonic bundle is unramifiedly good wild,
then $(\nbigp_{\ast}E,\theta)$ is unramifiedly good.
\end{prop}

Let $(\nbigp_{\ast}\nbigv,\theta)$ be a filtered Higgs bundle
on $(X,D)$.
Let $h$ be a harmonic metric of 
$(E,\theta):=(\nbigv,\theta)_{|X\setminus D}$.
We obtain a filtered bundle
$\nbigp_{\ast}E$ as above.
If the identity $E\simeq \nbigv_{|X\setminus D}$ 
on $X\setminus D$
is extended to 
an isomorphism
$\nbigp_{\ast}E\simeq\nbigp_{\ast}\nbigv$
on $(X,D)$,
we say that $h$ is adapted to
the filtered bundle $\nbigp_{\ast}\nbigv$.

\subsubsection{Harmonic bundles on projective curves}

Let $X$ be a projective curve with a finite subset
$D\subset X$.
We take small neighbourhoods $U_P$ of any $P\in D$.
Let $\harmonicbundle$ be a harmonic bundle
on $X\setminus D$.
It is called tame (resp. wild, unramifiedly good wild)
if $\harmonicbundle_{|U_P\setminus P}$ is tame
(resp. wild, unramifiedly good wild)
for any $P\in D$.
By applying the procedure in \S\ref{subsection;12.9.23.1}
to each $\harmonicbundle_{|U_P\setminus P}$,
we obtain a good filtered Higgs bundle
$(\nbigp_{\ast}E,\theta)$ on $(X,D)$.

\begin{prop}\mbox{{}}
$(\nbigp_{\ast}E,\theta)$ is polystable
and satisfies
$\pardeg(\nbigp_{\ast}E)=0$.
\end{prop}
\pf
The tame case was proved by Simpson
(Theorem 5 in \cite{Simpson90}).
The general case can be shown similarly.
\hfill\qed

\begin{prop}
Let $(\nbigp_{\ast}\nbigv,\theta)$
be a filtered Higgs bundle on $(X,D)$.
Let $h_j$ $(j=1,2)$ be two harmonic metrics of 
$(E,\theta):=
 (\nbigv,\theta)_{|X\setminus D}$
adapted to $\nbigp_{\ast}\nbigv$.
Then, there exists a decomposition
$(\nbigp_{\ast}\nbigv,\theta)
=\bigoplus_{i\in I}(\nbigp_{\ast}\nbigv_i,\theta_i)$
such that the following holds:
\begin{itemize}
\item
 We set $(E_i,\theta_i):=(\nbigv_i,\theta_i)_{|X\setminus D}$.
 Then, the decomposition
 $(E,\theta)=\bigoplus (E_i,\theta_i)$
 is orthogonal with respect to both of $h_j$ $(j=1,2)$.
\item
 There exist positive numbers $a_i$
 such that $h_{1|E_i}=a_i\,h_{2|E_i}$.
\end{itemize}
In particular,
if $(\nbigp_{\ast}\nbigv,\theta)$ is stable,
the adapted harmonic metric is unique up to
the multiplication of positive constants.
\end{prop}
\pf
The argument for the tame case is 
given in \S2.2 of \cite{mochi4}, for example.
The general case can be shown similarly.
\hfill\qed

\vspace{.1in}

Conversely, the following theorem holds.
The tame case is due to Simpson \cite{Simpson90}.
The wild case is due to Biquard-Boalch \cite{biquard-boalch},
although some inessential assumption is imposed 
on the residues.
We may also prove it as in \cite{Mochizuki-wild}
where a similar statement is established
for filtered flat bundles by directly 
using the method in \cite{Simpson90}.
We are planning to give more details
elsewhere,
including the higher dimensional case.

\begin{thm}
\label{thm;12.9.24.50}
Let $(\nbigp_{\ast}\nbigv,\theta)$
be a $\mu$-stable Higgs bundle
on $(X,D)$
with $\mu(\nbigp_{\ast}\nbigv)=0$.
Then, there exists a harmonic metric $h$ of
$(E,\theta):=(\nbigv,\theta)_{|X\setminus D}$
such that 
$\nbigp_{\veca}E=\nbigp_{\veca}\nbigv$
for any $\veca\in\real^D$.
Such a metric is unique up to 
the multiplication of a positive constant.
\end{thm}
We can apply it if the Higgs bundle is irreducible,
for example.

\subsubsection{Symmetry}

By the uniqueness in Kobayashi-Hitchin correspondence,
the symmetry of a filtered Higgs bundle
is inherited to a harmonic bundle.
For the explanation,
we give a typical statement.
Let $X$ be a smooth projective curve
with a finite subset $D\subset X$.
Let $f$ be a holomorphic automorphism 
of $(X,D)$.
Let $(\nbigp_{\ast}\nbigv,\theta)$ be a stable good 
filtered Higgs bundle on $(X,D)$,
with an isomorphism
$\ftilde:f^{\ast}\nbigp_{\ast}\nbigv\simeq\nbigp_{\ast}\nbigv$.
Suppose
$f^{\ast}\theta=\alpha\theta$ 
for a complex number $\alpha$ with $|\alpha|=1$
under the isomorphism.

\begin{prop}
There exists $a>0$
such that
$f^{\ast}h=a\,h$
under $\ftilde$.
\end{prop}
\pf
Put $(E,\theta):=(\nbigv,\theta)_{|X\setminus D}$.
Because $h$ and $f^{\ast}h$ 
are harmonic metrics
of $(E,\delbar_E,\alpha\theta)$
adapted to $\nbigp_{\ast}\nbigv$,
we obtain
$f^{\ast}h=a\,h$
by the uniqueness in Theorem \ref{thm;12.9.24.50}.
\hfill\qed

\vspace{.1in}

We can obtain the number $a$ by comparing
$\det(h)$ and $f^{\ast}\det(h)$.

\section{Toda-like harmonic bundles}
\label{section;12.12.16.3}
\subsection{Toda-like filtered Higgs bundles on $(\proj^1,\{0,\infty\})$}

\subsubsection{Meromorphic Higgs bundles}
\label{subsection;12.10.3.10}

Let $r$ and $m$ be positive integers.
Let $q$ be a variable.
Let $\nbigk(r,m)(q)$ be the matrix 
whose $(i,j)$-entries $\nbigk(r,m)_{i,j}(q)$ are 
$1$ if $i=j+1$,
$q^m$ if $(i,j)=(1,r)$,
or $0$ otherwise.
We may regard $\nbigk(r,m)$
as a matrix-valued holomorphic function on 
$\cnum_q=\{q\in\cnum\}$,
which is meromorphic at $q=\infty$.

We set $D:=\{0,\infty\}\subset \proj^1$.
Let $\nbigv_r$ be a free
$\nbigo_{\proj^1}(\ast D)$-module
$\bigoplus_{i=1}^r\nbigo_{\proj^1}(\ast D)e_i$.
Let $\theta_{r,m}$ be the Higgs field
determined by
$\theta_{r,m}\vece=\vece\,\nbigk(r,m)\,m\,dq/q$,
where $\vece=(e_1,\ldots,e_r)$.

\begin{rem}
Let $\nbigj$ be the $r$-square matrix
whose $(i,j)$-entries are
$\alpha_{i}q^{\ell_i}$ if $i-j\equiv 1$
modulo $r$,
or $0$ otherwise.
Let $\theta$ be the Higgs field of $\nbigv_r$
determined by
$\theta\vece=\vece\,\nbigj dq/q$.
If $m=\sum \ell_i>0$,
$(\nbigv_r,\theta)$
is transformed to
$(\nbigv_r,\theta_{r,m})$
by an appropriate gauge transform
and an appropriate coordinate change of $\proj^1$.
\end{rem}

Let $\tau=e^{2\pi\sqrt{-1}/r}$.
Let $\sigma$ be the endomorphism
of $\nbigv_r$ determined by
$\sigma(e_i)=\tau^{i}e_i$.
We have
$\sigma^{\ast}\theta:=
 \sigma\circ\theta_{r,m}\circ\sigma^{-1}=\tau\,\theta_{r,m}$.
Hence, $\sigma$ induces 
an isomorphism of meromorphic Higgs bundles
$\sigma:(\nbigv_r,\theta_{r,m})
 \simeq (\nbigv_{r},\tau\,\theta_{r,m})$.

We have a $\cnum^{\ast}$-action 
$\cnum^{\ast}\times\proj^1\lrarr\proj^1$
given by $(t,q)\longmapsto t^rq$.
The bundle $\nbigv_r$ is $\cnum^{\ast}$-equivariant
by the isomorphism
$\ttilde:t^{\ast}\nbigv_r\simeq\nbigv_r$,
$\ttilde(t^{\ast}e_i)=t^{im}e_i$.
We have $t^{\ast}\theta=t^m\theta$.

Let $\ell$ be a factor of $m$.
We set $\mu_{\ell}:=\{\kappa\in\cnum\,|\,\kappa^{\ell}=1\}$.
We have a natural $\mu_{\ell}$-action on $\proj^1$
given by $(\kappa,q)\longmapsto \kappa q$.
Let $\kappatilde:\kappa^{\ast}\nbigv_r\simeq \nbigv_r$
be the isomorphism given by
$\kappatilde(\kappa^{\ast}e_i)=e_i$.
It induces a $\mu_{\ell}$-action on $(\nbigv_r,\theta)$.
The descent is isomorphic to 
$(\nbigv_r,\theta_{r,m/\ell})$.

Let $\varphi:\proj^1_w\lrarr \proj^1_q$ be given by
$\varphi^{\ast}(q)=w^r$.
For $j=1,\ldots,r$,
we set
\begin{equation}
 \label{eq;12.12.5.20}
 v_j:=\sum_{i=1}^r(\tau^jw^m)^{r-i}\varphi^{\ast}e_i.
\end{equation}
We have the following formulas:
\begin{equation}
 \label{eq;12.9.24.10}
 \varphi^{\ast}(\theta_{r,m})v_j=
 \bigl(
 r\,m\,w^m\tau^j
 \bigr)\,v_j
 \frac{dw}{w}
\end{equation}
We also have the following:
\begin{equation}
 \label{eq;12.24.11}
 \sigma(v_j)
=v_{j-1}
\quad(1\leq j\leq r-1),
\quad
 \sigma(v_r)=v_1.
\end{equation}

We have a natural $\cnum^{\ast}$-action
on $\proj^1_w$ given by $(t,w)\longmapsto tw$.
The morphism $\varphi$ is $\cnum^{\ast}$-equivariant,
and $\varphi^{\ast}\nbigv_{r}$ is naturally
$\cnum^{\ast}$-equivariant.
We have 
$\ttilde(t^{\ast}v_j)=t^{mr}v_j$.

We have a $\mu_m$-action on $\proj^1_w$
given by $(\kappa,w)\longmapsto\kappa w$.
(The morphism $\varphi$ is not necessarily
$\mu_m$-equivariant.)
We have a $\mu_m$-action on 
$\varphi^{\ast}\nbigv_r$
by $\kappatilde(\kappa^{\ast}e_i)=e_i$.
We have 
$\kappatilde(\kappa^{\ast}v_j)=v_j$.
Indeed, we have
$\kappatilde(\kappa^{\ast}v_j)
=\sum (\tau^jw^m)^{r-i}
 \kappatilde(\kappa^{\ast}e_i)
=v_j$.

We have the description of $e_k$
in terms of $v_j$ as follows:
\begin{equation}
 \label{eq;12.9.24.101}
e_k=
\frac{1}{rw^{m(r-k)}}
\sum_{i=1}^r \tau^{ki}v_i
\end{equation}
Indeed, 
$\sum_{j=1}^r \tau^{kj}v_j
=\sum_{i=1}^r\sum_{j=1}^r
 \tau^{j(k-i)}(w^m)^{(r-i)}e_i
=rw^{m(r-k)}e_k$.

\subsubsection{Filtrations around $0$}
\label{subsection;12.9.24.30}

Let $U_0=\proj^1\setminus\{\infty\}$.
We shall consider filtered bundles on $(U_0,0)$
over $\nbigv:=\nbigv_{r,m|U_0}$.
The induced sections $e_{i|U_0}$ are denoted by
$e_i$ for simplicity.
If we are given a tuple of real numbers 
$\veca=(a_1,\ldots,a_r)\in\real^r$,
$\nbigv$ is equipped with
the filtration $\nbigp^{\veca}_{\ast}$
given by
$\nbigp^{\veca}_a\nbigv
=\bigoplus \nbigo([a-a_i]\,0)\,e_i$.
Here,
$[a]:=\max\bigl\{c\in\real\,\big|\,c\in\seisuu,c\leq a\bigr\}$
for any $a\in\real$.
The filtration is compatible with the decomposition
$\nbigv=\bigoplus\nbigo_{U_0}(\ast 0)e_i$.
Conversely, if a filtration $\nbigp_{\ast}$ is
compatible with the decomposition,
it is $\nbigp^{\veca}_{\ast}$
for some $\veca$.

\begin{lem}
\label{lem;12.9.24.1}
Let $\nbigp_{\ast}\nbigv$ be a filtered bundle
over $\nbigv$.
The automorphism $\sigma$ of $\nbigv$
gives an automorphism of $\nbigp_{\ast}\nbigv$,
i.e., $\sigma(\nbigp_a\nbigv)\subset\nbigp_a\nbigv$,
if and only if
$\nbigp_{\ast}\nbigv\simeq \nbigp_{\ast}^{\veca}\nbigv$ 
for some $\veca\in\real^r$.
\end{lem}
\pf
The ``if'' part is clear.
To prove the ``only if'' part,
put $a_i:=
 \min\bigl\{a\in\real\,\big|\,e_i\in\nbigp_a\nbigv\bigr\}$,
and $\veca:=(a_i)$.
We clearly have
$\nbigp^{\veca}_a\nbigv\subset\nbigp_a\nbigv$.
A section $f\in\nbigp_a\nbigv$
has the expression $\sum f_ie_i$.
Because $\sigma^j(f)=\sum \tau^{ij}f_ie_i\in\nbigp_a\nbigv$,
we obtain $f_ie_i\in\nbigp_a\nbigv$.
Hence, we obtain $f_i\in\nbigo([a-a_i])$,
which means $f\in\nbigp_a^{\veca}\nbigv$.
\hfill\qed

\vspace{.1in}

The following lemma is clear by the construction
of $\nbigp^{\veca}_{\ast}$.

\begin{lem}
\label{lem;12.10.6.1}
If the conditions in Lemma {\rm\ref{lem;12.9.24.1}}
hold, 
$\nbigp_{\ast}$ is $\cnum^{\ast}$-equivariant.
\end{lem}

Note that, in general,
the converse of the claim in Lemma \ref{lem;12.10.6.1}
does not hold.
For example, in the case $r=m=2$ and $0=a_1=a_2+1$,
$e_1$ and $q^{-1}e_2$ are $\cnum^{\ast}$-invariant sections.
By using any one dimensional subspace
of $\langle e_1,q^{-1}e_2\rangle$,
we obtain a filtration $\nbigp_{\ast}\nbigv_{2}$
such that
it is homogeneous with respect to the $\cnum^{\ast}$-action,
but it is not isomorphic to $\nbigp^{\veca}_{\ast}\nbigv_2$.

The following lemma is also clear by construction.
\begin{lem}
If the conditions in Lemma {\rm\ref{lem;12.9.24.1}}
hold,
$\nbigp_{\ast}$ is $\mu_m$-equivariant.
\end{lem}

Let $\theta:=\theta_{r,m|U_0}$.
We consider the condition
such that 
$\theta$ is logarithmic 
with respect to the filtered bundle
$\nbigp_{\ast}\nbigv$.
The following lemma is clear.
\begin{lem}
\label{lem;12.9.24.110}
Let $\nbigp_{\ast}\nbigv$ be a filtered bundle over $\nbigv$.
We set
$a_i:=\min\{a\,|\,e_i\in\nbigp_a\nbigv\}$.
If $\theta$ is logarithmic,
then we have
$a_i\geq a_{i+1}$ $(i=1,\ldots,r-1)$
and $a_r\geq a_1-m$.
If $\nbigp_{\ast}\nbigv=\nbigp^{\veca}_{\ast}\nbigv$
 for $\veca=(a_i)$,
 the converse is also true.
\end{lem}

\begin{lem}
\label{lem;12.12.6.10}
Assume that $(\nbigp_{\ast}\nbigv,\theta)$ is regular.
After a canonical change of frame,
we may assume that
$a_r>a_1-m$.
\end{lem}
\pf
Assume that $a_r=a_1-m$.
We have $j$ such that
$a_j>a_{j+1}=a_r=a_1-m$.
We set $e'_i:=q^{-m}e_{j+i}$
for $i\leq r-j$,
and 
$e'_i:=e_{i-r+j}$
for $i>r-j$.
Then, we have
$\theta\vece'=\vece'\nbigk(r,m)$,
and we have $a'_{r}>a'_1-m$.
\hfill\qed

\vspace{.1in}

The case $m=1$ is rather special.

\begin{lem}
\label{lem;12.10.6.3}
Assume that $m=1$,
and that $(\nbigp_{\ast}\nbigv,\theta)$ is regular.
Then, the conditions in Lemma 
{\rm\ref{lem;12.9.24.1}} hold.
In particular, $\nbigp_{\ast}\nbigv$ is $\cnum^{\ast}$-equivariant.
\end{lem}
\pf
Let us consider the case $a_r>a_1-1$.
Each $e_i$ gives a section of
$\nbigp_{a_i}\nbigv$.
The induced elements of
$\Gr^{\nbigp}_{a_i}(\nbigv)$ are denoted by
$u_i$.
For each $a_1-1<c\leq a_1$,
we set
$S(c):=\{i\,|\,a_i=c\}$.
We have only to prove that
$\bigl\{u_i\,\big|\,i\in S(c)\bigr\}$
is a base of $\Gr^{\nbigp}_c(\nbigv)$.

Assume that 
$u_i$ $(i\in S(c))$
are linearly dependent
in $\Gr^{\nbigp}_c(\nbigv)$,
and we shall deduce a contradiction.
Let $i_0:=\min\{i\in S(c)\}$
and $i_1:=\max\{i\in S(c)\}$.
We have
$\Res(\theta)u_i=u_{i+1}$
for $i_0\leq i<i_1$
and $\Res(\theta)u_{i_1}=0$.
By the assumption,
we have a non-trivial
linear relation
$\sum \alpha_j u_j=0$.
We have $i_2$  such that
$\alpha_{i_2}\neq 0$
and $\alpha_j=0$ for $j>i_2$.
We may assume that
$u_i$ $(i=i_0,\ldots,i_2-1)$
are linearly independent.
Then, because
$u_{i_2}=\sum_{j=i_0}^{i_2-1}\beta_j u_j$,
$\Res(\theta)$ preserves the subspace $H$
generated by 
$u_j$ $(i_0\leq j<i_2)$.
But, $\Res(\theta)_{|H}$
must have a non-zero eigenvalue.
On the contrary,
because the eigenvalues of
$\nbigk(r,1)_{|q=0}$ are $0$,
$\Res(\theta)$ is nilpotent on
$\Gr^{\nbigp}_c(\nbigv)$ for any $c\in\real$.
Thus, we arrive at a contradiction.
Hence, we obtain that
$u_j$ $(j\in S(c))$ are linearly independent.
Then, we can prove that
they give a base of $\Gr^{\nbigp}_c$,
by using $\sum_{a_1-1<c\leq a_1} \dim \Gr_c=r$.

Let us consider the case $a_r=a_1-1$.
We use the frame $\vece'$ as in the proof of Lemma 
\ref{lem;12.12.6.10}.
Then, by applying the result in the case $a_r>a_1-1$,
we obtain that $\nbigp_{\ast}\nbigv$
is compatible with the decomposition
$\nbigv=\bigoplus \nbigo_{U_0}(\ast 0)e_i$,
i.e.,
$\nbigp_{\ast}\nbigv\simeq
 \nbigp^{\veca}_{\ast}\nbigv$.
\hfill\qed

\begin{cor}
\label{cor;12.10.6.10}
Let $(\nbigp_{\ast}\nbigv,\theta)$
be a filtered Higgs bundle over
$(\nbigv_r,\theta_{r,m})_{|U_0}$.
The conditions in Lemma {\rm\ref{lem;12.9.24.1}}
hold,
if and only if
$\nbigp_{\ast}\nbigv$ is $\mu_m$-equivariant.
\end{cor}
\pf
The ``only if'' part is clear.
Let us prove the ``if'' part.
We have the descent 
$(\nbigp_{\ast}\nbigv',\theta')$
of $(\nbigp_{\ast}\nbigv,\theta)$
with respect to $\varphi:U_0\lrarr U_0$
given by $\varphi(q)=q^m$,
which is a filtered Higgs bundle over
$(\nbigv_r,\theta_{r,1})_{|U_0}$.
Hence, the condition in Lemma \ref{lem;12.9.24.1}
holds for $(\nbigp_{\ast}\nbigv',\theta')$.
Because $(\nbigp_{\ast}\nbigv,\theta)$
is isomorphic to 
$\varphi^{\ast}(\nbigp_{\ast}\nbigv',\theta')$,
the condition in Lemma \ref{lem;12.9.24.1}
holds also for $(\nbigp_{\ast}\nbigv,\theta)$.
\hfill\qed

\subsubsection{Filtration around $\infty$}
\label{subsection;12.9.24.31}

Let $U_{\infty}:=\proj^1\setminus\{0\}$.
Let $(\nbigp_{\ast}\nbigv,\theta)$ be a good filtered Higgs bundle
on $(U_{\infty},\infty)$ 
over $(\nbigv,\theta):=(\nbigv_{r},\theta_{r,m})_{|U_{\infty}}$.
The morphism $\varphi:\proj^1\lrarr\proj^1$
induces $\varphi:U_{\infty}\lrarr U_{\infty}$.
We have the induced filtered Higgs bundle
$\varphi^{\ast}(\nbigp_{\ast}\nbigv,\theta)$.
Because $(\nbigp_{\ast}\nbigv,\theta)$ is good,
the filtration of $\varphi^{\ast}\nbigv$
is compatible with the decomposition
$\varphi^{\ast}\nbigv=
 \bigoplus_{i=1}^r\nbigo_{U_{\infty}}(\ast \infty)v_i$,
where $v_i$ are as in \S\ref{subsection;12.10.3.10}.
Hence, $\nbigp_{\ast}\nbigv$ is determined by
$b_i:=\min\{b\,|\,v_i\in\nbigp_b\varphi^{\ast}\nbigv\}$
$(i=1,\ldots,r)$.

\begin{lem}
\label{lem;12.9.24.20}
Let $(\nbigp_{\ast}\nbigv,\theta)$ be a good filtered Higgs bundle
on $(U_{\infty},\infty)$ over $(\nbigv,\theta)$.
The following conditions are equivalent:
\begin{itemize}
\item
 The filtered bundle $\nbigp_{\ast}\nbigv$
 is compatible with the decomposition
 $\nbigv=\bigoplus\nbigo_{U_{\infty}}(\ast \infty)e_i$.
\item
 The automorphism $\sigma$ of $\nbigv$
 gives an automorphism of $\nbigp_{\ast}\nbigv$.
\item
 The automorphism of $\sigma$
 of $\varphi^{\ast}\nbigv$
 gives an automorphism
 of $\varphi^{\ast}\nbigp_{\ast}\nbigv$.
\item
 $b_i$ are independent of $i$.
\end{itemize}
\end{lem}
\pf
The equivalence of the first and second conditions
can be proved by the argument in the proof of
Lemma \ref{lem;12.9.24.1}.
The equivalence of the second and third conditions
follows from the construction
of the pull back and the descent.
The third condition imply the fourth condition
due to (\ref{eq;12.24.11}).
The fourth condition implies the third condition,
because $\varphi^{\ast}\nbigp_{\ast}\nbigv$
is compatible with the decomposition
$\nbigv=\bigoplus\nbigo_{U_{\infty}}(\ast \infty)v_i$
by the assumption that
$(\nbigp_{\ast}\nbigv,\theta)$ is good.
\hfill\qed

\vspace{.1in}

The following lemma is easy to see.
\begin{lem}
\label{lem;12.10.6.2}
Let $(\nbigp_{\ast}\nbigv,\theta)$ be a good filtered Higgs bundle
on $(U_{\infty},\infty)$ over $(\nbigv,\theta)$.
If the conditions in Lemma 
{\rm\ref{lem;12.9.24.20}} hold,
the filtration $\nbigp_{\ast}\nbigv$ 
is preserved by the $\cnum^{\ast}$-action.
\end{lem}

\begin{lem}
\label{lem;12.9.24.22}
Assume that $m$ and $r$ are relatively prime.
Let $(\nbigp_{\ast}\nbigv,\theta)$ be a good filtered Higgs bundle
on $(U_{\infty},\infty)$
over $(\nbigv,\theta)$.
Then, the conditions in Lemma {\rm\ref{lem;12.9.24.20}} hold.
\end{lem}
\pf
Because $\Gal(\varphi)$ acts on
$\{v_i\}_{i=1,\ldots,r}$ transitively,
$b_i$ are independent of $i$.
\hfill\qed

\begin{lem}
\label{lem;12.10.6.11}
The conditions in Lemma {\rm\ref{lem;12.9.24.20}} hold,
if and only if
$\nbigp_{\ast}\nbigv$ is $\mu_m$-equivariant
\end{lem}
\pf
The ``if'' part is clear.
Suppose that $\nbigp_{\ast}\nbigv$ is $\mu_m$-equivariant.
By Lemma \ref{lem;12.9.24.22},
the conditions in Lemma \ref{lem;12.9.24.20}
hold for the descent.
Hence, they hold for $\nbigp_{\ast}\nbigv$.
\hfill\qed

\subsubsection{Filtered Higgs bundle over $(\nbigv_r,\theta_{r,m})$}

Let $(\nbigp_{\ast}\nbigv_{r},\theta_{r,m})$
be a good filtered Higgs bundle on $(\proj^1,D)$
over $(\nbigv_r,\theta_{r,m})$.
Note that it is logarithmic at $0$,
according to Lemma \ref{lem;12.9.24.40}.
We have the following.
(See 
Lemmas \ref{lem;12.9.24.1},
\ref{lem;12.9.24.20},
\ref{lem;12.10.6.1}
and \ref{lem;12.10.6.2}.)
\begin{prop}
\label{prop;12.9.24.41}
 The automorphism $\sigma$ of $\nbigv_r$
 gives an automorphism of 
 $\nbigp_{\ast}\nbigv_r$,
if and only if 
 $\nbigp_{\ast}\nbigv_r$
 is compatible with the decomposition
 $\nbigv_r=\bigoplus \nbigo_{\proj^1}(\ast D)e_i$.
In that case, 
$\nbigp_{\ast}\nbigv_r$ is $\cnum^{\ast}$-equivariant.
\end{prop}

\begin{prop}
\label{prop;12.9.24.140}
If $m=1$, 
the conditions in
Proposition {\rm\ref{prop;12.9.24.41}} always hold.
If $m\neq 1$,
$(\nbigp_{\ast}\nbigv_r,\theta_{r,m})$ 
satisfies the conditions in 
Proposition {\rm\ref{prop;12.9.24.41}}
if and only if
it is $\mu_m$-equivariant. 
\end{prop}
\pf
The first claim follows from
Lemmas \ref{lem;12.10.6.3}
and \ref{lem;12.9.24.20}.
Note that the equivariance with respect to the
Galois covering $w\longmapsto w^r$
ensures the $4$-th condition in Lemma \ref{lem;12.9.24.20}.
The second claim follows from
Corollary \ref{cor;12.10.6.10}
and Lemma \ref{lem;12.10.6.11}.
\hfill\qed

\vspace{.1in}

A Toda-like filtered Higgs bundle is
a good filtered Higgs bundles
$(\nbigp_{\ast}\nbigv_r,\theta_{r,m})$
such that 
$\nbigp_{\ast}\nbigv_r$ is $\sigma$-invariant.
Let $\nbigf(r,m)$ denote the set of 
such Toda-like filtered Higgs bundles
over $(\nbigv_r,\theta_{r,m})$.
Let $\varphi_m:\proj^1\lrarr\proj^1$
be given by $\varphi_m(q)=q^m$.
We have the pull back
$\varphi_m^{\ast}:\nbigf(r,1)\lrarr\nbigf(r,m)$.
We obtain the following
from Proposition \ref{prop;12.9.24.140}.

\begin{cor}
The map
$\varphi_m^{\ast}:\nbigf(r,1)\lrarr\nbigf(r,m)$
is a bijection.
\end{cor}

If $m$ is relatively prime to $r$,
the Higgs bundle $(\nbigv_r,\theta_{r,m})$
is irreducible,
and hence
any good filtered Higgs bundle is stable.
In general, we have the following.

\begin{prop}
\mbox{{}}\label{prop;12.12.5.1}
Any $(\nbigp_{\ast}\nbigv_r,\theta_{r,m})\in\nbigf(r,m)$ 
is poly-stable.
(We shall give the decomposition into
the stable components in {\rm\S\ref{subsection;12.12.4.31}}.)
\end{prop}
\pf
It is easily reduced to the case 
$\pardeg(\nbigp_{\ast}\nbigv_r)=0$.
Let $(\nbigp_{\ast}\nbigv_{r},\theta_{r,m})
 \in \nbigf(r,m)$.
Because it is $\mu_m$-equivariant,
it is isomorphic to the pull back of
some object
$(\nbigp_{\ast}\nbigv_r,\theta_{r,1})$
in $\nbigf(r,1)$.
Then, the claim follows from
the Kobayashi-Hitchin correspondence.
(See Theorem \ref{thm;12.9.24.50}.)
\hfill\qed

\vspace{.1in}

Let $(\nbigp_{\ast}\nbigv_r,\theta)\in\nbigf(r,m)$.
The filtration around $0$
(resp. $\infty$) is denoted by
$\nbigp^{(0)}_{\ast}$
(resp. $\nbigp^{(\infty)}_{\ast}$).
Let $\varphi:\proj^1\lrarr\proj^1$ be given by
$\varphi(q)=q^r$.
Recall that the parabolic structure at $\infty$
is determined by 
$b=\min\{c\,|\,v_i\in\nbigp^{(\infty)}_c\varphi^{\ast}\nbigv\}$
$(\forall i)$.
It is easy to observe that,
if the condition $\pardeg(\nbigp_{\ast}\nbigv)=0$ is imposed,
$b$ is uniquely determined by
$a_i:=\min\{c\,|\,e_i\in\nbigp^{(0)}_c\nbigv\}$
$(i=1,\ldots,r)$.
Indeed, 
by (\ref{eq;12.9.24.101}),
we have the following
for each $k$:
\begin{equation}
\label{eq;12.12.23.10}
 \min\bigl\{
 c\,\big|\,
 e_k\in\nbigp^{(\infty)}_c\nbigv
 \bigr\}
=\frac{b}{r}-\frac{m(r-k)}{r}
\end{equation}
The condition
$\pardeg(\nbigp_{\ast}\nbigv)=0$
is equivalent to
\[
 \sum_{j=1}^r a_j
+\sum_{k=1}^r
 \left(
 \frac{b}{r}-\frac{m(r-k)}{r}
 \right)
=\sum_{j=1}^ra_j
+b-\frac{m(r-1)}{2}
=0.
\]
Hence, we have
$\pardeg(\nbigp_{\ast}\nbigv)=0$
if and only if
$b=-\sum a_j+m(r-1)/2$.
We have
\begin{equation}
 \label{eq;12.12.23.1}
b-m(r-k)=-\sum a_j-\frac{m(r+1)}{2}+mk.
\end{equation}

Let $\nbigf_0(r,m)\subset\nbigf(r,m)$ 
be determined by the condition that
the parabolic degree is $0$.
We obtain the following proposition
by the above consideration.
\begin{prop}
\label{prop;12.12.21.4}
$\nbigf_0(r,m)$ naturally bijective to
the following:
\[
 \gbigr_{r,m}:=
 \bigl\{
 (a_1,\ldots,a_r)\in\real^r\,\big|\,
a_1\geq a_2\geq\cdots\geq a_r\geq a_1-m
 \bigr\}
\]
\end{prop}
In the following,
let $\nbigp_{\ast}^{\veca}\in\nbigf_0(r,m)$
denote the filtration corresponding to
$\veca\in\gbigr_{r,m}$.

\subsubsection{Decomposition into the stable components}
\label{subsection;12.12.4.31}

Let $\veca=(a_i)\in\gbigr_{r,m}$.
We set $b_i:=ra_i+mi$.
We define $b_{\ell}$ for any $\ell\in\seisuu$
by using the convention $b_{r+\ell}=b_{\ell}$.
We set $r_1:=r/\gcd(m,r)$.
We define
\[
 \nbigs(\veca,m):=
 \bigl\{
 1\leq j\leq r\,\big|\,
 b_{i+j}=b_i\,\,(\forall i),\,\,
 j\equiv 0 \mod r_1
 \bigr\}. 
\]
\begin{prop}
\label{prop;12.12.4.110}
$(\nbigp^{\veca}_{\ast}\nbigv_r,\theta_{r,m})$
is stable,
if $\nbigs(\veca,m)=\{r\}$.
\end{prop}
\pf
According to Proposition \ref{prop;12.12.5.1},
it is polystable.
Hence, we have only to show that
any endomorphism of $(\nbigp_{\ast}\nbigv_r,\theta_{r,m})$
is a scalar multiplication.
Let $\varphi:\proj_w^1\lrarr\proj_q^1$
be given by $\varphi(w)=w^r$.
We set 
$u_i:=w^{-mi}\varphi^{\ast}e_i$ for $i=1,\ldots,r$.
We define  $u_i$ for any $i\in\seisuu$
by the convention $u_i=u_{r+i}$.
We have the parabolic structure $\nbigp_{\ast}$
on $\varphi^{\ast}\nbigv_r$
obtained as the pull back of
$\nbigp^{\veca}_{\ast}\nbigv_r$,
that is
$\nbigp_{\ast}\varphi^{\ast}\nbigv_r
=\bigoplus_{i=1}^r
 \nbigo_{\proj^1}\bigl(
 [c_0-b_i]\cdot 0
+[c_{\infty}-d]\cdot\infty
 \bigr)\,u_i$
for any $(c_0,c_{\infty})\in\real^2$,
where $d\in\real$ is independent of $i$.

Let $g$ be the endomorphism of
$\varphi^{\ast}\nbigv_r$
determined by
$g(u_i)=u_{i+1}$.
We have
$\varphi^{\ast}\theta_{r,m}=g\cdot (rm w^{m-1}dw)$.
Let $\tau$ be an $n$-th primitive root.
Let $\iota_{\tau}$ denote the map $\proj^1_w\lrarr\proj^1_w$
given by $\iota_{\tau}(w)=\tau w$.
We have
$\iota_{\tau}^{\ast}g=\tau^{-m}g$.

Any endomorphism $F$
of $\varphi^{\ast}(\nbigv_r,\theta_{r,m})$
is expressed as
$F=\sum_{j=0}^{r-1} F_j(w)\,g^j$.
If $F$ preserves the parabolic structure at $\infty$,
then we obtain $F_j(w)\in \cnum[w^{-1}]$.
If it preserves the parabolic structure at $0$,
then each $F_jg^j$ also preserves the parabolic
structure at $0$.

Assume $F_j\neq 0$ for some $1\leq j\leq r-1$,
and we shall derive a contradiction.
Let $\ell=\deg_{z^{-1}}F_{j}\geq 0$.
We obtain that $z^{-\ell}g^j$ preserves 
the parabolic structure at $0$.
Then, we obtain
$b_i\geq b_{i+j}+\ell$ for any $i$,
which implies 
that $\ell=0$
and that $b_{i+j}=b_i$ for any $i$.
Note that 
endomorphisms of
$(\nbigp^{\veca}_{\ast}\nbigv_r,\theta_{r,m})$
bijectively correspond
to endomorphisms $F$ of
$\varphi^{\ast}(\nbigp^{\veca}_{\ast}\nbigv_r,\theta_{r,m})$
such that $\iota_{\tau}^{\ast}F=F$.
Hence, 
we have $mj\equiv 0$ modulo $r$,
which implies $j\equiv 0$ modulo $r_1$,
i.e., $j$ is an element of $\nbigs(\veca,m)$.
It contradicts with $\nbigs(\veca,m)=\{r\}$.
Thus, we obtain $F_j=0$ for any $j\neq 0$.
\hfill\qed

\vspace{.1in}

Let $r_0$ be the minimum of
$\nbigs(\veca,m)$.
If we regard $\nbigs(\veca,m)$
as a subset of $\mu_r$
by $j\longmapsto e^{2\pi\sqrt{-1}j/r}$,
then $\nbigs(\veca,m)$ is a subgroup of $\mu_r$.
Hence, $r_0$ is a divisor of any $i\in\nbigs(\veca,m)$,
in particular, a divisor of $r$.
By construction, we have $r_0m\equiv 0$ modulo $r$.

We put $j_0:=r/r_0\in\seisuu$
and $m_0:=m/j_0=mr_0/r\in\seisuu$.
Let $\tau_0$ be a primitive $j_0$-th root of $1$.
For $i=1,\ldots,r_0$
and $s=0,\ldots,j_0-1$,
we set
$x^{(s)}_i:=\sum_{k=0}^{j_0-1}
 \tau_0^{sk}q^{-km_0}e_{i+kr_0}$.
Let $f\in\End(\nbigv_r)$ be determined by
$f\vece=\vece\,\nbigk(r,m)$.
We have
\[
  f(x_i^{(s)})=x^{(s)}_{i+1}
\,\,\,(i=1,\ldots,r_0-1),
\quad\quad
 f(x_{r_0}^{(s)})
=q^{m_0}\tau_0^{-s}x_1^{(s)}.
\]
For $s=0,\ldots,j_0-1$,
we define
$\nbigv^{(s)}:=
 \bigoplus_{j=1}^{r_0}
 \nbigo_{\proj_q^1}(\ast D)x_j^{(s)}$.
We obtain a decomposition
$\nbigv_r=\bigoplus_{s=0}^{j_0-1}\nbigv^{(s)}$
preserved by $\theta_{r,m}$.
Let $\theta^{(s)}:=\theta_{r,m|\nbigv^{(s)}}$.
We have
$\theta^{(s)}(x_i^{(s)})=x_{i+1}^{(s)}\,mdq/q$
for $i=1,\ldots,r_0-1$, 
and 
$\theta^{(s)}(x_{r_0}^{(s)})
=q^{m_0}\tau_0^{-s}\,x_{1}^{(s)}\,mdq/q$.
It is easy to see that
there is no common eigenvalues
of $\theta^{(s)}$ and $\theta^{(s')}$
if $s\neq s'$.
In particular,
any morphism
$(\nbigv^{(s)},\theta^{(s)})
\lrarr
 (\nbigv^{(s')},\theta^{(s')})$
is $0$.

We set $\veca_0=(a_1,\ldots,a_{r_0})\in\real^{r_0}$.
Because $a_{1+r_0}=a_1-m_0$,
we have $a_1\geq\cdots \geq a_{r_0}\geq a_1-m_0$.
We have the parabolic structure 
$\nbigp_{\ast}^{\veca_0}$ of $\nbigv^{(s)}$
for the frame $\vecx^{(s)}=(x_1^{(s)},\ldots,x^{(s)}_{r_0})$.
By construction,
we have the following lemma.
\begin{lem}
\label{lem;12.12.5.10}
The parabolic structure $\nbigp_{\ast}^{\veca}\nbigv_r$
is compatible with the decomposition
$\nbigv_r=\bigoplus \nbigv^{(s)}$.
The induced parabolic structure of $\nbigv^{(s)}$
is the same as
$\nbigp_{\ast}^{\veca_0}\nbigv^{(s)}$.
\end{lem}

For each $s$,
we can choose a complex number $\beta_s$
such that
$\psi_s^{\ast}(\nbigp^{\veca_0}_{\ast}\nbigv^{(s)},\theta^{(s)})
\simeq
 (\nbigp^{\veca_0}_{\ast}\nbigv_{r_0},\theta_{r_0,m_0})$,
where 
$\psi_s:\proj^1\lrarr\proj^1$
is given by $\psi_s(q)=\beta_sq$.

According to Lemma \ref{lem;12.12.5.10},
we have the decomposition
\begin{equation}
\label{eq;12.12.4.120}
(\nbigp_{\ast}^{\veca}\nbigv_r,\theta_{r,m})
=\bigoplus_{s=0}^{j_0-1}
 (\nbigp_{\ast}^{\veca_0}\nbigv^{(s)},\theta^{(s)}).
\end{equation}
By construction,
we have $\nbigs(\veca_0,m_0)=\{r_0\}$.
By Proposition \ref{prop;12.12.4.110}
and the above remark,
each $(\nbigp^{\veca_0}\nbigv^{(s)},\theta^{(s)})$
is stable.
We also have the description of the endomorphisms.

\begin{prop}
\label{prop;12.12.20.110}
Let $f$ be as above, i.e.,
$f\vece=\vece\nbigk(r,m)$.
Then,
the space of endomorphisms of
$(\nbigp^{\veca}_{\ast}\nbigv_r,\theta_{r,m})$
is 
$\Bigl\{
 \sum_{s=0}^{j_0-1}
 \beta_s\,(q^{-m_0}f^{r_0})^{s}\,\Big|\,
 \beta_s\in\cnum
 \Bigr\}$.
\end{prop}
\pf
Let $g$ be as in the proof of 
Proposition {\rm\ref{prop;12.12.4.110}}.
We have
$g=w^{-m}\varphi^{\ast}f$,
and $g^{r_0}=\varphi^{\ast}(q^{-m_0}f^{r_0})$.
Then, the claim is clear from the proof of Proposition
\ref{prop;12.12.4.110}.
\hfill\qed

\subsubsection{Symmetric pairing}
\label{subsection;12.12.21.101}

Let $\nbigv_r^{\lor}$ denote 
the dual of $\nbigv_r$,
i.e.,
$\nhom_{\nbigo_{\proj^1}(\ast D)}
 \bigl(\nbigv_r,\nbigo_{\proj^1}(\ast D)\bigr)$.
It is naturally equipped with the dual $\theta_{r,m}^{\lor}$
of $\theta_{r,m}$.
It is also equipped with the dual frame $e_i^{\lor}$
$(i=1,\ldots,r)$.
Let $\Psi:\nbigv_r\simeq\nbigv_r^{\lor}$
be given by
$\Psi(e_i)=e_{r+1-i}^{\lor}$.
It gives an isomorphism of meromorphic Higgs bundles
$\Psi:(\nbigv_r,\theta_{r,m})\simeq
 (\nbigv_r^{\lor},\theta_{r,m}^{\lor})$.
Let $f\in\End(\nbigv_r)$ be determined by
$\theta_{r,m}=f\,(mdq/q)$.
Let $\veca\in\gbigr_{r,m}$.
Let $\nbigh(\nbigp^{\veca}_{\ast}\nbigv_r,\theta_{r,m})$
denote the space of
morphisms
$(\nbigp^{\veca}_{\ast}\nbigv_r,\theta_{r,m})
\lrarr \bigl((\nbigp^{\veca}_{\ast}\nbigv_r)^{\lor},
 \theta_{r,m}^{\lor}\bigr)$.

\begin{prop}
\label{prop;12.12.21.10}
We have $\nbigh(\nbigp^{\veca}_{\ast}\nbigv_r,\theta_{r,m})
\neq 0$
if and only if
there exist $k_0,\ell\in\seisuu$
with $0\leq k_0\leq r-1$
such that the following holds:
\begin{equation}
 \label{eq;12.12.21.1}
 a_i+a_j=
 \left\{
 \begin{array}{ll}
 -\ell & (i+j=r+1-k_0)\\
 -\ell-m &(i+j=2r+1-k_0)
 \end{array}
 \right.
\end{equation}
In that case,
$\nbigh(\nbigp^{\veca}_{\ast}\nbigv_r,\theta_{r,m})
=\Bigl\{
 \sum_{s=0}^{j_0-1}
 \beta_sq^{\ell-sm_0}\Psi\circ f^{k_0+sr_0}
 \,\Big|\,
 \beta_s\in\cnum
 \Bigr\}$.
\end{prop}
\pf
We use the symbols $\nbigv$ and $\theta$
instead of $\nbigv_r$ and $\theta_{r,m}$,
respectively.
We argue the ``only if'' part of the first claim.
We put $a^{\lor}_i:=-a_{r+1-i}$.
They give $\veca^{\lor}\in\gbigr_{r,m}$,
and $\Psi$ gives an isomorphism
$(\nbigp^{\veca^{\lor}}_{\ast}\nbigv,\theta)
\simeq
 \bigl((\nbigp^{\veca}_{\ast}\nbigv)^{\lor},\theta^{\lor}\bigr)$.
Any morphism
$F:(\nbigv,\theta)
\lrarr 
 \bigl(\nbigv^{\lor},\theta^{\lor}\bigr)$
is described as 
$F=\sum_{k=0}^{r-1}\alpha_k(q)\Psi\circ f^k$.
It is easy to see that
$F$ gives a morphism
$\nbigp^{\veca}_{\ast}\nbigv
\lrarr
 (\nbigp^{\veca}_{\ast}\nbigv)^{\lor}$
if and only if 
each $\alpha_k(q)\Psi\circ f^k$
gives a morphism
$\nbigp^{\veca}_{\ast}\nbigv
\lrarr
(\nbigp^{\veca}_{\ast}\nbigv)^{\lor}$.

Suppose that $(\nbigp^{\veca}_{\ast}\nbigv,\theta)$ is stable.
If $\nbigh(\nbigp^{\veca}_{\ast}\nbigv,\theta)\neq 0$,
by Lemma \ref{lem;12.12.21.2},
there exist $k_0$ with $0\leq k_0\leq r-1$
such that any non-zero element of 
$\nbigh(\nbigp^{\veca}_{\ast}\nbigv,\theta)$
is described as
$\alpha(q)\Psi\circ f^{k_0}$.
If $\alpha(q)\neq 0$, 
then $\alpha(q)\Psi\circ f^{k_0}$ has to be an isomorphism
by Lemma \ref{lem;12.12.21.2} again.
Hence, we obtain that there exists an integer $\ell$
such that
any morphism is described as
$\beta\,q^{\ell}\Psi\circ f^{k_0}$
for some $\beta\in\cnum$,
and moreover
$q^{\ell}\Psi\circ f^{k_0}$
is an isomorphism.
Then, we obtain (\ref{eq;12.12.21.1}),
and thus the ``only if'' part 
in the case that
$(\nbigp^{\veca}_{\ast}\nbigv,\theta)$ is stable.

Let us consider the case that
$(\nbigp^{\veca}_{\ast}\nbigv,\theta)$
is not necessarily stable.
We have the decomposition 
(\ref{eq;12.12.4.120})
into the stable components.
Any morphism 
$(\nbigp^{\veca}_{\ast}\nbigv,\theta)
\lrarr
 \bigl((\nbigp^{\veca}_{\ast}\nbigv)^{\lor},\theta^{\lor}\bigr)$
is the direct sum of
morphisms
$(\nbigp^{\veca_0}_{\ast}\nbigv^{(s)},\theta^{(s)})
\lrarr
 \bigl(
 (\nbigp^{\veca_0}_{\ast}\nbigv^{(s)})^{\lor},(\theta^{(s)})^{\lor}
 \bigr)$.
By using the ``only if'' part in the stable case,
we obtain that 
there exist integers $k_0$ and $\ell_0$
with $0\leq k_0\leq r_0-1$
such that the following holds
for $1\leq i,j\leq r_0$:
\begin{equation}
 \label{eq;12.12.21.3}
 a_i+a_j=\left\{
 \begin{array}{ll}
 -\ell_0 & (i+j=r_0+1-k_0)\\
 -\ell_0-m_0 & (i+j=2r_0+1-k_0)
 \end{array}
 \right.
\end{equation}
If (\ref{eq;12.12.21.3}) holds, (\ref{eq;12.12.21.1}) holds
for $k_0$ and $\ell=\ell_0+r-m_0$.
Thus, the ``only if'' part is proved.

Conversely,
if (\ref{eq;12.12.21.1}) holds
for $k_0$ and $\ell$,
we can easily check that
$q^{\ell}\Psi\circ f^{k_0}$
gives a morphism
$(\nbigp^{\veca}_{\ast}\nbigv,\theta)
\lrarr
 \bigl(
 (\nbigp^{\veca}_{\ast}\nbigv)^{\lor},\theta^{\lor}
 \bigr)$.
Hence, we obtain the ``if'' part.
Note that
$q^{\ell}\Psi\circ f^{k_0}$
is an isomorphism.
The second claim 
follows from Proposition \ref{prop;12.12.20.110}.
\hfill\qed

\vspace{.1in}

The condition in Proposition {\rm\ref{prop;12.12.21.10}}
is equivalent to 
that the following holds
for any $s\in\seisuu$ satisfying $0\leq k_0+sr_0\leq r-1$:
\begin{equation}
\label{eq;12.12.21.100}
 a_i+a_j=\left\{
 \begin{array}{ll}
 -\ell+sm_0 & \bigl(i+j=r+1-(k_0+sr_0)\bigr)\\
 -\ell-m+sm_0 & \bigl(i+j=2r+1-(k_0+sr_0)\bigr)
 \end{array}
 \right.
\end{equation}
It is also equivalent to that
the following holds 
if $i+j-1+k_0\equiv 0$ modulo $r_0$:
\begin{equation}
 \label{eq;12.12.22.100}
  a_i+a_j=-\ell+m-\frac{m_0}{r_0}(i+j-1+k_0)
\end{equation}

\begin{cor}
\label{cor;12.12.22.200}
We have
$\nbigh(\nbigp^{\veca}_{\ast}\nbigv_r,\theta_{r,m})\neq 0$,
if and only if there exists a half-integer $\nu$
such that
$\prod_{i=1}^{r}(T-e^{2\pi\sqrt{-1}(a_i+\nu)/m})
\in\real[T]$.
\end{cor}
\pf
If $\nbigh(\nbigp^{\veca}_{\ast}\nbigv_r,\theta_{r,m})\neq 0$,
we obtain
$\prod_{i=1}^{r}(T-e^{2\pi\sqrt{-1}(a_i+\nu)/m})
\in\real[T]$
for $\nu=\ell/2$
from (\ref{eq;12.12.21.1}) some $k_0$ and $\ell$.
Conversely, suppose 
$\prod_{i=1}^{r}(T-e^{2\pi\sqrt{-1}(a_i+\nu)/m})
\in\real[T]$ for a half-integer $\nu$.
We put $c_i:=a_i+\nu$.
We may assume that
$0\geq c_1\geq\cdots c_r>-m$
by a successive use of the gauge transform 
and the inverse transform 
in the proof of Lemma \ref{lem;12.12.6.10}.
There exists $s$ such that
$0=c_1=\cdots =c_s>c_{s+1}$.
We have
$c_{s+j}+c_{r+1-j}=-m$ for $j=1,\ldots,r-s$.
By Proposition \ref{prop;12.12.21.10},
we have
$\nbigh(\nbigp^{\veca}_{\ast}\nbigv_r,\theta_{r,m})\neq 0$.
\hfill\qed

\vspace{.1in}

For $\nbigo_{\proj^1}(\ast D)$-modules $\nbigf_i$ $(i=1,2)$,
a pairing $\nbigf_1\times\nbigf_2\lrarr\nbigo_{\proj^1}(\ast D)$
means a bi-$\nbigo_{\proj^1}(\ast D)$-homomorphism.
The natural pairing
$\nbigv_r\times\nbigv_r^{\lor}\lrarr\nbigo_{\proj^1}(\ast D)$
is denoted by
$\langle\cdot,\cdot\rangle$.
A pairing $C:\nbigv_r\times\nbigv_r\lrarr\nbigo_{\proj^1}(\ast D)$ 
induces a morphism
$\Phi_C:\nbigv_r\lrarr\nbigv_r^{\lor}$
by $\langle\Phi_C(u),v\rangle=C(u,v)$
for local sections $u,v$ of $\nbigv_r$.
We have
$C(\theta_{r,m} u,v)=C(u,\theta_{r,m} v)$
if and only if
$\Phi_C\circ\theta_{r,m}=\theta_{r,m}^{\lor}\circ\Phi_C$.

Let $C_{k}$ denote the pairing of 
$\nbigv_r\times\nbigv_r\lrarr\nbigo_{\proj^1}(\ast D)$
given by 
$C_k(u,v):=\langle
 \Psi\circ f^k(u),v
 \rangle$.
For $0\leq k\leq r-1$,
we have the following:
\[
 C_k(e_i,e_j)=
 \left\{
 \begin{array}{ll}
 1 & (i+j=r+1-k)\\
 q^{m} &(i+j=2r+1-k)\\
 0  & \mbox{\rm otherwise}
 \end{array}
 \right.
\]
In particular,
$C_k$ are symmetric,
i.e., $C_k(u,v)=C_k(v,u)$.

Let $\gbigp(\nbigp^{\veca}_{\ast}\nbigv_r,\theta_{r,m})$ 
denote the space of pairings 
$C:\nbigv_r\times\nbigv_r\lrarr\nbigo_{\proj^1}(\ast D)$
such that 
$\Phi_C$ gives a morphism
$(\nbigp^{\veca}_{\ast}\nbigv_r,\theta_{r,m})
\lrarr
 \bigl(
 (\nbigp^{\veca}_{\ast}\nbigv_r)^{\lor},\theta_{r,m}^{\lor}
 \bigr)$.
\begin{cor}
We have 
$\gbigp(\nbigp^{\veca}_{\ast}\nbigv_r,\theta_{r,m})
 \neq 0$,
if and only if 
there exist integers $k_0$ and $\ell$
with $0\leq k_0\leq r-1$
such that {\rm(\ref{eq;12.12.21.1})} holds.
In that case, 
{\rm(\ref{eq;12.12.22.100})} holds
if $i+j-1+k_0\equiv 0$ modulo $r_0$,
and we have
\[
 \gbigp(\nbigp^{\veca}_{\ast}\nbigv_r,\theta_{r,m})
=\Bigl\{
 \sum_{s=0}^{j_0-1} \beta_sq^{\ell-sm_0}C_{k_0+sr_0}
 \,\Big|\,
 \beta_s\in\cnum
 \Bigr\}.
\]
In particular,
any $C\in \gbigp(\nbigp^{\veca}_{\ast}\nbigv_r,\theta_{r,m})$
satisfies $t^{\ast}C=t^{r\ell-m(r+1-k_0)}C$.
\end{cor}
\pf
The first two claims follows from
Proposition \ref{prop;12.12.21.10}.
The last claim can be checked directly.
\hfill\qed

\begin{cor}
Suppose $(\nbigp^{\veca}_{\ast}\nbigv_r,\theta_{r,m})$ is stable.
We have $\gbigp(\nbigp^{\veca}_{\ast}\nbigv_r,\theta_{r,m})\neq 0$,
if and only if
there exist integers $k_0$ and $\ell$
with $0\leq k_0\leq r-1$
such that {\rm(\ref{eq;12.12.21.1})} holds.
In that case, we have
$\gbigp(\nbigp^{\veca}_{\ast}\nbigv_r,\theta_{r,m})
=\bigl\{
 \beta q^{\ell}C_{k_0}
 \,\big|\,
 \beta\in\cnum
 \bigr\}$.
\end{cor}

As remarked in Corollary \ref{cor;12.12.22.200},
the existence of $k_0$ and $\ell$ as above 
can be replaced with the existence of
a half-integer $\nu$
such that 
$\prod_{i=1}^{r}(T-e^{2\pi\sqrt{-1}(a_i+\nu)/m})\in\real[T]$.

\subsection{Toda-like Harmonic bundle}

\subsubsection{Preliminaries}

We set $E_{r}:=\nbigv_{r|\proj^1\setminus D}$,
which is equipped with the base $e_{i|\proj^1\setminus D}$
and the Higgs field $\theta_{r,m|\proj^1\setminus D}$.
We will omit ``$|\proj^1\setminus D$''.
Let $h$ be any harmonic metric of
$(E_{r},\theta_{r,m})$.
Clearly,
$(E_{r},\theta_{r,m},h)$ is tame at $0$,
and wild at $\infty$.
We have the associated good filtered Higgs bundle
denoted by $(\nbigp^{h}_{\ast}E_{r},\theta_{r,m})$
and the meromorphic Higgs bundle
$(\nbigp^hE_r,\theta_{r,m})$.
(We add the superscript $h$ to emphasize the dependence
on $h$.)

\subsubsection{Toda-like harmonic metric}
\label{subsection;13.1.8.1}

In this paper,
a harmonic metric $h$ of $(E_{r},\theta_{r,m})$
is called Toda-like,
if $h(e_i,e_j)=0$ for $i\neq j$.
In that case,
$(E_r,\theta_{r,m},h)$ is called 
a Toda-like harmonic bundle.
Clearly,
the condition is equivalent to 
 $\sigma^{\ast}h=h$.
The following lemma is clear.
\begin{lem}
\label{lem;12.10.17.2}
If $h$ is a Toda-like harmonic metric of
$(E_r,\theta_{r,m})$,
then $\sigma$ gives an automorphism of
$\nbigp^h_{\ast}E$.
\end{lem}

A Toda-like harmonic metric
is uniquely determined
by the parabolic structure
in the sense of the following lemma.

\begin{lem}
\label{lem;12.10.7.1}
Let $h_i$ $(i=1,2)$ be Toda-like harmonic metrics
of $(E_r,\theta_{r,m})$
such that
$\nbigp^{h_1}_{\ast}E=\nbigp^{h_2}_{\ast}E$.
Then, $h_1=a\,h_2$ for some $a>0$.
\end{lem}
\pf
We have an automorphism
$g$ of $(\nbigp^{h_1}_{\ast}E_r,\theta_{r,m})$
such that 
(i) $g$ is self-adjoint with respect to both of $h_i$ $(i=1,2)$,
(ii) $h_1=h_2\,g$.
Because $\sigma^{\ast}h_i=h_i$,
we have $\sigma\circ g=g\circ \sigma$.
We also have $g\circ \theta=\theta\circ g$.
Then, it is easy to prove that 
$g=a\,\id_{E_{r}}$.
\hfill\qed

\subsubsection{Adaptedness}

A harmonic metric $h$ of $(E_r,\theta_{r,m})$
is called adapted to $\nbigv_{r}$,
if $\nbigp^hE_r=\nbigv_r$.
The case $m=1$ is special.

\begin{prop}
\label{prop;12.9.24.150}
Any harmonic metric  of 
$(E_{r},\theta_{r,1})$ adapted to $\nbigv_r$
is \!Toda-like.
\end{prop}
\pf
Because 
$(E_r,\theta_{r,1})$ is irreducible,
the filtered Higgs bundle $(\nbigp^h_{\ast}E_r,\theta_{r,1})$
is stable.
Because $\sigma^{\ast}\theta=\tau\theta$
with $|\tau|=1$,
$\sigma^{\ast}h$ is also a harmonic metric of
$(E_r,\theta_{r,1})$.
By Proposition \ref{prop;12.9.24.140},
we have
$\nbigp^{\sigma^{\ast}h}_{\ast}E_r
\simeq
 \sigma^{\ast}\nbigp^{h}_{\ast}E_r
\simeq
 \nbigp^h_{\ast}E_r$.
Hence, we have $\sigma^{\ast}h=a\,h$ for some $a>0$
by the uniqueness in Theorem \ref{thm;12.9.24.50}.
Because we clearly have
$\det(\sigma^{\ast}h)=\det(h)$,
we obtain $\sigma^{\ast}h=h$.
\hfill\qed

\vspace{.1in}

We give a complement.
We set $\Omega:=e_1\wedge\cdots\wedge e_r$.

\begin{lem}
\label{lem;12.9.24.130}
A harmonic metric $h$ of $(E_r,\theta_{r,m})$
is adapted to $\nbigv_r$, if and only if
$|e_i|_h=O\bigl((|q|+|q|^{-1})^{N}\bigr)$ 
for some $i$ and some $N$.
Suppose moreover that $h$ is Toda-like.
Then, $h$ is adapted to $\nbigv_r$,
if and only if
$\bigl|\Omega\bigr|_h
=O\bigl((|q|+|q|^{-1})^{N}\bigr)$
for some $N$.
\end{lem}
\pf
If $h$ is adapted, 
any $e_i$ are sections of $\nbigp^hE_r$,
and hence 
we have
$|e_i|_h=O\bigl((|q|+|q|^{-1})^{N}\bigr)$
for some $N$.
Conversely,
assume 
$|e_1|_h=O\bigl((|q|+|q|^{-1})^{N}\bigr)$.
For the expression $\theta=f\cdot \bigl(m\,dq/q\bigr)$,
$f$ naturally gives an endomorphism of 
$\nbigp^hE_r$,
according to 
\cite{biquard-boalch}, \cite{Mochizuki-wild}
and \cite{Simpson90}.
Because $e_j=f^{j-1}(e_1)$ for $j=1,\ldots,r$,
we obtain that $e_j$ are also sections of
$\nbigp^hE_r$.
Then, we obtain
$\nbigv_r\subset\nbigp^hE_r$,
which implies
$\nbigv_r=\nbigp^hE_r$
because both of them are $\nbigo_{\proj^1}(\ast D)$-locally free.
Because $f(e_i)=e_{i+1}$ $(i\leq r-1)$ and $f(e_r)=q^me_1$,
we may replace the role of $e_1$
with the other $e_i$.
Hence, we obtain the first claim.

Let us prove the second claim.
The ``only if'' part is clear.
Assume that $h$ is Toda-like,
and that 
$\bigl|\Omega\bigr|_h
=O\bigl((|q|+|q|^{-1})^{N}\bigr)$
for some $N$.
We have
$\bigl|
 \Omega
 \bigr|_h=
C\prod_{i=1}^r|e_i|_h$
for some $C>0$.
By the relation $e_r=f^{r-j}(e_j)$,
we have $|e_j|_h\geq (|q|+|q|^{-1})^{-N}|e_r|_h$.
Hence,
$|e_r|_h^r\leq (|q|+|q|^{-1})^{M}$ for some $M$.
We obtain that $h$ is adapted.
\hfill\qed

\vspace{.1in}

It is easy to see that, if $h$ is adapted,
we have
$\bigl|\Omega\bigr|_h
=C|q|^{a}$
for some $C>0$ and $a\in\real$.
Indeed, 
$\Omega$ is a section of
$\nbigp\det(E)$,
and we have 
$\del_q\del_{\qbar}\log|\Omega|_h=0$.
Because the parabolic degree of
$\nbigp_{\ast}\det(E)$ is $0$,
there exists $a\in\real$ such that
$\log|\Omega|_h-a\log|q|$
is bounded from above on $\proj^1\setminus D$.
Then, we obtain that 
$\log|\Omega|_h-a\log|q|$ is constant.

\begin{cor}
\label{cor;12.9.24.131}
Let $(E_{r},\theta_{r,m},h)$
be any Toda-like harmonic bundle.
There exists a holomorphic function $g$ on $\cnum^{\ast}$
such that
the Toda-like harmonic metric $e^{\Re g}h$ 
of $(E_r,\theta_{r,m})$
is adapted to $\nbigv_r$.
\end{cor}
\pf
We have a holomorphic function $g_1$ on $\cnum^{\ast}$
such that
(i) $g_1$ is nowhere vanishing on $\cnum^{\ast}$,
(ii) $g_1\,\Omega$
is a section of $\det\nbigp^h E$.
It means
$\det h(\Omega,\Omega)|g_1|^2=
 (|q|+|q|^{-1})^N$ for some $N$.
We have an integer $\ell$
such that
$\int_{|q|=1}(dg_1/g_1+\ell dq/q)=0$.
We can find 
a holomorphic function $g_2$ on $\cnum^{\ast}$
such that $e^{g_2}=g_1\,q^{\ell}$.
We obtain that 
$e^{2\Re(g_2)}\det h(\Omega,\Omega)=
O\bigl((|q|+|q|^{-1})^N\bigr)$ for some $N$.
Hence, 
$e^{2\Re(g_2)/r}h$ is adapted to $\nbigv_r$,
according to Lemma \ref{lem;12.9.24.130}.
\hfill\qed

\subsubsection{Homogeneity and $\mu_m$-equivariance}

We have the symmetry of Toda-like harmonic bundles,
inherited to the symmetry of the underlying filtered Higgs bundles.
By the restriction of the $\cnum^{\ast}$-action $\rho$,
we obtain the $S^1$-action on $E_r$,
i.e.,
$\rho:S^1\times \cnum^{\ast}\lrarr \cnum^{\ast}$
is given by $(t,q)\longmapsto t^rq$,
and $\ttilde:t^{\ast}E_r\simeq E_r$ is given by
$\ttilde(t^{\ast}e_i)=t^{im}e_i$.
We have $t^{\ast}\theta=t^m\theta$.

\begin{lem}
\label{lem;12.10.7.3}
Let $h$ be a Toda-like harmonic metric of
$(E_{r},\theta_{r,m})$
adapted to $\nbigv_r$.
Then, $t^{\ast}h=h$
for $t\in S^1$.
\end{lem}
\pf
According to Lemma \ref{lem;12.10.17.2},
the conditions in Proposition \ref{prop;12.9.24.41}
hold.
In particular, $\nbigp_{\ast}E_r$ 
is $\cnum^{\ast}$-equivariant.
We clearly have
$\det(t^{\ast}h)=\det(h)$
for any $t\in S^1$.
Because $t^{\ast}h$ is also Toda-like,
the claim follows from Lemma \ref{lem;12.10.7.1}.
\hfill\qed

\vspace{.1in}

We have the natural action of $\mu_m$
on $(E_{r},\theta_{r,m})$
given by
$\kappatilde(\kappa^{\ast}e_i)=e_i$
for $\kappa\in\mu_m$.
The following lemma follows from
Lemma \ref{lem;12.10.7.3},
or it can be also proved easily
by the argument in the proof of
Lemma \ref{lem;12.10.7.3}.

\begin{lem}
\label{lem;12.10.17.10}
Any Toda-like harmonic metric of
$(E_r,\theta_{r,m})$ 
is $\mu_m$-equivariant,
if it is adapted to $\nbigv_r$.
\end{lem}

\subsubsection{Correspondence}
\label{subsection;12.10.18.10}

Let $a\in\real$.
Let $\gbigh_{r,m}(a)$ be the set of Toda-like harmonic metrics
of $(E_r,\theta_{r,m})$ such that $|\Omega|_h=|q|^{-a}$.
Let $\gbigr_{r,m}(a)$ be the set of
$\veca=(a_i)\in\real^{r}$
 such that
 $a_1\geq a_2\geq\cdots \geq a_r\geq a_1-m$
 and $\sum a_i=a$.
From $h\in\gbigh_{r,m}(a)$,
we obtain 
$\Phi_{r,m}(h)=\bigl(a_i(h)\bigr)\in\real^r$
such that
$a_i(h):=\min\big\{
 c\,\big|\,e_i\in\nbigp^h_cE \,\,\mbox{\rm at }0
 \}$.
They satisfy
$a_i(h)\geq a_{i+1}(h)$,
$a_r(h)\geq a_{1}(h)-m$
and $\sum a_i(h)=a$.
Thus, we obtain a map
$\Phi_{r,m}:\gbigh_{r,m}(a)\lrarr\gbigr_{r,m}(a)$.

\begin{thm}
\label{thm;12.9.24.200}
$\Phi_{r,m}$ is bijective.
\end{thm}
\pf
It is injective by Lemma \ref{lem;12.10.7.1}.
Let us consider the case $m=1$.
Let $\veca\in\gbigr_{r,1}(a)$.
Let $(\nbigp_{\ast}^{\veca}\nbigv_{r},\theta_{r,1})$ be 
the associated filtered Higgs bundle.
It is always stable.
By the Kobayashi-Hitchin correspondence,
we have a unique harmonic metric $h$
of $(E_{r},\theta_{r,1})$
such that
$\nbigp^h_{\ast}E_r\simeq
 \nbigp^{\veca}_{\ast}\nbigv_r$
and $|\Omega|_h=|q|^{-a}$.
It is automatically Toda-like
as remarked in Proposition  \ref{prop;12.9.24.150}.
It means that $\Phi_{r,1}$ is a bijection.

Let us consider the general case.
Let $\veca\in\gbigr_{r,m}(a)$.
Let $(\nbigp_{\ast}^{\veca}\nbigv,\theta_{r,m})$ be 
the associated filtered Higgs bundle.
It is $\mu_m$-equivariant.
Any Toda-like harmonic metrics are 
also $\mu_m$-equivariant
as remarked in Lemma \ref{lem;12.10.17.10}.
Then, the bijectivity of $\Phi_{r,m}$
follows from the bijectivity of $\Phi_{r,1}$.
Thus, we obtain Theorem \ref{thm;12.9.24.200}.
\hfill\qed

\paragraph{Asymptotic behaviour around $0$ and $\infty$}

We have the norm estimate for wild harmonic bundles
(\cite{biquard-boalch},
\cite{mochi2}, \cite{Mochizuki-wild},
\cite{Simpson90}).
In our situation,
it is described as follows.
Suppose that $h\in\gbigh_{r,m}(a)$
corresponds to
$\veca\in \gbigr_{r,m}(a)$.
We obtain the following around $\infty$
from 
 (\ref{eq;12.12.23.10}) and (\ref{eq;12.12.23.1}):
\begin{equation}
 \label{eq;12.12.23.11}
 \log|e_i|_h
+\frac{1}{r}
 \left(
 a+\frac{m(r+1)}{2}-mi 
 \right)
 \times\log|q|
=O(1)
\end{equation}
Let us describe the asymptotic behaviour
of $\log|e_i|_h$ around $0$.
For $c\in \real$,
we set $S(\veca,c)=\bigl\{i\,|\,a_i=c\bigr\}$.
First, let us consider the case $a_r>a_1-m$.
For $i$,
we define
\begin{equation}
 \label{eq;12.12.22.101}
 k_i:=
 \bigl|
 S(\veca,a_i)
 \bigr|
-2\bigl|\{j\,|\,a_j=a_i,\,j\leq i\}\bigr|
+1.
\end{equation}
Then, we have the following around $0$:
\begin{equation}
 \label{eq;12.9.24.200}
 \log|e_i|_h
+a_i\log|q|
-\frac{k_i}{2}\log\bigl(-\log|q|\bigr)
=O(1)
\end{equation}
Let us consider the case $a_r=a_1-m$.
If $a_i\neq a_j$ $(j=1,r)$,
then the behaviour of $\log|e_i|_h$ around $0$
is given by (\ref{eq;12.9.24.200}).
We have $j_0$ such that 
$a_{j_0}>a_{j_0+1}=a_r$.
We have $j_1$ such that
$a_1=a_{j_1}>a_{j_1+1}$.
We put
\[
 S'(\veca,a_1):=
S(\veca,a_1-1)\cup S(\veca,a_1)
=\bigl\{j_0+1,j_0+2,\ldots,r,1,\ldots,j_1\bigr\}.
\]
Let $\leq'$ be the order given by
$j_0+1,j_0+2,\ldots,r,1,\ldots,j_1$.
For $i\in S'(\veca,a_1)$, we put
\[
 k_i':=
 \bigl|S'(\veca,a_1)\bigr|
-2\bigl|
 \{j\in S'(\veca,a_1)\,|\, j\leq'i \}
 \bigr|
+1.
\]
Then, for $i\in S'(\veca,a_1)$,
we have the following around $0$:
\begin{equation}
\label{eq;12.12.2.1}
  \log|e_i|_h
+a_i\log|q|
-\frac{k_i'}{2}\log\bigl(-\log|q|\bigr)
=O(1)
\end{equation}

\subsubsection{Additional symmetry}
\label{subsection;12.11.28.2}

We use the notation in 
\S\ref{subsection;12.12.4.31}--\ref{subsection;12.12.21.101}.
Let $\veca\in\gbigr_{r,m}(a)$.
Suppose
$\gbigp(\nbigp^{\veca}_{\ast}\nbigv_r,\theta_{r,m})\neq 0$.
Let $k_0$ and $\ell$ be as in Proposition \ref{prop;12.12.21.10}.
See also the remark right after the proposition.

\begin{prop}
\label{prop;12.10.18.42}
Let $h\in\gbigh_{r,m}(a)$ be the adapted
Toda-like harmonic metric of
$(\nbigp^{\veca}_{\ast}\nbigv_r,\theta_{r,m})$.
If $i+j+k_0-1\equiv 0$ modulo $r_0$,
then we have
\begin{equation}
 \label{eq;12.12.21.130}
\log |e_i|_h+\log |e_j|_h
=\Bigl(
 \ell-m+\frac{m_0}{r_0}(i+j+k_0-1) \Bigr)
 \log|q|.
\end{equation}
\end{prop}
\pf
Let $h^{\lor}$ denote the induced Toda-like harmonic metric
of $(E_r^{\lor},\theta_{r,m}^{\lor})$,
which is adapted to
$(\nbigp^{\veca}_{\ast}\nbigv_r)^{\lor}$.
Recall that we have an isomorphism
$\Phi:=q^{\ell-sm_0}\Psi\circ f^{k_0+sr_0}:
 (\nbigp^{\veca}_{\ast}\nbigv_r,\theta_{r,m})
\lrarr
 \bigl((\nbigp^{\veca}_{\ast}\nbigv_r)^{\lor},
 \theta^{\lor}_{r,m}\bigr)$.
By the uniqueness of the Toda-like harmonic metric,
$\Phi:(E_r,h)\lrarr (E_r^{\lor},h^{\lor})$
is isometric up to the multiplication
of the positive constant.
We have
$\det\Phi(e_1\wedge\cdots\wedge e_r)
=\pm q^{b_1}e_1^{\lor}\wedge\cdots\wedge
 e_r^{\lor}$ for some $b_1\in\seisuu$.
We have
$\bigl|e_1^{\lor}\wedge\cdots\wedge e_r^{\lor}\bigr|_{h^{\lor}}
=|q|^{b_2}$ for some $b_2\in\seisuu$.
Hence, we obtain that
$\det\Phi:\det(E_r,h)\lrarr
\det(E_r,h)^{\lor}$ is an isometry,
and thus $\Phi$ is an isometry.
It implies
\[
 h(e_i,e_i)
=\left\{
 \begin{array}{ll}
 h^{\lor}(e_j^{\lor},e_j^{\lor})|q|^{2(\ell-sm_0)}
 & \bigl(i+j=r+1-(k_0+sr_0)
 \bigr) \\
  h^{\lor}(e_j^{\lor},e_j^{\lor})|q|^{2(\ell+m-sm_0)}
 & \bigl(i+j=2r+1-(k_0+sr_0)
 \bigr).
 \end{array}
 \right.
\]
Thus, we obtain (\ref{eq;12.12.21.130}).
\hfill\qed

\subsection{Toda lattice with opposite sign}
\label{subsection;13.1.8.2}

According to {\rm\cite{Guest-Lin}},
the Toda lattice with opposite sign 
is the following equation
for $\real^r$-valued function
$\vecw=(w_i\,|\,i=1,\ldots,r)$
on $\cnum^{\ast}$:
\[
 2\delbar_z\del_zw_i
-e^{2(w_i-w_{i-1})}
+e^{2(w_{i+1}-w_i)}=0
\]
Here, we use the convention $w_{r+i}:=w_i$ for any $i$.
By the coordinate change $z=rq$,
it is transformed to
\begin{equation}
 \label{eq;12.9.24.201}
 2\delbar_q\del_qw_i
-r^2e^{2(w_i-w_{i-1})}
+r^2e^{2(w_{i+1}-w_i)}
=0
\end{equation}

Recall that
the solutions of the equation (\ref{eq;12.9.24.201})
naturally correspond to
the Toda-like harmonic metrics of
the Higgs bundle $(E_r,\theta_{r,r})$.
Let $h$ be a Toda-like harmonic metric of
$(E_r,\theta_{r,r})$.
We set $u_i:=q^{-(i-1)}e_i$.
Note that $\theta$ is represented by $A\,dq$,
where $A_{ij}=r$ if $i-j\equiv 1$ modulo $r$,
and $A_{ij}=0$ otherwise.
Let $w_i$ be the functions determined by
$e^{2w_i}=h(u_i,u_i)$.
Then, it is direct to check that
$(w_i\,|\,i=1,\ldots,r)$ satisfies
the equation (\ref{eq;12.9.24.201}).
Conversely,
if we are given a solution of $(w_i\,|\,i=1,\ldots,r)$
of (\ref{eq;12.9.24.201}),
we define the metric $h$ by
$h(u_i,u_i):=e^{2w_i}$ and $h(u_i,u_j)=0$ $(i\neq j)$.
Then, $h$ is a Toda-like harmonic metric of
$(E_r,\theta_{r,r})$.

\vspace{.1in}

Let $\gbigs_r(a)$ be the set of solutions 
$\vecw=(w_i)$
of (\ref{eq;12.9.24.201})
such that 
\[
 \sum_{i=1}^r w_i
=\left(
 -\frac{r(r-1)}{2}-a
 \right)\log|q|.
\]

We obtain the following theorem
from Theorem \ref{thm;12.9.24.200}
and Lemma \ref{lem;12.10.7.3}.
\begin{thm}
\label{thm;12.10.6.100}
We have a natural bijection $\gbigs_r(a)$
and $\gbigr_{r,r}(a)$.
For any $\vecw\in \gbigs_r(a)$,
the functions $w_i$ depend only on $|q|$.
\end{thm}

Let $\vecw=(w_i)$ be any solution of the equation
{\rm(\ref{eq;12.9.24.201})}.
Then, $\sum w_i$ is an $\real$-valued
harmonic function on $\cnum^{\ast}$.
Hence, it is described as
$A\log |z|^2+\Re f$ 
for a real number $A$ and
a holomorphic function $f$.
If we put $\wtilde_i:=w_i-\Re(f)/r$,
then $\vecwtilde_i$ is a solution of
{\rm(\ref{eq;12.9.24.201})}
such that $\sum \wtilde_i=A\log|z|^2$.
In this sense,
Theorem {\rm\ref{thm;12.10.6.100}}
gives a classification of the $\real$-valued solutions
of {\rm(\ref{eq;12.9.24.201})}.

\paragraph{Asymptotic behaviour}
Let $\veca\in\gbigr_{r,r}(a)$.
Let $\vecw\in\gbigs_r(a)$ 
and $h\in\gbigh_{r,r}(a)$
corresponding to $\veca$.
Let us describe the asymptotic behaviour of $\vecw$.
Note 
$w_i=\log|u_i|_h=\log|e_i|_h-(i-1)\log|q|$.
We obtain the following around $q=\infty$
from (\ref{eq;12.12.23.11}):
\[
 w_i=\log|q|\times
 \left(
 -\frac{a}{r}-\frac{r-1}{2}
 \right)
+O(1)
\]

Let $k_i$ be as in (\ref{eq;12.12.22.101}).
Around $q=0$,
if $a_r>a_1-r$,
we obtain the following
for any $i=1,\ldots,r$:
\[
 w_i=   
-(a_i+i-1)\log|q|
+\frac{k_i}{2}
 \log\bigl(-\log|q|\bigr)
+O(1)
\]
If $a_r=a_1-r$,
the numbers $k_i$ should be modified to $k_i'$
as in (\ref{eq;12.12.2.1}).

\paragraph{Additional symmetry}

Let $\vecw\in\gbigs_r(a)$ corresponding to
$\veca\in\gbigr_{r,r}$.
Suppose that there exist integers $k_0$ and $\ell$
with $0\leq k_0\leq r-1$
such that the following holds:
\begin{equation}
 \label{eq;12.12.21.15}
 a_i+a_j=
 \left\{
 \begin{array}{ll}
 -\ell & (i+j=r+1-k_0)\\
 -\ell-r &(i+j=2r+1-k_0)
 \end{array}
 \right.
\end{equation}
We use the notation in 
\S\ref{subsection;12.12.4.31}--\ref{subsection;12.12.21.101}.

\begin{thm}
$w_i+w_j=(\ell-r+1+k_0)\log|q|$
if $i+j+k_0-1\equiv 0$ modulo $r_0$.
\end{thm}
\pf
Let $(E_r,\theta_{r,r},h)$ be the corresponding
Toda-like harmonic bundle.
If $i+j+k_0-1\equiv 0$ modulo $r_0$,
we have
$w_i+w_j=
 \log |e_i|_h-(i-1)\log|q|
+\log|e_j|_h-(j-1)\log|q|
=\bigl(
 \ell-m+k_0+1
 \bigr)\log|q|$.
\hfill\qed

\vspace{.1in}

We state some special cases as corollaries.
\begin{cor}
Assume that there exist
an integer $\ell$ satisfying the following:
\begin{equation}
 \label{eq;12.10.19.1}
 a_i+a_j=\left\{
 \begin{array}{ll}
 -\ell & (\mbox{\rm if }i+j=2+\ell)\\
 -\ell-r & (\mbox{\rm if }i+j=2+\ell+r)
 \end{array}
 \right.
\end{equation}
Then, we have
$w_i+w_j=0
\quad
 (i+j\equiv 2+\ell\,\,\modulo r)$.
In particular,
if $a_i+a_j=-r+1$ holds for 
$(i,j)$ with $i+j=r+1$,
we have
$w_i+w_j=0$.
\end{cor}

\part{}

\section{Stokes factor and integral structure}
\label{section;12.10.8.10}
\subsection{Preliminary}

\subsubsection{Stokes structure}
\label{subsection;12.10.8.4}

Let us recall the concept of Stokes structure
in the local and unramified case.
Let $X:=\bigl\{z\in\cnum\,\big|\,|z|<1\bigr\}$
and $D:=\{0\}$.
Let $(V,\nabla)$ be a meromorphic flat bundle
on $(X,D)$.
For simplicity, we assume that it is unramified,
i.e., its formal completion
$(V,\nabla)\otimes_{\nbigo_X}\cnum[\![z]\!]$
has a formal decomposition
\begin{equation}
 \label{eq;12.10.8.1}
 (V,\nabla)\otimes_{\nbigo_X}\cnum[\![z]\!]
=
 \bigoplus_{\gminia\in z^{-1}\cnum[z^{-1}]}
 (\Vhat_{\gminia},\nablahat_{\gminia}),
\end{equation}
where $\nablahat_{\gminia}-d\gminia$ are regular singular.
We set $\nbigi(\nabla):=
 \bigl\{\gminia\in z^{-1}\cnum[z^{-1}]\,\big|\,
 \Vhat_{\gminia}\neq 0
 \bigr\}$.

\begin{rem}
\label{rem;12.10.8.3}
In general, according to the Hukuhara-Levelt-Turrittin theorem,
we have such a decomposition
of $(V,\nabla)\otimes_{\nbigo_X}\cnum[\![z^{1/e}]\!]$.
\hfill\qed
\end{rem}

In this paper,
a sector $S\subset X\setminus D$ means a subset
of the form 
$\bigl\{
 z\in X\setminus D\,\big|\,
 0<|z|<R_0,\,\,
 \theta_0\leq \arg(z)\leq\theta_1
 \bigr\}$.
We set $\del S=\{z\in S\,|\,\arg(z)=\theta_0\,\,\mbox{\rm or }
 \theta_1\}$.
For each $\gminia,\gminib\in\nbigi(V)$
with $\gminia\neq\gminib$,
put $S(\gminia,\gminib):=\bigl\{
 z\in S\,\big|\,
 \Re(\gminia-\gminib)(z)=0
 \bigr\}$.
A sector $S$ is called adapted to $(V,\nabla)$,
if the following holds
 for each $\gminia,\gminib\in \nbigi(V)$ with $\gminia\neq\gminib$:
\begin{itemize}
\item
$S(\gminia,\gminib)$ is connected or empty,
and we have $S(\gminia,\gminib)\cap \del S=\emptyset$.
\end{itemize}
Such a sector is called strictly adapted to $(V,\nabla)$,
if moreover each $S(\gminia,\gminib)$ is non-empty.

According to the classical asymptotic analysis for meromorphic
flat bundles on curves,
if $S$ is adapted, 
then there exists a flat decomposition
\begin{equation}
\label{eq;12.10.8.2}
 (V,\nabla)_{|S}
=
 \bigoplus_{\gminia\in\nbigi(V)}
 (V_{\gminia,S},\nabla_{\gminia,S})
\end{equation}
which is asymptotic to (\ref{eq;12.10.8.1}).
If $S$ is strictly adapted,
then the decomposition (\ref{eq;12.10.8.2}) is unique.

A flat frame $\vecy=(y_i)$ of $(V,\nabla)_{|S}$ 
is called adapted,
if it is compatible with a decomposition (\ref{eq;12.10.8.2})
in the sense that,
for each $y_i$, there exists $\gminia\in\nbigi(V)$
such that $y_i\in V_{\gminia,S}$.
Let $S^{(i)}$ $(i=1,2)$ be adapted sectors 
such that $S^{(1)}\cap S^{(2)}\neq\emptyset$.
Suppose we are given adapted flat frames $\vecy^{(i)}$ of
$(V,\nabla)_{|S^{(i)}}$.
We obtain the matrix $A$ 
determined by
$\vecy^{(1)}=\vecy^{(2)}A$.
The matrix $A$ is called the Stokes matrix
or Stokes factor with respect to $\vecy^{(i)}$ $(i=1,2)$.
(If $S^{(1)}\cap S^{(2)}$ is not connected,
we consider the matrix for each connected component.)

\begin{rem}
Although we fix a coordinate $z$,
the above notions are independent of $z$.
It is better to replace
the index set $z^{-1}\cnum[z^{-1}]$
with $\nbigo_X(\ast D)/\nbigo_X$.
\hfill\qed
\end{rem}

\begin{rem}
For a general meromorphic flat bundle on $(X,D)$,
we take an appropriate ramified covering
$\varphi:(X,D)\lrarr (X,D)$ given by $\varphi(z)=z^e$
such that the pull back $\varphi^{\ast}(V,\nabla)$
is unramified,
for which the Stokes structure is considered as above.
\hfill\qed
\end{rem}

\begin{rem}
Let $C$ be a general complex curve
with a discrete subset $Z\subset C$.
Let $(V,\nabla)$ be a meromorphic flat bundle
on $(C,Z)$.
We can consider the Stokes structure
of $(V,\nabla)$ on a neighbourhood of
each point of $Z$.
\hfill\qed
\end{rem}

\paragraph{Lift of a formal frame}

Let $(V,\nabla)$ be an unramified meromorphic flat bundle
with the decomposition (\ref{eq;12.10.8.2}).
For each $\gminia\in\nbigi(V)$,
let $r(\gminia):=\rank(\Vhat_{\gminia})$.
We have a frame
$\vecvhat_{\gminia}$ of $\Vhat_{\gminia}$
with the following property:
\[
 \nabla\vecvhat_{\gminia}
=\vecvhat_{\gminia} 
 \Bigl(d\gminia I_{\gminia}+A_{\gminia}\frac{dz}{z}\Bigr)
\]
Here, $I_{\gminia}$ denote 
the $r(\gminia)$-th identity matrix,
and $A_{\gminia}\in M_{r(\gminia)}(\cnum)$.
They give a frame 
$\vecvhat=(\vecvhat_{\gminia}\,|\,\gminia\in\nbigi(V))$
of $V\otimes\cnum[\![z]\!]$
for which we have
$\nabla\vecvhat=\vecvhat
 \Bigl(
 \bigoplus\bigl(
 d\gminia I_{\gminia}+A_{\gminia}dz/z
 \bigr)
 \Bigr)$.
Recall the following.
\begin{lem}
Let $S$ be an adapted sector.
Take a flat decomposition {\rm(\ref{eq;12.10.8.2})}.
We have a holomorphic frame $\vecv_{\gminia,S}$
of $V_{\gminia,S}$ such that
(i) $\vecv_{\gminia,S}$ is asymptotic to $\vecvhat_{\gminia}$,
(ii) $\nabla_{\gminia,S}\vecv_{\gminia,S}=
  \vecv_{\gminia,S}(d\gminia I_{\gminia}+A_{\gminia}dz/z)$.
Such $\vecv_{\gminia,S}$ is uniquely determined,
once we fix the decomposition
{\rm(\ref{eq;12.10.8.2})}.
\hfill\qed
\end{lem}

In particular,
if $S$ is strictly adapted,
there uniquely exists a frame $\vecv_S$ of $V_{|S}$
such that
(i) $\vecv_S$ is asymptotic to $\vecvhat$,
(ii)  $\nabla\vecv_S=\vecv_S
 \Bigl(
 \bigoplus\bigl(
 d\gminia I_{\gminia}+A_{\gminia}dz/z
 \bigr)
 \Bigr)$.
If we choose a branch of $\log z$ on $S$,
we obtain an adapted flat frame
\begin{equation}
\label{eq;12.12.30.1}
 \vecvtilde_S:=
 \vecv_S\,\exp\Bigl(
 \bigoplus\bigl(
-\gminia\,I
-A_{\gminia}\log z
 \bigr)
 \Bigr).
\end{equation}

\paragraph{Pure slope}
Let $(V,\nabla)$ be an unramified meromorphic flat bundle
on $(X,D)$ with the decomposition (\ref{eq;12.10.8.2}).
Let $d$ be a positive integer.
In this paper, 
we say that $(V,\nabla)$ has pure slope $d$,
if (i) $\deg_{z^{-1}}(\gminia)=d$ for each $\gminia\in\nbigi(\nabla)$,
(ii) $\deg_{z^{-1}}(\gminia-\gminib)=d$ for 
each pair $\gminia,\gminib\in\nbigi(\nabla)$ 
with $\gminia\neq\gminib$.
We shall use the following lemma implicitly.
\begin{lem}
\label{lem;12.12.3.20}
Suppose that $(V,\nabla)$ has pure slope $d$.
Let $S$ be a sector adapted to $(V,\nabla)$
whose angle is larger than $\pi/d$.
Then, $S$ is strictly adapted.
In particular,
there exist strictly adapted sectors $S_i$
$(i=1,\ldots,N)$
such that $\bigcup S_i=X\setminus D$.
\hfill\qed
\end{lem}

Let $(V,\nabla)$ be a meromorphic flat bundle
on $(X,D)$,
which is not necessarily unramified.
Let $\mu$ be a positive rational number.
In this paper,
we say that $(V,\nabla)$ has pure slope $\mu$,
if there exists a ramified covering 
$\varphi:(X,D)\lrarr (X,D)$, $\varphi(z)=z^e$
such that 
$\varphi^{\ast}(V,\nabla)$ is unramified
with pure slope $e\mu\in\seisuu$.

\subsubsection{$R$-structure}

Let $X$ and $D$ be as in \S\ref{subsection;12.10.8.4}.
Let $(V,\nabla)$ be a meromorphic flat bundle on $(X,D)$.
Let $\nbigl$ be the local system on
$X\setminus D$ associated to $(V,\nabla)$,
i.e., it is the local system of flat sections of 
$(V,\nabla)_{|X\setminus D}$.
Let $R$ be a subring of $\cnum$.
(We will be mainly interested in the cases 
that $R$ is $\seisuu$ or $\rnum$.)
An $R$-structure of $(V,\nabla)_{|X\setminus D}$ is
an $R$-local system $\nbigl_{R}$ 
with an isomorphism
$\nbigl_{R}\otimes_R\cnum \simeq\nbigl$.
In this paper, an $R$-local system means
a locally constant sheaf of 
finitely generated projective $R$-modules.
We shall often regard
$\nbigl_R$ 
as a subsheaf of $\nbigl$.
We recall a compatibility condition 
of $R$-structure with Stokes structure
in the case that $(V,\nabla)$ has pure slope.

\begin{df}
Let $(V,\nabla)$ be a meromorphic flat bundle on $(X,D)$
such that it is unramified with pure slope.
We say that
an $R$-structure of $(V,\nabla)_{|X\setminus D}$
is compatible with the Stokes structure,
if the following holds:
\begin{itemize}
\item
Let $S$ be a sector strictly adapted to $(V,\nabla)$.
Then the decomposition {\rm(\ref{eq;12.10.8.2})}
of $(V,\nabla)_{|S}$ is induced by
a decomposition of $\nbigl_{R|S}$.
\end{itemize}
In that case,
the $R$-structure of $(V,\nabla)_{|X\setminus D}$
is called an $R$-structure of $(V,\nabla)$.
\hfill\qed
\end{df}

For each $\gminia\in\nbigi(\nabla)$,
we can take a meromorphic flat bundle
$\Gr^{\nbigf}_{\gminia}(V,\nabla)$ on $(X,D)$
with an isomorphism
$\Gr^{\nbigf}_{\gminia}(V,\nabla)
 \otimes_{\nbigo_X}\cnum[\![z]\!]
\simeq
 (\Vhat_{\gminia},\nablahat_{\gminia})$,
where $(\Vhat_{\gminia},\nablahat_{\gminia})$
are as in (\ref{eq;12.10.8.1}).
Recall that an $R$-structure of $(V,\nabla)$
induces an $R$-structure of
$\Gr^{\nbigf}_{\gminia}(V,\nabla)$.
Indeed, for any strictly adapted sectors $S$,
{\rm(\ref{eq;12.10.8.2})}
gives a unique isomorphism
$(V,\nabla)_{|S}
 \simeq
 \bigoplus
 \Gr^{\nbigf}_{\gminia}(V,\nabla)_{|S}$,
which induces an $R$-structure
of $\Gr^{\nbigf}_{\gminia}(V,\nabla)_{|S}$.
It is easy to observe that the $R$-structures
can be glued when the strictly adapted sectors
are varied.

\begin{lem}
\label{lem;12.12.6.220}
Let $(V,\nabla)$ be a meromorphic flat bundle on $(X,D)$
which is unramified with pure slope.
The following conditions are equivalent.
\begin{itemize}
\item
$(V,\nabla)$ has an $R$-structure $\nbigl_R$
such that 
the induced $R$-structure $\nbigl_{\gminia,R}$
of $\Gr^{\nbigf}_{\gminia}(V,\nabla)$
are local systems of free $R$-modules.
\item
There exists a frame $\vecvhat$
of $V\otimes\cnum[\![z]\!]$ such that,
for any strictly adapted sectors $S_i$ $(i=1,2)$
with $S_1\cap S_2\neq\emptyset$,
the entries of 
the Stokes matrix with respect to $\vecvtilde_{S_i}$ $(i=1,2)$
are elements of $R$.
(See {\rm(\ref{eq;12.12.30.1})} for $\vecvtilde_{S_i}$.)
\item
We can associate an adapted flat frame $\vecvtilde_S$
of $(V,\nabla)_{|S}$
to each strictly adapted sector $S$
such that,
if $S\cap S'\neq\emptyset$,
then the entries of the Stokes matrix
with respect to $\vecvtilde_S$ and $\vecvtilde_{S'}$
are elements of $R$.
\end{itemize}
\end{lem}
\pf
The second condition implies the third.
Suppose the third.
By using the adapted flat frames $\vecvtilde_{S}$ 
for any strictly adapted sectors  $S$,
we can define an $R$-structure of the local system.
It is clearly compatible with the Stokes structure
by construction,
i.e.,
the first condition is satisfied.

Let us prove that the first condition implies the second.
We take a holomorphic frame
$\vecv_{\gminia}$ of
$\Gr^{\nbigf}_{\gminia}(V,\nabla)$
such that
 $\vecv_{\gminia}=\vecv_{\gminia}
 \Bigl(
 d\gminia\,I_{r(\gminia)}
+A_{\gminia}\,dz/z
 \Bigr)$,
where $A_{\gminia}$ is a constant matrix.
We may assume that the multi-valued flat frame
$\vecv_{\gminia}\exp\bigl(
 -\gminia I_{r(\gminia)}-A_{\gminia}\log z
 \bigr)$
is a frame of 
$\nbigl_{\gminia\,R}$.
We have a natural isomorphism
$\Gr^{\nbigf}(V)\otimes\cnum[\![z]\!]
 \simeq
 V\otimes\cnum[\![z]\!]$.
Let $\vecvhat_{\gminia}$ be induced by
$\vecv_{\gminia}$
and the isomorphism.
For any strictly adapted sector $S$,
we take a flat frame $\vecvtilde_S$
from $\vecvhat=(\vecvhat_{\gminia})$.
It naturally gives a frame of 
$\nbigl_{R|S}$.
Then, the entries of 
the Stokes matrices with respect to
the lift of $\vecvhat$ are elements of $R$.
\hfill\qed

\begin{rem}
The conditions in Lemma {\rm\ref{lem;12.12.6.220}}
are equivalent to the existence of 
a covering $X\setminus D=\bigcup_{i=1}^N S_i$
by strictly adapted sectors,
such that 
(i) each $S_i$ is equipped with an adapted flat frame
$\vecv_{S_i}$ of $(V,\nabla)_{|S_i}$
 $\vecv_{S_i}$,
(ii) if $S_i\cap S_j\neq\emptyset$,
 the entries of the Stokes matrix 
 with respect to $\vecvtilde_{S_i}$ and $\vecvtilde_{S_j}$
 are elements of $R$.
It will be explained elsewhere.
\hfill\qed
\end{rem}

Let us consider the case where
$(V,\nabla)$ is not necessarily unramified
but has pure slope.
We have an appropriate ramified covering
$\varphi:(X',D')\lrarr (X,D)$
such that
$\varphi^{\ast}(V,\nabla)$ is unramified
with pure slope.
\begin{df}
We say that
an $R$-structure of $(V,\nabla)_{|X\setminus D}$
is an $R$-structure of $(V,\nabla)$,
if the induced $R$-structure of
$\varphi^{\ast}(V,\nabla)_{|X\setminus D}$
is an $R$-structure of
$\varphi^{\ast}(V,\nabla)$.
\hfill\qed
\end{df}

Suppose we have a decomposition
$(V,\nabla)=\bigoplus_{i=1}^s (V_i,\nabla_i)$
such that 
$\nbigi(\nabla_i)\cap\nbigi(\nabla_j)=\emptyset$
if $i\neq j$.
Let $\nbigl_i$ denote the local systems
corresponding to $(V_i,\nabla_i)$.

\begin{prop}
\label{prop;12.12.4.131}
If $(V,\nabla)$ has an $R$-structure $\nbigl_R$,
it is compatible with the decomposition
$\nbigl=\bigoplus\nbigl_i$.
Namely,
each $\nbigl_i$ has an $R$-structure $\nbigl_{Ri}$
compatible with the Stokes structure
such that $\nbigl_R=\bigoplus\nbigl_{Ri}$.
\end{prop}
\pf
Let $S$ be any sector strictly adapted to $(V,\nabla)$.
It is also adapted to each $(V_i,\nabla_i)$.
The decomposition (\ref{eq;12.10.8.2})
is compatible with the $R$-structure of
$\nbigl_R$.
Because (\ref{eq;12.10.8.2})
is a refinement of
the decomposition $V=\bigoplus V_i$,
we obtain that 
$\nbigl_R$ is compatible with
$V=\bigoplus V_i$.
\hfill\qed

\vspace{.1in}

Let us give a remark on the uniqueness 
of $R$-structure
in the case that
$R$ is a subfield of $\cnum$.
Let $(V,\nabla)$ be a meromorphic flat bundle
on $(X,D)$ with pure slope.
Let $\nbigl$ be the corresponding local system.
\begin{prop}
\label{prop;12.12.4.101}
Suppose that $(V,\nabla)$ is irreducible,
and that it has $R$-structures
$\nbigl_{Ri}$ $(i=1,2)$.
Then, 
there exists a complex number $\alpha$
such that
$\nbigl_{R1}=\alpha\,\nbigl_{R2}$
in $\nbigl$.
The number $\alpha$ is well defined in
$\cnum^{\ast}/\rnum^{\ast}$,
where $k^{\ast}$ denotes $k\setminus\{0\}$
for a field $k$.
\end{prop}
\pf
Let $\varphi:(X',D')\lrarr (X,D)$
be a ramified covering $\varphi(z)=z^e$
such that 
$\varphi^{\ast}(V,\nabla)$ is unramified.
It has the $R$-structures
$\varphi^{-1}\nbigl_{Ri}$,
which are compatible with 
the Stokes structure of $\varphi^{\ast}(V,\nabla)$.
Let $\Hom^{\St}(\nbigl_{R1},\nbigl_{R2})$
denote the space of $R$-homomorphisms
$f:\nbigl_{R1}\lrarr\nbigl_{R2}$
such that,
for any sectors $S\subset X'\setminus D'$
adapted $\varphi^{\ast}(V,\nabla)$,
the restriction $\varphi^{-1}f_{|S}$ 
preserves the decompositions of $\varphi^{-1}\nbigl_{Ri|S}$.
It is naturally an $R$-vector space.
Similarly,
let $\Hom^{\St}(\nbigl,\nbigl)$
denote the space of $\cnum$-homomorphisms
$f:\nbigl\lrarr\nbigl$
such that,
for any strictly adapted sectors $S\subset X'\setminus D'$,
the restriction $\varphi^{-1}f_{|S}$ preserves the decompositions
of $\varphi^{-1}\nbigl_{|S}$.
By the Riemann-Hilbert correspondence,
$\Hom^{\St}(\nbigl,\nbigl)$
is naturally isomorphic to the space of
flat morphisms
$(V,\nabla)\lrarr(V,\nabla)$.
By the irreducibility of $(V,\nabla)$,
we have
$\Hom^{\St}(\nbigl,\nbigl)=\cnum$.
The following naturally induced morphism 
is an isomorphism:
\[
 \Hom^{\St}(\nbigl_{R1},\nbigl_{R2})\otimes_R\cnum
\lrarr
 \Hom^{\St}(\nbigl,\nbigl)
\]
We take a non-zero element
$f\in\Hom^{\St}(\nbigl_{R1},\nbigl_{R2})$.
Then, there exists a complex number $\beta$
such that $\beta\,f=\id_{\nbigl}$.
It implies the claim of the proposition.
\hfill\qed

\subsubsection{Appendix: Stokes filtration}

We also have a formulation of $R$-structure
 in terms of the Stokes filtrations,
which works even in the case that $(V,\nabla)$
does not have pure slope,
although we shall not use it in this paper.
For an adapted sector $S$,
which is not necessarily strictly adapted,
the decomposition {\rm(\ref{eq;12.10.8.2})}
is not uniquely determined.
For simplicity, suppose $(V,\nabla)$ is unramified.
We consider the filtration $\nbigf^S$
given by
$\nbigf^{S}_{\gminia}(V_{|S})
=\bigoplus_{\gminib\leq_S\gminia}
 V_{\gminib,S}$.
Here, $\gminia\leq_S\gminib$ is defined to be
$-\Re(\gminia)\leq -\Re(\gminib)$ on $S$.
The filtration is given also by a growth order of
flat sections of $(V,\nabla)$.
In particular, it is independent of the choice of
the decomposition {\rm(\ref{eq;12.10.8.2})}.
It is called the Stokes filtration of $(V,\nabla)$.
It induces a filtration of $\nbigl_{|S}$
denoted by $\nbigf^S$.

\begin{rem}
Suppose that 
$(V,\nabla)$ is unramified with pure slope.
Then, 
$(V,\nabla)$ has an $R$-structure,
if and only if
the associated local system has an $R$-structure
$\nbigl_{R}$ satisfying the following conditions:
\begin{itemize}
\item
 For each adapted sector $S$,
 the filtration $\nbigf^S(\nbigl_{|S})$
 is induced by a filtration
 $\nbigf^S(\nbigl_{R|S})$ 
 of $\nbigl_{R|S}$.
\item
 The filtration $\nbigf^S(\nbigl_{R|S})$
 has a splitting,
 i.e.,
 $\nbigf^S_{\gminia}(\nbigl_{R|S})
\lrarr
 \nbigf^S_{\gminia}(\nbigl_{R|S})
 \big/
 \nbigf^S_{<\gminia}(\nbigl_{R|S})$
has a splitting.
\end{itemize}
Indeed, the ``only if'' part is clear.
We shall explain the ``if'' part elsewhere.
The condition for the Stokes filtration
makes sense,
even if $(V,\nabla)$ does not have pure slope.
\end{rem}

\subsection{Integral structure of
a certain meromorphic flat bundle}
\label{subsection;12.10.6.110}

\subsubsection{Statement}

Let $\veca=(a_i\,|\,i=1,\ldots,r)\in\cnum^r$.
We set
$P_{\veca}(T):=\prod_{i=1}^r(T-e^{2\pi\sqrt{-1}a_i})$.
Let $\diag[a_1,\ldots,a_r]$ denote
the diagonal matrix whose $(i,i)$-entries are $a_i$.
We consider the following meromorphic connection 
$\nabla_{\veca}$ of $\nbigv_r$:
\begin{equation}
\label{eq;12.10.2.1}
 \nabla_{\veca}\vece
=\vece\,\Bigl(
 \nbigk(r,1)-\diag[a_1,\ldots,a_r]
 \Bigr)
\frac{dq}{q}
\end{equation}
It is irregular singular at $\infty$,
and regular singular at $0$.
Let $R$ be any subring of $\cnum$,
and let $R^{\ast}$ denote the set of
the invertible elements in $R$.

\begin{thm}
\label{thm;12.10.6.40}
The meromorphic flat bundle
$(\nbigv_r,\nabla_{\veca})$
has an $R$-structure
if and only if 
$P_{\veca}(T)\in R[T]$
and $P_{\veca}(0)\in R^{\ast}$.
In that case,
the induced $R$-structure on
$\Gr^{\nbigf}_{\gminia}(\nbigv_r,\nabla_{\veca})$
are local systems of free $R$-modules.
(See the remark right before Lemma {\rm\ref{lem;12.12.6.220}}
for the induced $R$-structure of
$\Gr^{\nbigf}_{\gminia}(\nbigv_r,\nabla_{\veca})$.)
\end{thm}

Let us prove the ``only if'' part.
Suppose $(\nbigv_r,\nabla_{\veca})$ has
an $R$-structure $\nbigl_R$.
Let $R_f:=R[f^{-1}]$.
We may assume that $R$ is finitely generated over $\seisuu$.
Take any $f\in R$ such that
$\nbigl_R\otimes_RR_f$
is a locally constant sheaf of free $R_f$-module.
We have
$P_{\veca}(T)\in R_f[T]$,
because the polynomial
$P_{\veca}(T)$
is the characteristic polynomial of
the monodromy along the loop around $q=0$
in the counter clock-wise direction.
Because $\nbigl_R$ is a locally constant sheaf of
finitely generated projective $R$-modules,
we can take $f_i\in R$ $(i=1,\ldots,N)$
as above such that $\Spec(R)=\bigcup \Spec(R_{f_i})$.
Hence, we have
$P_{\veca}(T)\in R[T]$.

To prove the ``if'' part,
we shall observe that
the non-zero entries of the Stokes factors
with respect to some natural frames
can be described
in terms of the coefficients of the polynomial
$P_{\veca}(T)$
(Proposition \ref{prop;12.12.5.210}).

\vspace{.1in}

Before going to the proof,
we shall give some consequences.

\begin{cor}
\label{cor;12.12.26.1}
Let $\nbigr$ denote the subring of $\cnum$
generated by $e^{\pm 2\pi\sqrt{-1}a_i}$ $(i=1,\ldots,r)$
over $\seisuu$.
Then, 
$(\nbigv_r,\nabla_{\veca})$
has an $\nbigr$-structure.
\hfill\qed
\end{cor}

\begin{cor}
For the expression
$P_{\veca}(T)=T^r+\sum_{j=0}^{r-1} P_jT^j$,
let $\nbigr$ denote the subring of $\cnum$
generated by 
$P_i$ $(i=0,\ldots,r-1)$ and $P_0^{-1}$
over $\seisuu$,
Then, 
$(\nbigv_r,\nabla_{\veca})$
has an $\nbigr$-structure.
\end{cor}

\begin{cor}
\label{cor;12.12.4.33}
Suppose $\veca\in\rnum^r$.
Then, $(\nbigv_r,\nabla_{\veca})$
has a $\rnum$-structure
if and only if 
$P_{\veca}(T)\in\seisuu[T]$.
In that case,
$(\nbigv_r,\nabla_{\veca})$
has a $\seisuu$-structure.
\end{cor}
\pf
If $(\nbigv_r,\nabla_{\veca})$
has a $\rnum$-structure,
we have
$P_{\veca}(T)\in\rnum[T]$.
Because $e^{2\pi\sqrt{-1}a_i}$ are algebraic integers,
we obtain $P_{\veca}(T)\in\seisuu[T]$.
The converse is clear from Theorem 
\ref{thm;12.10.6.40}.
\hfill\qed

\subsubsection{Preliminary}
\label{subsection;12.12.6.1}

In the following,
$\nabla_{\veca}$ will be denoted by $\nabla$.
Let $\varphi:\proj_w^1\lrarr\proj^1_q$
be the ramified covering given by
$\varphi(w)=w^r$.
The pull back 
$\varphi^{\ast}(\nbigv_r,\nabla)$
is denoted by $(\nbigvcirc,\nablacirc)$.
We put $\tau:=\exp(2\pi\sqrt{-1}/r)$.
We set
$v_j:=\sum_i \tau^{ij}w^{-i+1}\varphi^{\ast}e_i$.
\begin{lem}
We have
$\nablacirc\vecv=\vecv\left(
 \diag[\tau^{-1},\ldots,\tau^{-r}]\,r\,dw
+B\,dw/w
 \right)$
for some constant matrix $B$.
For any $i$,
the $(i,i)$-entry of $B$ is
$b:=-\sum a_j-(r-1)/2$.
\end{lem}
\pf
It can be checked by a direct computation.
We give an indication.
Let $v_i':=w^{-i+1}e_i$.
Let $C$ be a cyclic matrix
such that
the $(i,j)$-entry is $1$ if $i-j\equiv 1$ modulo $r$,
or $0$ otherwise.
Then, we have
\[
\nablacirc\vecv'=\vecv'
\Bigl(
 Crdw
-\bigl(
 \diag[a_1,\ldots,a_r]r+
 [0,1,\ldots,r-1]
\bigr)dw/w
\Bigr).
\]
Let $\Upsilon$ be the $r\times r$-matrix whose 
$(i,j)$-entries are $\tau^{ij}$.
We have
$\vecv=\vecv'\,\Upsilon$.
We have
$\Upsilon^{-1}=r^{-1}\,(\tau^{-ij})$.
We set
$B:=-\Upsilon^{-1}
 \bigl(
r\diag[a_1,\ldots,a_r]
+[0,1,\ldots,r-1]
 \bigr)
 \Upsilon$.
Then, we obtain the formula for $\nablacirc\vecv$.
Moreover,
the $(i,i)$-entries of $B$ are 
$\sum r^{-1}\tau^{-ij}(-a_jr+1-j)\tau^{ji}
=-\sum a_j-(r-1)/2$.
\hfill\qed

\vspace{.1in}

By a standard argument as in \cite{levelt},
at the formal completion $\inftyhat$,
we have a formal frame $\vecphat$ such that
\begin{equation}
 \label{eq;12.12.5.112}
\nablacirc \vecphat=\vecphat
 \left(
 \diag[\tau^{-1},\ldots,\tau^{-r}]r\,dw
+\diag[b,\ldots,b]\frac{dw}{w}
 \right)
\end{equation}
In particular, we have
$\nbigi(\nablacirc)=
 \bigl\{
 r\tau^{-i}w\,\big|\,
 i=1,\ldots,r
 \bigr\}$.
We set 
$\gminia_i:=r\tau^{-i}w$.
We use the convention
$\gminia_{i}=\gminia_{r+i}$.

Let $\mu_r:=\bigl\{\kappa\in\cnum\,\big|\,\kappa^{r}=1\bigr\}$.
It is regarded as the Galois group of 
the ramified covering $\varphi$
with the action
$(\kappa,w)\longmapsto \kappa w$.
The meromorphic flat bundle
$(\nbigvcirc,\nablacirc)$ is naturally $\mu_r$-equivariant.
We have the ambiguity of the choice of $\vecphat$,
i.e.,
the scalar multiplication of each $\phat_i$.
We have
$(\tau^{-1})^{\ast}\phat_j=\xi_{j}\phat_{j+1}$
for some $\xi_j\in\cnum^{\ast}$
with respect to the action.
We have $\prod \xi_j=1$.

\subsubsection{Flat frames}
\label{subsection;12.10.2.2}

We consider the following sectors:
\begin{equation}
 \label{eq;12.12.5.100}
S_1:=\bigl\{ 
 w\,\big|\,
 \bigl|\arg(w)-(\pi/2-\epsilon)\bigr|<\pi/2 
 \bigr\} 
\quad\quad
S_2:=\bigl\{ 
 w\,\big|\,
 \bigl|\arg(w)-(\pi/2+2\pi/r-\epsilon)\bigr|<\pi/2 
 \bigr\}
\end{equation}
Here, $\epsilon$ is a sufficiently small positive number.
We have $\tau(S_1)=S_2$.
We put $\vecS=S_1\cup S_2\cup\{\arg(w)=2\pi/r-\epsilon\}$.
The last term is necessary only in the case $r=2$.

It is easy to see that
$S_i$ $(i=1,2)$ are strictly adapted to
the connection $\nablacirc$.
We have the following flat decompositions
as in (\ref{eq;12.10.8.2}):
\begin{equation}
\label{eq;12.12.5.110}
 \nbigvcirc_{|S_i}
=\bigoplus_{\gminia\in\nbigi(\nablacirc)} \nbigvcirc_{S_i,\gminia}
\end{equation}

We set $\omega:=e^{2\pi\sqrt{-1}b}\in R^{\ast}$.
Let $C_1=\bigl((C_1)_{ij}\bigr)$ be the $r\times r$-matrix
such that $(C_1)_{ij}=1$ if $i-j=1$,
$(C_1)_{ij}=\omega$ if $(i,j)=(r,1)$,
and $(C_1)_{ij}=0$ otherwise.

\begin{lem}
\label{lem;12.12.5.120}
There exist flat frames
$\vecptilde_{S_i}=(\ptilde_{j,S_i}\,|\,j=1,\ldots,r)$
of $\nbigvcirc_{|S_i}$
such that
(i)  $\ptilde_{j,S_i}$ are flat sections of
$\nbigv^{\circ}_{S_i,\gminia_j}$,
(ii) $(\tau^{-1})^{\ast}\vecptilde_{S_1}=\vecptilde_{S_2}C_1$.
\end{lem}
\pf
We put $\omega_1:=e^{2\pi\sqrt{-1}b/r}$.
Let $\vecphat$ be a formal frame of
$\nbigvcirc$ at $\infty$ satisfying (\ref{eq;12.12.5.112}).
By adjusting the scalar multiplication,
we may assume that
$\phat_j=\omega_1^{j-1}\cdot (\tau^{1-j})^{\ast}\phat_1$.
Then, we have
$(\tau^{-1})^{\ast}\vecphat
=\vecphat C_1\omega_1^{-1}$.

On any strictly adapted sector $S$,
we have a unique holomorphic frame
$\vecp_S$ of $V$ 
such that
$\vecp_S$ is asymptotic to $\vecphat$,
and 
$\nabla\vecp_S
=\vecp_S\Bigl(
 \diag[\tau^{-1},\ldots,\tau^{-r}]\,r\,dw
+\diag[b,\ldots,b]dw/w
 \Bigr)$ 
holds.
We fix a branch of $\log w$ on $S_1\cup S_2$.
If $r=2$, we impose that
it is analytic on $\vecS$.
We define
\[
 \vecptilde_{S_i}:=
 \vecp_{S_i}\,\exp\bigl(
 \diag[-\tau^{-1},\ldots,-\tau^{-r}]rw
 \bigr)\exp(-b\log w),
\]
which gives a flat frame on $S$.
Then, we have the following equalities:
\begin{multline}
 (\tau^{-1})^{\ast}\vecptilde_{S_1}
=(\tau^{-1})^{\ast}\vecp_{S_1}
 (\tau^{-1})^{\ast}
 \Bigl(
 \exp\bigl(
 \diag[-\tau^{-1}rw,\ldots,-\tau^{-r}rw]
 \bigr)
 \exp(-b\log w)
\Bigr)
 \\
=\vecp_{S_2}C_1\omega_1^{-1}
   \exp\bigl(
 \diag[-\tau^{-2}rw,\ldots,-\tau^{-r-1}rw]
 \bigr)
  \exp(-b\log w)\omega_1
\\
=\vecptilde_{S_2}
  \exp\bigl(
 \diag[\tau^{-1}rw,\ldots,\tau^{-r}rw]
 \bigr)\,
  C_1\,
  \exp\bigl(
 \diag[-\tau^{-2}rw,\ldots,-\tau^{-r-1}rw]
 \bigr)
=\vecptilde_{S_2}\,C_1
\end{multline}
Thus, the proof of Lemma \ref{lem;12.12.5.120}
is finished.
\hfill\qed

\vspace{.1in}

To adjust the signature in the later formulas,
we set 
\[
 y_{j,S_i}:=
 \left\{
 \begin{array}{ll}
 \omega\,\ptilde_{j,S_i} & (j=1,\ldots,[(r-1)/2])\\
 \ptilde_{j,S_i} & (j=[(r-1)/2]+1,\ldots,r).
 \end{array}
 \right.
\]
We define $y_{j,S_i}$ for any $j\in\seisuu$
by using the convention
$y_{j+\ell r,S_i}=y_{j,S_i}$.

\begin{lem}
\label{lem;12.12.5.201}
If $r=2m$,
there exist complex numbers
$\alpha_j$ $(j=1,\ldots,m-1)$
and $\beta_j$ $(j=1,\ldots,m)$
such that the following holds:
\[
 y_{j,S_1}
=\left\{
 \begin{array}{ll}
 y_{j,S_2} 
 & (0\leq j\leq m-1) \\
 y_{j,S_2}+\alpha_{-j}y_{-j,S_2}
 +\beta_{-j}y_{-j-1,S_2} 
 & (-m<j<0) \\
 y_{-m,S_2}+\beta_my_{m-1,S_2}
 & (j=-m)
 \end{array}
 \right.
\]
If $r=2m+1$,
there exist complex numbers
$\alpha_j$ 
and $\beta_j$
$(j=1,\ldots,m)$
such that the following holds:
\[
 y_{j,S_1}
=\left\{
 \begin{array}{ll}
 y_{j,S_2} 
 & (0\leq j\leq m)
 \\
  y_{j,S_2}+\alpha_{-j}y_{-j,S_2}
 +\beta_{-j}y_{-j-1,S_2}
 & (-m\leq j<0)
 \end{array}
 \right.
\]
\end{lem}
\pf
On the half line $\arg(w)=0$, we have
\[
 \Re(-\tau^{-\ell}w)=\Re(-\tau^{\ell}w)
\quad
(\ell=1,\ldots,[(r-1)/2])
\]
On the half line $\arg(w)=\pi/r$, we have
\[
 \Re(-\tau^{-\ell}w)
=\Re(-\tau^{\ell+1}w)
\quad
(\ell=0,1,\ldots,[r/2]-1)
\]
On $S_1\cap S_2$,
we have
$\Re(-\tau^{-\ell}w)<\Re(-\tau^{\ell}w)$
and
$\Re(-\tau^{-\ell}w)<\Re(-\tau^{\ell+1}w)$
for $\ell$ as above.
Then, the claim of the lemma follows.
\hfill\qed

\subsubsection{Stokes factor}

We set
$x_{j,S_i}:=y_{-[r/2]+j-1,S_i}$
and $\vecx_{S_i}=(x_{1,S_i},\ldots,x_{r,S_i})$.
By construction,
we have
\begin{equation}
\label{eq;12.12.5.200}
(\tau^{-1})^{\ast}\vecx_{S_1}=\vecx_{S_2}C_1.
\end{equation}
We can naturally regard
$\vecx_{S_i}$ as tuples of flat sections
on $\vecS$.
Let $A$ be the matrix determined by
$\vecx_{S_1}=\vecx_{S_2}A$.
If $r=2m$, we have the following:
\begin{equation}
 A_{k\ell}=
 \left\{
 \begin{array}{ll}
 1  & (k=\ell)\\
 \beta_j & (k,\ell)=(m+j,m-j+1),\quad (1\leq j\leq m) \\
 \alpha_j & (k,\ell)=(m+j+1,m-j+1),\quad(1\leq j\leq m-1) \\
 0 & \mbox{(otherwise)}
 \end{array}
 \right.
\end{equation}
If $r=2m+1$,
we have the following:
\begin{equation}
 A_{k\ell}=
 \left\{
 \begin{array}{ll}
 1  & (k=\ell)\\
 \beta_j & (k,\ell)=(m+j,m-j+1),\quad(1\leq j\leq m) \\
 \alpha_j & (k,\ell)=(m+j+1,m-j+1),\quad (1\leq j\leq m) \\
 0 & \mbox{(otherwise)}
 \end{array}
 \right.
\end{equation}

For example,
if $r=4$,
we have the following relations:
\[
 \vecx_{S_1}
=\vecx_{S_2}
 \left(
 \begin{array}{cccc}
 1 & 0 & 0 & 0 \\
 0 & 1 & 0 & 0 \\
 0 & \beta_1 & 1 &0 \\
 \beta_2 &\alpha_1 & 0 & 1
 \end{array}
 \right)
\quad\quad
 (\tau^{-1})^{\ast}\vecx_{S_1}
=\vecx_{S_2}
 \left(
 \begin{array}{cccc}
 0 & 0 & 0 & \omega\\
 1 & 0 & 0 & 0 \\
 0 & 1 & 0 & 0 \\
 0 & 0 & 1 & 0
 \end{array}
 \right)
\]
If $r=5$,
we have the following relation:
\[
\vecx_{S_1}
=\vecx_{S_2}
 \left(
 \begin{array}{ccccc}
 1 & 0 & 0 & 0 & 0\\
 0 & 1 & 0 & 0 & 0 \\
 0 & \beta_1 & 1 &0 & 0\\
 \beta_2 &\alpha_1 & 0 & 1 & 0 \\
 \alpha_2 & 0 & 0 & 0 & 1
 \end{array}
 \right)
\quad\quad
 (\tau^{-1})^{\ast}\vecx_{S_1}
=\vecx_{S_2}
 \left(
 \begin{array}{ccccc}
 0 & 0 & 0 & 0 & \omega\\
 1 & 0 & 0 & 0 & 0 \\
 0 & 1 & 0 & 0 & 0 \\
 0 & 0 & 1 & 0 & 0 \\
 0 & 0 & 0 & 1 & 0
 \end{array}
 \right)
\]

\subsubsection{The Stokes factor and the monodromy}
\label{subsection;12.12.6.2}

The following proposition means that
the non-zero entries of the Stokes factor $A$
is described by the coefficients
of the polynomial $P_{\veca}(T)$.

\begin{prop}
\label{prop;12.12.5.210}
We have the following formula:
\[
 P_{\veca}(T)=T^r
-\sum_{j=1}^{[(r-1)/2]}\alpha_jT^{r-2j}
-\sum_{j=1}^{[r/2]}\beta_jT^{r-2j+1}
-\omega
\]
\end{prop}
\pf
We naturally regard 
$\vecx_{S_1}$ and $(\tau^{-1})^{\ast}\vecx_{S_1}$
as tuples of flat sections of
$(\nbigv^{\circ},\nabla^{\circ})_{|\vecS}$.
The monodromy of $(\nbigv_r,\nabla)$
along the loop in $\cnum_q^{\ast}$ around $0$
in the counter-clockwise direction is the relation between
$\vecx_{S_1}$ and $(\tau^{-1})^{\ast}\vecx_{S_1}$,
and it is expressed by 
the matrix $M$ determined by
$\vecx_{S_1}=
(\tau^{-1})^{\ast}\vecx_{S_1}M$.
The characteristic polynomial of $M$
is $P_{\veca}(T)$.
By the definition of $A$ and the relation (\ref{eq;12.12.5.200}),
we obtain $M=C_1^{-1}A$.
In particular,
the characteristic polynomial 
$\nbigp_{C_1^{-1}A}(T)$ of
$C_1^{-1}A$ is equal to $P_{\veca}(T)$.
It is elementary to prove the following formulas:
\[
 \frac{\del}{\del\beta_j}
 \nbigp_{C_1^{-1}A}(T)
=-T^{r+1-2j}
\quad\quad
 \frac{\del}{\del\alpha_j}
 \nbigp_{C_1^{-1}A}(T)
=-T^{r-2j}
\]
If $\alpha_j=0$ and $\beta_j=0$ for any $j$,
we have $\nbigp_{C_1^{-1}A}(T)=T^r-\omega$.
Then, we obtain the claim of
Proposition  \ref{prop;12.12.5.210}.
\hfill\qed

\subsubsection{End of the proof of Theorem \ref{thm;12.10.6.40}}

Suppose $P_{\veca}(T)\in R[T]$.
We obtain that the entries of $A$ 
are elements of $R$
by Proposition \ref{prop;12.12.5.210}.
We also have $\omega\in R^{\ast}$.
We set $S_{3/2}:=
 \bigl\{
 w\,\big|\,
 \bigl|
 \arg(w)-(\pi/2+\pi/r-\epsilon)
 \bigr|<\pi/2
 \bigr\}$.
It is standard to construct an adapted flat frame
$\vecx_{S_{3/2}}$ by using the factorization of
the matrix $A$.
Let $A^{(1)}$ and $A^{(2)}$ be the matrix as follows:
\[
 A^{(1)}_{k\ell}=\left\{
 \begin{array}{ll}
 1 & (k=\ell)\\
 \beta_j & (k,\ell)=([r/2]+j,[r/2]-j+1),\quad
(1\leq j\leq [r/2])\\
 0 &\mbox{\rm otherwise}
 \end{array}
 \right.
\]
\[
 A^{(2)}_{k\ell}=\left\{
 \begin{array}{ll}
 1 & (k=\ell)\\
 \alpha_j & (k,\ell)=([r/2]+j+1,[r/2]-j+1),\quad(1\leq j\leq [(r-1)/2])\\
 \end{array}
 \right.
\]
We set $\vecx_{S_{3/2}}:=\vecx_{S_2}A^{(2)}$.
Then, we have $\vecx_{S_1}=\vecx_{S_{3/2}}A^{(1)}$.
The entries of $A^{(i)}$ and $(A^{(i)})^{-1}$
are integers.

For $\ell=j+1$ or $j+3/2$ $(j=0,\ldots,r-1)$,
we set $S_{\ell}:=\tau^{j}(S_{\ell-j})$.
They are strictly adapted to $(\nbigvcirc,\nablacirc)$.
We have adapted frames
$\vecx_{S_{\ell}}:=
 (\tau^{-j})^{\ast}\vecx_{S_{\ell-j}}C_1^{-j}$ 
of $(\nbigvcirc,\nablacirc)_{|S_{\ell}}$.
The Stokes factors with respect to the frames
$\vecx_{S_{j+3/2}}$ and $\vecx_{S_{j+1}}$
is $C_1^{j}A^{(1)}C_1^{-j}$.
The Stokes factors with respect to the frames
$\vecx_{S_{j+2}}$ and $\vecx_{S_{j+3/2}}$
is $C_1^{j}A^{(2)}C_1^{-j}$.
We remark that
$(\tau^{-r})^{\ast}\vecx_{S_{\ell}}C_1^{-r}=
\vecx_{S_{\ell}}$,
where we care the analytic continuation of $\log w$.
Any strictly adapted sector can be 
deformed to one of $S_{\ell}$ preserving
the strictly adaptedness.
Therefore, we obtain Theorem \ref{thm;12.10.6.40},
according to Lemma \ref{lem;12.12.6.220}.
\hfill\qed

\subsection{Integral structure on the pull back}

\subsubsection{Statement}

We continue to use the notation in 
\S\ref{subsection;12.10.6.110}.
Suppose that $e^{\pm 2\pi\sqrt{-1}a_i}$ are algebraic integers.
Let $m$ be a positive integer.
We consider the following meromorphic connection
$\nabla^{(m)}_{\veca}$ of $\nbigv_r$:
\begin{equation}
\label{eq;12.12.26.3}
 \nabla^{(m)}_{\veca}\vece
=\vece\,\Bigl(
 \nbigk(r,m)-\diag[a_1,\ldots,a_r]
 \Bigr)
m\frac{dq}{q}
\end{equation}
Let $K\subset\cnum$ be an algebraic number field.
Let $R$ be the ring of integers in $K$.

\begin{thm}
\label{thm;12.12.26.12}
If $(\nbigv_r,\nabla^{(m)}_{\veca})$
is irreducible,
the following conditions are equivalent.
\begin{itemize}
\item
 $(\nbigv_r,\nabla^{(m)}_{\veca})$ has an $R$-structure.
\item
 $(\nbigv_r,\nabla^{(m)}_{\veca})$ has a $K$-structure.
\item
There exists $\gamma\in\cnum^{\ast}$ such that
(i) $\gamma^m\in K$,
(ii) $\prod_{i=1}^r(T-\gamma e^{2\pi\sqrt{-1}a_i})\in K[T]$.
\end{itemize}
\end{thm}
\pf
The first condition clearly implies the second.
Let us show that the second condition implies the first.
Let $\nbigk\subset\cnum$ be 
a sufficiently large algebraic number field
which contains $K$ and $e^{2\pi\sqrt{-1}a_i}$ $(i=1,\ldots,r)$.
Let $\nbigr$ denote the ring of integers in $\nbigk$.
Because
$\prod_{i=1}^r(T-e^{2\pi\sqrt{-1}a_i})\in\nbigr[T]$,
according to Theorem \ref{thm;12.10.6.40},
$(\nbigv_r,\nabla_{\veca})$ has an $\nbigr$-structure.
Because $(\nbigv_r,\nabla^{(m)}_{\veca})$
is isomorphic to the pull back of
$(\nbigv_r,\nabla_{\veca})$ by
$\varphi_m(q)=q^m$,
$(\nbigv_r,\nabla^{(m)}_{\veca})$
also has an $\nbigr$-structure,
denoted by $\nbigl_{\nbigr}$.
Let $\nbigl_{K}$ be a $K$-structure of
$(\nbigv_r,\nabla^{(m)}_{\veca})$
in the second condition.
By Proposition \ref{prop;12.12.4.101},
we may assume 
$\nbigl_{\nbigr}\otimes_{\nbigr}\nbigk
=\nbigl_{K}\otimes_{K}\nbigk$
by adjusting $\nbigl_{K}$.
Then, we obtain an $R$-structure 
$\nbigl_{R}$ of 
$(\nbigv_r,\nabla^{(m)}_{\veca})$
by using the following lemma.
\begin{lem}
\label{lem;12.12.26.10}
Let $V_{K}$ be a finite dimensional $K$-vector space.
Suppose that we are given 
an $\nbigr$-lattice $V_{\nbigr}$ of
$V_{\nbigk}:=V_{K}\otimes_{K}\nbigk$.
Let $\iota:V_{K}\lrarr V_{\nbigk}$ denote 
the natural inclusion.
Then, $V_{R}:=\iota^{-1}(V_{\nbigr})$
is an $R$-lattice of $V_{K}$.
Moreover,
if we are given decompositions
$V_{K}=\bigoplus_{i\in \nbigi} V_{K i}$
and 
$V_{\nbigr}=\bigoplus_{i\in\nbigi} V_{\nbigr i}$
such that
$V_{K i}\otimes_{K} \nbigk
=V_{\nbigr i}\otimes_{\nbigr}\nbigk$,
then we have an induced decomposition 
$V_{R}
=\bigoplus_{i\in I} V_{R i}$
such that
$V_{R i}\otimes_R K=V_{K i}$
and 
$V_{R i}\otimes \nbigr=V_{\nbigr i}$.
\end{lem}
\pf
If we regard $V_{\nbigr}$ as an $R$-module
by $R\subset\nbigr$,
then it is finitely generated over the Noetherian ring $R$.
Because $V_{R}$ is naturally an $R$-submodule
of $V_{\nbigr}$,
it is finitely generated over $R$.
It is torsion-free.
Because $R$ is a one dimensional regular ring,
$V_{R}$ is projective.
Let $j:V_{R}\lrarr V_{K}$ be the natural inclusion.
The induced morphism
$j_{K}:V_{R}\otimes_{R} K\lrarr V_{K}$
is injective.
For any $f\in V_{K}$,
there exists an integer $n_f$ such that
$n_f\iota(f)\in V_{\nbigr}$.
Hence, $j_{K}$ is surjective.
Thus, we obtain the first claim.
We can easily obtain the second claim 
by using the first claim.
Thus, Lemma \ref{lem;12.12.26.10}
is proved.
\hfill\qed

\vspace{.1in}

Let us show the equivalence of the second and third
conditions in Theorem \ref{thm;12.12.26.12}.
We put $\kappa:=\exp(2\pi\sqrt{-1}/m)$
and $\mu_m:=\bigl\{\kappa^j\,|\,j=0,1,\ldots,m-1\bigr\}$.
We consider the $\mu_m$-action on $\proj^1$
given by $(\kappa^j,q)\longmapsto \kappa^j q$.
We have a $\mu_m$-action on $\nbigv_r$
given by
$\kappa^{\ast}e_i\longleftrightarrow e_i$.
The connection $\nabla^{(m)}_{\veca}$
is equivariant with respect to the action.

For any $\nu\in\cnum$,
let $L(\nu)$ be the meromorphic flat bundle
$\nbigo_{\proj^1}(\ast D)\,y_{\nu}$
with $\nabla y_{\nu}=y_{\nu}\,\nu dq/q$.
The multi-valued flat section is given by
$\rho_{\nu}:=y_{\nu}\,\exp(-\nu\log q)$.
We consider a $\mu_m$-action on $L(\nu)$
given by $\kappa^{\ast}y_{\nu}=y_{\nu}$.
The tensor product
$(\nbigv_r,\nabla^{(m)}_{\veca})\otimes L(\nu)$
is equipped with an induced $\mu_m$-action.
Its descent with respect to the $\mu_m$-action is isomorphic to 
$(\nbigv,\nabla')$,
where $\nabla'$ is given as follows:
\begin{equation}
 \label{eq;12.12.26.11}
 \nabla'\vece
=\vece\,
 \Bigl(
 \nbigk(r,1)-\diag[a_1-\nu/m,\ldots,a_r-\nu/m]
 \Bigr)\frac{dq}{q}
\end{equation}

Suppose the second condition.
Let $\nbigl_{K}$ be a $K$-structure
of $(\nbigv_r,\nabla^{(m)}_{\veca})$.
By Proposition \ref{prop;12.12.4.101},
there exists a $\gamma\in\cnum^{\ast}$ such that
$(\kappa^{-1})^{\ast}\nbigl_{K}
=\gamma\nbigl_{K}$.
We have $\gamma^m\in K^{\ast}$.
We take $\nu$ such that 
$\exp(-2\pi\sqrt{-1}\nu/m)=\gamma$.
Then, $L(\nu)$ has a $K$-structure
$\nbigl_{\nu}$ given by 
$K\,\rho_{\nu}\subset\nbigl_{\nu}$,
and we have
$(\kappa^{-1})^{\ast}\nbigl_{\nu}
=\gamma^{-1}\nbigl_{\nu}$.
Because the induced $K$-structure 
$\nbigl_{K}\otimes\nbigl_{\mu}$ is $\mu_m$-equivariant,
we obtain a $K$-structure of
$(\nbigv_r,\nabla')$.
Applying Theorem \ref{thm;12.10.6.40},
we obtain
$\prod_{i=1}^r(T-\gamma e^{2\pi\sqrt{-1}a_i})\in K[T]$,
i.e.,
the third condition holds.

Suppose the third condition.
We take $\nu\in\cnum$ such that
$\exp(-2\pi\sqrt{-1}\nu/m)=\gamma$.
By Theorem \ref{thm;12.10.6.40},
$(\nbigv_r,\nabla')$ has a $K$-structure.
It induces a $K$-structure on
$(\nbigv_r,\nabla^{(m)}_{\veca})\otimes L(\nu)$,
and hence a $K$-structure on $(\nbigv_r,\nabla^{(m)}_{\veca})$.
Thus, Theorem \ref{thm;12.12.26.12} is proved.
\hfill\qed

\vspace{.1in}
We have a special case of Theorem \ref{thm;12.12.26.12}.
\begin{cor}
Suppose that $\veca\in\rnum^r$.
Then, the following conditions are equivalent.
\begin{itemize}
\item
 $(\nbigv_r,\nabla^{(m)}_{\veca})$
 has a $\seisuu$-structure.
\item
 $(\nbigv_r,\nabla^{(m)}_{\veca})$
 has a $\rnum$-structure.
\item
 There exists $\gamma\in\cnum^{\ast}$ such that
 (i) $\gamma^m\in \rnum$,
 (ii) $\prod_{i=1}^r(T-\gamma e^{2\pi\sqrt{-1}a_i})\in \rnum[T]$.
\hfill\qed
\end{itemize}
\end{cor}

\subsubsection{Examples}

Let us consider the conditions in the case that 
$r=m$ is a prime number,
and that $K=\rnum$ and $R=\seisuu$.
We begin with a preliminary.
Take $\gamma\in\real_{>0}\setminus\rnum$
such that $\gamma^r\in\rnum$.
Let $P$ be the minimal polynomial of
$\eta:=\gamma e^{2\pi\sqrt{-1}v/(r^mu)}$,
where $u$ and $v$ are integers
such that
$\gcd(u,v)=\gcd(r,u)=\gcd(r,v)=1$.

\begin{lem}
\label{lem;12.12.27.10}
Suppose $\deg(P)\leq r$.
\begin{itemize}
\item
If $r$ is an odd prime,
$P=T^r\pm \gamma^r$.
\item
If $r=2$,
we have $P=T^2\pm\gamma^2$
unless $(u,m)=(1,3),(3,2)$.
\item
If $r=2$ and $(u,m)=(1,3)$,
we have $\gamma\in\sqrt{2}\rnum_{>0}$
and $P(T)\in \bigl\{
 T^2+2bT+2b^2\,\big|\, 
 b\in\rnum
 \bigr\}$.
\item
If $r=2$ and $(u,m)=(3,2)$,
we have $\gamma\in\sqrt{3}\rnum_{>0}$
and $P(T)\in 
 \bigl\{
 T^2+3b T+3b^2\,\big|\,b\in\rnum
 \bigr\}$.
\end{itemize}
\end{lem}
\pf
We set 
$K:=\rnum(\eta)$.
By the hypothesis,
we have $[K:\rnum]\leq r$,
where $[K_1:K_2]$ denotes the degree
for a finite extension $K_1$ over $K_2$.
Because
$\eta^{r^{m+1}}
=\gamma^{r^{m+1}}e^{2\pi\sqrt{-1}(rv)/u}
 \in K$,
we have
$e^{2\pi\sqrt{-1}/u}\in K$.

Let us consider the case that $r$ is an odd prime.
Suppose $m\geq 2$,
and we shall derive a contradiction.
We have
$\eta^{r^{m-1}u}
=\gamma^{r^{m-1}u}e^{2\pi\sqrt{-1}v/r}\in K$,
and hence
$e^{2\pi\sqrt{-1}/r}\in K$.
Because $[\rnum(e^{2\pi\sqrt{-1}/r}):\rnum]=r-1$
and $[K:\rnum]\leq r$,
we obtain
$K=\rnum(e^{2\pi\sqrt{-1}/r})$
and $\deg P=r-1$.
Because
$\prod_{i=0}^{r^mu-1}\bigl(
 T-\gamma e^{2\pi\sqrt{-1}i/r^mu}
\bigr)
=T^{r^mu}-\gamma^{r^mu}\in\rnum[T]$,
a number conjugate to $\eta$ over $\rnum$
is of the form
$\gamma e^{2\pi\sqrt{-1}c}$
for $c\in\rnum$.
We have
$P=\prod_{j=1}^{r-1}(T-\gamma e^{2\pi\sqrt{-1}c_j})$
for some $c_j\in\rnum$ $(j=1,\ldots,r-1)$,
and hence
$P(0)=(-1)^{r-1}\gamma^{r-1}e^{2\pi\sqrt{-1}\sum c_j}$.
If $e^{2\pi\sqrt{-1}\sum c_j}\not\in\real$,
we have
$P(0)\not\in\rnum$ by considering the polar part.
If $e^{2\pi\sqrt{-1}\sum c_j}\in\real$,
then 
$e^{2\pi\sqrt{-1}\sum c_j}=\pm 1$,
and hence $P(0)\not\in\rnum$.
It contradicts with the hypothesis
$P(T)\in\rnum[T]$.
Thus, we obtain $m\leq 1$.

If $m\leq 1$,
we have 
$\gamma e^{2\pi\sqrt{-1}t/r}\in K$
for some integer $t$.
Recall that
$T^r-\gamma^r$ is irreducible in $\rnum[T]$.
Let $\overline{\rnum}$ be an algebraic closure of $\rnum$,
and let $f:\overline{\rnum}\lrarr\overline{\rnum}$
be an automorphism of the field over $\rnum$
such that
$f(\gamma e^{2\pi\sqrt{-1}t/r})=\gamma$.
Then,
$f(e^{2\pi\sqrt{-1}/u})\in K(\gamma)\subset\real$.
We obtain $u=1,2$,
and $\eta=\pm \gamma e^{2\pi\sqrt{-1}v/r}$.
Then, its minimal polynomial is
$T^r\pm \gamma^r$.

\vspace{.1in}
Let us consider the case $r=2$.
If $m=0$,
we have $\gamma\in K$
and $e^{2\pi\sqrt{-1}v/u}\in K$.
Hence, we have 
$K=\rnum(\gamma)$,
and 
$e^{2\pi\sqrt{-1}v/u}\in\rnum(\gamma)\subset\real$.
We obtain $\eta=\pm\gamma$,
and hence $P(T)=T^2-\gamma^2$.
Suppose $m\geq 1$.
We have
$e^{2\pi\sqrt{-1}v/(2^{m-1}u)}\in K$.
If 
$e^{2\pi\sqrt{-1}v/(2^{m-1}u)}\in\rnum$,
then
we have $e^{2\pi\sqrt{-1}v/(2^mu)} \in \{\pm 1,\pm\sqrt{-1}\}$.
Hence, $P=T^2\pm\gamma^2$.
Suppose 
$e^{2\pi\sqrt{-1}v/(2^{m-1}u)}\not\in\rnum$.
We have
$K=\rnum(e^{2\pi\sqrt{-1}v/(2^{m-1}u)})$.
Because $[K:\rnum]=2$,
we have $u=1$ or $3$.

In the case $u=1$,
we have $m=3$,
because
$e^{2\pi\sqrt{-1}v/2^{m-1}}\not\in \rnum$
with degree $2$.
If $(u,m)=(1,3)$,
we have $K=\rnum[\sqrt{-1}]$
and 
$\eta=\gamma e^{2\pi\sqrt{-1}v/8}\in\rnum[-1]$,
where $v=1,3,5,7$.
Then, we obtain
$\gamma=\sqrt{2}a$ for some $a\in\rnum_{>0}$.
The minimal polynomial of $\eta$ is
either one of the following
($\zeta_8:=e^{2\pi\sqrt{-1}/8}$):
\[
 \bigl(T-\sqrt{2}a \zeta_8
\bigr)
 \bigl(T-\sqrt{2}a\zeta_8^7
 \bigr)
=T^2-2aT+2a^2,
\quad\quad
 \bigl(T-\sqrt{2}a\zeta_8^{3}
\bigr)
 \bigl(T-\sqrt{2}a\zeta_8^5
\bigr)
=T^2+2aT+2a^2
\]

Similarly, if $u=3$,
we have $m=1$ or $2$.
If $(u,m)=(3,1)$,
we have $K=\rnum(e^{2\pi\sqrt{-1}/3})=\rnum(\sqrt{-3})$.
Because
$\gamma e^{2\pi\sqrt{-1}v/6}\in\rnum(\sqrt{-3})$,
we obtain $\gamma\in\rnum_{>0}$,
which contradicts with our assumption
$\gamma\not\in\rnum$.
If $(u,m)=(3,2)$,
we have
$\gamma e^{2\pi\sqrt{-1}v/12}\in\rnum(\sqrt{-3})$,
where $v=1,5,7,11$.
Hence, we obtain
$\gamma\in \sqrt{3}\rnum_{>0}$.
The minimal polynomial of $\eta$ is
either one of the following
($\zeta_{12}=e^{2\pi\sqrt{-1}/12}$):
\[
 \Bigl(
 T-\sqrt{3}a\zeta_{12}
 \Bigr)
 \Bigl(
 T-\sqrt{3}a\zeta_{12}^{11}
 \Bigr)
=T^2-3aT+3a^2,
\quad\quad
 \Bigl(
 T-\sqrt{3}a\zeta_{12}^{5}
 \Bigr)
 \Bigl(
 T-\sqrt{3}a\zeta_{12}^7
 \Bigr)
=T^2+3aT+3a^2
\]
Thus, we obtain Lemma \ref{lem;12.12.27.10}
\hfill\qed

\vspace{.1in}
If $r$ is an odd prime,
we have no exceptional cases
in the following sense.

\begin{prop}
\label{prop;12.12.28.1}
Suppose
(i) $r$ is an odd prime,
(ii) $\veca\in\rnum^r$,
(iii) $(\nbigv_r,\nabla^{(r)}_{\veca})$ is irreducible.
Then, $(\nbigv_r,\nabla^{(r)}_{\veca})$
has a $\seisuu$-structure,
if and only if
there exists a half-integer $\ell$
such that 
$\prod_{i=1}^r(T-e^{2\pi\sqrt{-1}(a_i+\ell/r)})\in\rnum[T]$.
\end{prop}
\pf
By Theorem \ref{thm;12.12.26.12},
$(\nbigv_r,\nabla^{(r)}_{\veca})$
has a $\seisuu$-structure,
if and only if
there exist
$\gamma\in\real_{>0}$ and a half-integer $\ell$
such that 
(i) $\gamma^r\in\rnum$,
(ii) $\prod_{i=1}^r(T-\gamma e^{2\pi\sqrt{-1}(a_i+\ell/r)})
 \in \rnum[T]$.
Suppose $\gamma\not\in\rnum$.
By using Lemma \ref{lem;12.12.27.10},
we obtain 
$\prod_{i=1}^r(T-\gamma e^{2\pi\sqrt{-1}(a_i+\ell/r)})
=T^r\pm\gamma^r$.
It implies that
$\prod_{i=1}^r(T- e^{2\pi\sqrt{-1}(a_i+\ell/r)})\in\rnum[T]$.
Hence, we obtain the proposition.
\hfill\qed

\vspace{.1in}

If $r=2$,
we obtain the following proposition from 
Lemma \ref{lem;12.12.27.10}.
We have the exceptional cases.
\begin{prop}
Suppose (i) $\veca\in\rnum^2$,
(ii) $(\nbigv_r,\nabla^{(2)})$ is irreducible.
Then, $(\nbigv_r,\nabla^{(2)})$ has a $\seisuu$-structure
if and only if there exists an integer $n$
such that one of the following holds:
\begin{itemize}
\item
 $\prod_{i=1}^2(T-e^{2\pi\sqrt{-1}(a_i+n/4)})\in\rnum[T]$.
\item
$\{a_1,a_2\}=
 \{-1/8+n/4,-7/8+n/4\}$
\item
$\{a_1,a_2\}=
 \{-1/12+n/4,\,-11/12+n/4\}
 \,\,\mbox{\rm or }
 \{-5/12+n/4,-7/12+n/4\}$.
\hfill\qed
\end{itemize}
\end{prop}

Let us compute Stokes matrices in a case.
The other cases can be computed similarly.
Let $(a_1,a_2)=(-1/8,-7/8)$.
We have $\gamma=\sqrt{2}$
and 
\[
 (T-\gamma e^{2\pi\sqrt{-1}a_1})
 (T-\gamma e^{2\pi\sqrt{-1}a_2})
=T^2-2T+2.
\]
We have
$b=-(a_1+a_2)-(r-1)/2=1/2$.
We take a formal frame $\vecp=(\phat_1,\phat_2)$
compatible with the Hukuhara-Levelt-Turrittin decomposition
such that
\[
(-1)^{\ast}\vecphat
 =\vecphat
 \frac{1}{\sqrt{2}}
 \left(
 \begin{array}{cc}
 0 & -2 \\ 1 & 0
 \end{array}
 \right)
\]
We put 
$\vecphat':=(\phat_2,\phat_1)$.

Let $S_1:=\bigl\{
 -\epsilon<\arg q<\pi+\epsilon
 \bigr\}$
and 
$S_2:=\bigl\{
 -\pi-\epsilon<\arg q<\epsilon
 \bigr\}$.
We put
$T_1:=\{\pi-\epsilon<\arg q<\pi+\epsilon\}$
and 
$T_2:=\{-\epsilon<\arg q<\epsilon\}$.
We have
$S_1\cap S_2=T_1\sqcup T_2$.
We take flat frames
$\vecptilde_{S_i}$,
which are lifts of $\vecphat$.
On $T_1$,
we have
\[
 \vecptilde_{S_1}
=\vecptilde_{S_2}A,
\quad
 A=\left(
 \begin{array}{cc}
 1 & 0 \\
 \alpha & 1
 \end{array}
 \right)
\]
We have 
\[
 (-1)^{\ast}\vecptilde_{S_1}
=\vecptilde_{S_2}
 \frac{1}{\sqrt{2}}
 \left(
 \begin{array}{cc}
 0 & -2 \\ 1 & 0
 \end{array}
 \right)
\]
Let $M$ be the monodromy matrix
given by
$\vecptilde_{S_1}=(-1)^{\ast}\vecptilde_{S_1}$.
Then, we have
\[
 \sqrt{2}M=
 \left(
 \begin{array}{cc}
 0 & 2 \\ -1 & 0
 \end{array}
 \right)
 \left(
 \begin{array}{cc}
 1 & 0 \\ \alpha & 1
 \end{array}
 \right)
=\left(
 \begin{array}{cc}
 2\alpha & 2 \\
 -1 & 0
 \end{array}
 \right)
\]
We obtain
$T^2-2T+2
=(T-2\alpha)T+2$
and hence $\alpha=1$.

Let $\vecptilde'_{S_i}$ be the lift of $\vecphat'$.
On $T_2$, we have 
\[
 \vecptilde'_{S_2}
=\vecptilde'_{S_1}A',
\quad
 A'=\left(
 \begin{array}{cc}
 1 & 0 \\
 \beta & 1
 \end{array}
 \right)
\]
We have 
\[
 (-1)^{\ast}\vecptilde'_{S_2}
=\vecptilde_{S_1}
 \frac{1}{\sqrt{2}}
 \left(
 \begin{array}{cc}
 0 & 1 \\ -2 & 0
 \end{array}
 \right)
\]
Let $M'$ be the monodromy matrix
given by
$\vecptilde'_{S_2}=(-1)^{\ast}\vecptilde'_{S_2}$.
Then, we have
\[
 \sqrt{2}M'=
 \left(
 \begin{array}{cc}
 0 & -1 \\ 2 & 0
 \end{array}
 \right)
 \left(
 \begin{array}{cc}
 1 & 0 \\ \beta & 1
 \end{array}
 \right)
=\left(
 \begin{array}{cc}
 -\beta & -1\\
 2 & 0
 \end{array}
 \right)
\]
We obtain
$T^2-2T+2
=(T+\beta)T+2$
and hence $\beta=-2$.

\section{Harmonic bundle with homogeneity}
\label{section;12.12.14.100}
\subsection{Preliminary}

\subsubsection{Holomorphic vector fields}

Let $X$ be a complex manifold.
The tangent bundle $TX$ is equipped with
a complex structure,
i.e., an automorphism $J$
with $J^2=-1$.
We have the canonical decomposition
$TX\otimes_{\real}\cnum
=T^{1,0}X\oplus T^{0,1}X$,
i.e.,
$T^{1,0}X$ and $T^{0,1}X$
denote the eigen spaces of $J$
corresponding to $\sqrt{-1}$
and $-\sqrt{-1}$.
The bundle $T^{1,0}X$ is naturally a holomorphic
vector bundle on $X$,
and $T^{0,1}X$ is naturally a holomorphic
vector bundle on the conjugate $X^{\dagger}$.
For the natural conjugation on $TX\otimes\cnum$,
we have $T^{0,1}X=\overline{T^{1,0}X}$.

A section $\gminiv$ of $TX$
has the decomposition
$\gminiv=\gminiv^{1,0}+\gminiv^{0,1}$,
and we have
$\overline{\gminiv^{1,0}}
=\gminiv^{0,1}$.
The section $\gminiv$ is called holomorphic,
if $\gminiv^{1,0}$ is a holomorphic section of $T^{1,0}X$,
which is equivalent to
that $\gminiv^{0,1}$ is a holomorphic section of
$T^{1,0}X^{\dagger}=T^{0,1}X$ on $X^{\dagger}$.
If we are given a local coordinate $(z_1,\ldots,z_n)$,
we have the expression
$\gminiv=\sum (f_i\del_{z_i}+\fbar_i\del_{\zbar_i})$.
It is holomorphic if and only if $f_i$ are holomorphic.

\begin{lem}
Let $L_{\gminiv}$ denote the Lie derivative
on $\Omega^{\bullet}_X$.
The following conditions are equivalent.
\begin{itemize}
\item
The vector field $\gminiv$ is holomorphic.
\item
$L_{\gminiv}\Omega_X^{p,q}\subset \Omega_X^{p,q}$
for any $p,q$.
\item
$[L_{\gminiv},\delbar_X]=0$.
\item
$[L_{\gminiv},\del_X]=0$.
\end{itemize}
\end{lem}
\pf
If $\gminiv=\sum (f_i\del_{z_i}+\fbar_i\del_{\zbar})$,
we have
$L_{\gminiv}(dz_i)=df_i$,
from which we can deduce the equivalence
of the first and second conditions.
Because $[L_{\gminiv},d_X]=0$,
the second implies the third and fourth.
If the third condition holds,
we have
$L_{\gminiv}(dz_i)=L_{\gminiv}\delbar_X z_i
=\delbar_X(L_{\gminiv}z_i)$.
Hence, we have
$L_{\gminiv}\Omega_X^{1,0}\subset\Omega_X^{1,0}$.
Similarly, we have
$L_{\gminiv}\Omega_X^{0,1}\subset\Omega_X^{0,1}$.
Then, we obtain the second.
Clearly, the third and the fourth are equivalent.
\hfill\qed

\subsubsection{An action on vector bundle}

Let $E$ be a $C^{\infty}$-vector bundle on $X$.
We identify it with the sheaf of $C^{\infty}$-sections.
Let $\gminiv$ be a vector field on $X$.
An action of $\gminiv$ on $E$
is a $\cnum$-linear map  of sheaves
$L^E_{\gminiv}:E\lrarr E$
such that 
$L_{\gminiv}^E(fs)
=\gminiv(f)\,s+fL_{\gminiv}^E(s)$.
It is standard that,
an action of $\gminiv$ on $E$
uniquely induces an action of $\gminiv$
on $\Omega^{\bullet}_X\otimes E$
satisfying
$L_{\gminiv}^E(\omega\otimes s)
=L_{\gminiv}(\omega)\otimes s
+\omega\otimes L_{\gminiv}^E(s)$.
A section $f$ of $E$ is called homogeneous 
of degree $\alpha\in\cnum$ with respect to $\gminiv$,
if $L^E_{\gminiv}f=\alpha f$.
An action of $\gminiv$ on $E_i$ $(i=1,2)$
naturally induces an action on
$E_1\oplus E_2$,
$E_1\otimes E_2$
and $\Hom(E_1,E_2)$.

\vspace{.1in}
Suppose that $E$ is a holomorphic vector bundle,
and that $\gminiv$ is holomorphic.
An action $L^E_{\gminiv}$ of $\gminiv$ on $E$ is 
called holomorphic,
if $[L^E_{\gminiv},\delbar_E]=0$.
The following lemma is clear.
\begin{lem}
Let $E^{\hol}$ denote the sheaf of holomorphic sections.
A holomorphic action of $\gminiv$ on $E$
induces a $\cnum$-linear map
$L_{\gminiv}^E:E^{\hol}\lrarr E^{\hol}$
such that
$L_{\gminiv}^E(f\,s)=
 L_{\gminiv}(f)\,s+f\,L_{\gminiv}^E(s)$
for local sections $f\in\nbigo_X$
and $s\in E^{\hol}$.
Conversely,
such a $\cnum$-linear map on $E^{\hol}$
induces a holomorphic action of $\gminiv$
on $E$.
\hfill\qed
\end{lem}

\subsection{Homogeneous harmonic bundle
and integrable variation of twistor structure}
\label{subsection;12.12.13.10}

Let $X$ be a complex manifold
with a holomorphic vector field $\gminiv$.
A harmonic bundle
$(E,\delbar_E,\theta,h)$ is called
homogeneous of degree $m\in\real$,
with respect to $\gminiv$,
if $E$ is equipped with a $\gminiv$-action $L^E_{\gminiv}$
such that
$[L^E_{\gminiv},\delbar_E]=0$,
$L^E_{\gminiv}(h)=0$
and
$[L^E_{\gminiv},\theta]=\sqrt{-1}m\theta$.
Note that we have
$[L^E_{\gminiv},\del_{E,h}]=0$
and 
$[L^E_{\gminiv},\theta^{\dagger}_h]
=-\sqrt{-1}m\theta^{\dagger}$.
In the following,
we impose that $m\neq 0$,
and we shall observe that
the associated variation of twistor structure
is naturally integrable,
if $m\neq 0$.

\begin{rem}
Suppose $[L_{\gminiv}^E,\theta]=\alpha\theta$
for some complex number $\alpha\not\in\sqrt{-1}\real$,
instead of
$[L_{\gminiv}^E,\theta]=\sqrt{-1}m\theta$.
Then, the local structure of the harmonic bundle
is easy.
See Proposition {\rm\ref{prop;12.12.20.11}} below.
Hence, we exclude the case.
\hfill\qed
\end{rem}

\subsubsection{The associated meromorphic flat connection}
\label{subsection;12.12.14.20}

We set $\nbigx:=\cnum_{\lambda}\times X$
and $\nbigxlambda:=\{\lambda\}\times X$.
Let $p_{\lambda}:\nbigx\lrarr X$ denote the projection.
We have the holomorphic vector bundle
$\nbige=p_{\lambda}^{-1}E$
with the holomorphic structure
$\delbar_E+\lambda\theta^{\dagger}+\delbar_{\lambda}$.
It is equipped with
the relative flat connections
\[
 \nabla^{\lambda}:
 \nbige\lrarr\nbige\otimes
 \Omega^1_{\nbigx/\cnum_{\lambda}}
\otimes\nbigo_{\nbigx}(\nbigx^{0})
\]
given by 
$\nabla^{\lambda}:=
 (\delbar_E+\lambda\theta^{\dagger})
 +(\del_E+\lambda^{-1}\theta)$.

Let $\gminivtilde:=
 \gminiv+
m\sqrt{-1}\bigl(
 \lambda\del_{\lambda}
-\lambdabar\del_{\lambdabar}
\bigr)$.
We have natural actions of 
$m\sqrt{-1}\bigl(
 \lambda\del_{\lambda}
-\lambdabar\del_{\lambdabar}
\bigr)$
and 
$\gminiv$ on $p^{-1}(E)$.
They give an action of $\gminivtilde$
on $p^{-1}(E)$,
denoted by
$L^{\nbige}_{\gminivtilde}$.
We have
$[L^{\nbige}_{\gminivtilde},\nabla^{\lambda}]=
[L^{\nbige}_{\gminivtilde},\delbar_{\nbige}]=0$.
We can check by an elementary computation
that 
$\nbige$ has a unique meromorphic flat connection
\[
 \nablatilde:
 \nbige\lrarr
 \nbige\otimes
 \Omega^{1}_{\nbigx}
 \bigl(\log\nbigx^{0}\bigr)
 \otimes\nbigo(\nbigx^{0}),
\]
such that
(i) $\nablatilde_v=\nabla^{\lambda}_v$
for any vector fields $v$ in the $X$-direction,
(ii) $\nablatilde_{\gminivtilde}=L^{\nbige}_{\gminivtilde}$.

\subsubsection{Conjugate}
\label{subsection;12.12.3.11}

Let $X^{\dagger}$ denote the conjugate of $X$.
We set $\nbigx^{\dagger}:=\cnum_{\mu}\times X^{\dagger}$.
The projection
$\nbigx^{\dagger}\lrarr X^{\dagger}$ 
is denoted by $p_{\mu}$.
We obtain a holomorphic vector bundle
$\nbige^{\dagger}:=
 \bigl(
 p_{\mu}^{-1}E,\,p_{\mu}^{\ast}\del_E+\mu\theta
 \bigr)$.
It is equipped with the relative flat connection
$\nabla^{\dagger\mu}=
 \delbar_E+\del_E+\mu\theta+\mu^{-1}\theta^{\dagger}$.

Note that $\gminiv$ naturally gives a holomorphic 
vector field on $X^{\dagger}$.
The harmonic bundle
$(E,\del_E,\theta^{\dagger},h)$
is homogeneous of degree $-m$
with respect to $\gminiv$.
By the same procedure,
we obtain a meromorphic connection
\[
 \nablatilde^{\dagger}:
 \nbige^{\dagger}
\lrarr
 \nbige^{\dagger}\otimes
 \Omega^{1}_{\nbigx^{\dagger}}
 \bigl(\log\nbigx^{0 \dagger}\bigr)
\otimes
 \nbigo(\nbigx^{0\dagger}).
\]

We identify
$\cnum_{\lambda}^{\ast}$
and $\cnum_{\mu}^{\ast}$
by $\lambda=\mu^{-1}$.
We have a natural identifications
$\nbige_{|\cnum_{\lambda}^{\ast}\times X}
=\nbige_{|\cnum_{\mu}^{\ast}\times X^{\dagger}}$.

\begin{lem}
\label{lem;12.12.3.1}
On $\cnum_{\lambda}^{\ast}\times X$,
we have $\nablatilde=\nablatilde^{\dagger}$.
\end{lem}
\pf
We have only to prove that
$\nablatilde_{\del_{\lambda}}
=\nablatilde^{\dagger}_{\del_{\lambda}}$.
We remark
$\sqrt{-1}m(\lambda\del_{\lambda}-\lambdabar\del_{\lambdabar})
=-\sqrt{-1}m(\mu\del_{\mu}-\mubar\del_{\mubar})$.
By the constructions,
we have the equalities
$\nablatilde_{V}=\nablatilde^{\dagger}_V$,
if $V$ is 
$\gminiv+\sqrt{-1}m
 (\lambda\del_{\lambda}-\lambdabar\del_{\lambdabar})$,
or if $V$ is a vector field in the $X$-direction.
We obtain
$\nablatilde_{\lambda\del_{\lambda}-\lambdabar\del_{\lambdabar}}
=\nablatilde^{\dagger}_{\lambda\del_{\lambda}-\lambdabar\del_{\lambdabar}}$.
Because the $d\lambdabar$-parts are the same,
we obtain
$\nablatilde_{\del_{\lambda}}
=\nablatilde^{\dagger}_{\del_{\lambda}}$.
\hfill\qed

\subsubsection{Polarized integrable variation of pure twistor structure}
\label{subsection;12.12.14.50}

Let $p:\proj^1\times X\lrarr X$
be the projection.
Due to Simpson \cite{s3},
we have the polarized variation of pure twistor structure
$(p^{-1}(E),\DD^{\sankaku})$ of weight $0$,
associated to $(E,\delbar_E,\theta,h)$.
(See \cite{mochi8} for a review.)
Lemma \ref{lem;12.12.3.1} implies that
the operators
$\nablatilde$ and $\nablatilde^{\dagger}$
gives a $\TTtilde E$-structure $\DDtilde^{\sankaku}$
of $p^{-1}(E)$.
(\cite{Hertling}. See \S2.1 of \cite{mochi8}
for a review of $\TTtilde E$-structure.)
Hence, we obtain 
an integrable variation of pure twistor structure
$(p^{-1}(E),\DDtilde^{\sankaku})$,
which is an enrichment of 
$(p^{-1}(E),\DD^{\sankaku})$.

\vspace{.1in}

Let $\sigma:\proj^1\lrarr\proj^1$
be given by $\sigma(\lambda)=(-\lambdabar)^{-1}$.
The induced map
$\proj^1\times X\lrarr\proj^1\times X$
is also denoted by $\sigma$.
A polarization $S_h$ of $(p^{-1}(E),\DD^{\sankaku})$
is a pairing of 
$(p^{-1}(E),\DD^{\sankaku})$
and 
$\sigma^{\ast}(p^{-1}(E),\DD^{\sankaku})$
taking values in 
the trivial line bundle
on $\proj^1\times X$,
satisfying some conditions.
In this case,
it is given by
$S_h(u,\sigma^{\ast}v):=h(u,\sigma^{\ast}v)$
for sections $u$ and $v$ of $p^{-1}(E)$
on $\nbigu$ and $\sigma(\nbigu)$,
respectively,
where $\nbigu$ is any open subset of $\proj^1\times X$.
\begin{lem}
The polarization is compatible with
the integrability,
i.e., for local sections $u,v$ of $p^{-1}(E)$,
we have
\begin{equation}
 \label{eq;12.12.13.2}
 \lambda\del_{\lambda}S_h(u,\sigma^{\ast}v)
=S_h\bigl(
 \nablatilde_{\lambda\del_{\lambda}}u,\sigma^{\ast}v
 \bigr)
+S_h\bigl(
 u,\sigma^{\ast}(\nablatilde_{\lambda\del_{\lambda}}v)
 \bigr).
\end{equation}
\end{lem}
\pf
We have 
$d_Xh(u,v)
=h\bigl(\nabla^{\lambda}u,\sigma^{\ast}v\bigr)
+h\bigl(u,\sigma^{\ast}(\nabla^{\lambda}v)\bigr)$.
Hence, for a vector field $V$ in the $X$-direction,
we have
\[
 Vh(u,\sigma^{\ast}v)
=h(\nablatilde_Vu,\sigma^{\ast}v)
+h\bigl(u,\sigma^{\ast}(\nablatilde_{\Vbar}v)\bigr).
\]
We have
$\delbar_{\lambda}h(u,\sigma^{\ast}v)d\lambdabar
=h\bigl(\delbar_{\lambda}u\,d\lambdabar,\,\sigma^{\ast}v\bigr)
+h\bigl(u,\sigma^{\ast}(\delbar_{\mu}v\,d\mubar)\bigr)$,
and hence
\[
 \lambdabar\delbar_{\lambda}
 h(u,\sigma^{\ast}v)
=h\bigl(\lambdabar\nablatilde_{\lambdabar}u,
 \sigma^{\ast}v\bigr)
+h\bigl(u,\sigma^{\ast}(\mubar\nablatilde_{\mubar}v)\bigr).
\]
Note $\sigma^{\ast}(\gminivtilde)=\gminivtilde$.
We have
$\gminivtilde h(u,\sigma^{\ast}v)
=h\bigl(L^{\nbige}_{\gminivtilde}u,\sigma^{\ast}v\bigr)
+h\bigl(u,\sigma^{\ast}(L^{\nbige}_{\gminivtilde}v)\bigr)$.
We deduce
\[
 \lambda\del_{\lambda}
 h(u,\sigma^{\ast}v)
=h(\lambda\nablatilde_{\lambda}u,\sigma^{\ast}v)
+h(u,\sigma^{\ast}(\mu\nablatilde_{\mu}v)),
\]
i.e.,
we obtain (\ref{eq;12.12.13.2})
\hfill\qed

\subsubsection{Equivalence}

Let $(V^{\sankaku},\DDtilde^{\sankaku})$
be an integrable variation of twistor structure on $X$.
If $\DDtilde^{\sankaku}_{\gminivtilde}(s)$ is $C^{\infty}$
for any $C^{\infty}$-section $s$ of $V^{\sankaku}$,
then $(V^{\sankaku},\DDtilde^{\sankaku})$ is called 
homogeneous of degree $m$ with respect to $\gminiv$.
If we describe $(V^{\sankaku},\DDtilde^{\sankaku})$
as the gluing of 
$TE$-structure
$(\nbigv,\nablatilde)$
and $\Ttilde E$-structure
$(\nbigv^{\dagger},\nablatilde^{\dagger})$,
the condition is equivalent to 
$\nablatilde_{\gminivtilde}\nbigv
\subset\nbigv$
and 
$\nablatilde^{\dagger}_{\gminivtilde}\nbigv^{\dagger}
\subset\nbigv^{\dagger}$.

In \S\ref{subsection;12.12.14.20}--\ref{subsection;12.12.14.50},
we explained that
a harmonic bundle with homogeneity of degree $m$
with respect to $\gminiv$
induces a polarized integrable variation of pure twistor structure
of weight $0$ which is homogeneous of degree $m$
with respect to $\gminiv$.
The construction gives an equivalence,
which is a variant of the fundamental result due to C. Simpson
\cite{s3}.
(See \cite{Hertling} and \cite{sabbah2} 
for the TERP case and the integrable case.)

\begin{prop}
Suppose that $(V^{\sankaku},\DDtilde^{\sankaku})$
is pure of weight $0$ and equipped with a polarization $S$.
Then, the corresponding harmonic bundle is 
naturally homogeneous of degree $m$ with respect to $\gminiv$.
The constructions are mutually converse.
\end{prop}
\pf
Let $p:\proj^1\times X=\gbigx\lrarr X$ denote the projection.
Let $\harmonicbundle$ be the harmonic bundle
corresponding to the polarized variation of 
pure twistor structure of weight $0$.
Let $\DD^{\sankaku}$ and
$\DDtilde^{\sankaku}_{\proj^1}$
denote the restriction of $\DDtilde^{\sankaku}$
to the $X$-direction and $\proj^1$-direction, respectively.
Let $d''_{\proj^1}$ denote the $(0,1)$-part of
$\DDtilde^{\sankaku}_{\proj^1}$.
We have
$V^{\sankaku}=p^{-1}(E)$
and $\DD^{\sankaku}=
 \delbar_E+\del_E
+\lambda\theta^{\dagger}+\lambda^{-1}\theta$.
Note that $[d''_{\proj^1},\DDtilde^{\sankaku}_{\gminivtilde}]=0$.
Hence, for any $C^{\infty}$-section $s$ of $E$,
we can define a $C^{\infty}$-section
$L^E_{\gminiv}s$ of $E$
by the condition
$\DDtilde^{\sankaku}_{\gminivtilde} p^{-1}(s)
=p^{-1}\bigl(
 L^E_{\gminiv}s
 \bigr)$.
It is easy to check that
$L^E_{\gminiv}$ gives an action of $\gminiv$ on $E$.

We have another description of
$L^E_{\gminiv}$.
Let $\gbigx^0:=\{0\}\times X\subset\gbigx$.
For any $C^{\infty}$-section $s$ of $E=V^{\sankaku}_{|\gbigx^0}$,
we take any $C^{\infty}$-section $\stilde$ of $V^{\sankaku}$
such that $\stilde_{|\gbigx^0}=s$.
Then, we have
$L^E_{\gminiv}(s)=
 \DDtilde^{\sankaku}_{\gminivtilde}(\stilde)_{|\gbigx^0}$.

Let $\DD^{\sankaku}=\DD^{1,0}+\DD^{0,1}$
denote the decomposition into the $(1,0)$-part
and $(0,1)$-part.
We have
$\bigl[
 \DDtilde^{\sankaku}_{\gminivtilde},\DD^{0,1}
\bigr]=0$
and
$\bigl[
 \DDtilde^{\sankaku}_{\gminivtilde},\lambda\DD^{1,0}
\bigr]=\sqrt{-1}m\lambda\DD^{1,0}$.
By using the above description,
we obtain
$[L^E_{\gminiv},\delbar_E]=0$
and
$[L^E_{\gminiv},\theta]=\sqrt{-1}m\theta$.

For $C^{\infty}$-sections $u$ and $v$ of $E$,
we have
$S(p^{-1}u,\sigma^{\ast}p^{-1}v)=h(u,v)$.
Because
\[
\gminivtilde
 S(u,\sigma^{\ast}v)=
 S\bigl(\DDtilde^{\sankaku}_{\gminivtilde}u,\sigma^{\ast}v\bigr)
+S\bigl(u,\sigma^{\ast}(\DDtilde^{\sankaku}_{\gminivtilde}v) \bigr),
\quad\quad
(\lambda\del_{\lambda}-\lambdabar\del_{\lambdabar})h(u,v)=0,
\]
we obtain 
$L^E_{\gminiv}h=0$.
It is clear that the constructions are mutually converse.
\hfill\qed

\subsubsection{Appendix}
\label{subsection;12.12.20.2}

Let $\harmonicbundle$ be a harmonic bundle on $X$.
Suppose that $E$ is equipped with an action of
a holomorphic vector field $\gminiv$
such that
$[L^E_{\gminiv},\delbar_E]=0$,
$L^E_{\gminiv}(h)=0$
and 
$[L^E_{\gminiv},\theta]=\alpha\theta$
for a complex number $\alpha\not\in\sqrt{-1}\real$.

\begin{lem}
\label{lem;12.12.20.10}
We have $[\delbar_E,\del_E]=0$
and $[\theta,\theta^{\dagger}]=0$.
\end{lem}
\pf
For each point of $P\in X$,
there exists $\epsilon>0$,
such that 
$(E,\delbar_E,e^{s\alpha}\theta,h)$
are harmonic bundles 
for any $-\epsilon<s<\epsilon$
on a neighbourhood $X_P$ of $P$.
We obtain
$[\delbar_E,\del_E]+e^{2s\Re(\alpha)}[\theta,\theta^{\dagger}]=0$.
Then, the claim of the lemma follows.
\hfill\qed

\begin{prop}
\label{prop;12.12.20.11}
Let $P$ be any point of $X$.
There exist a simply connected neighbourhood $X_P$ of $P$,
holomorphic functions $\gminia_i$ $(i=1,\ldots,\ell)$
and an orthogonal decomposition
$(E,\theta)_{|X_P}=\bigoplus (E_{i},d\gminia_i\id_{E_i})$.
\end{prop}
\pf
The Higgs bundle $(E,\theta)$ naturally induces
a coherent $\nbigo_{T^{\ast}X}$-module.
Its support with the reduced structure
is called the spectral variety of the Higgs bundle,
and denoted by $\Sigma_{E,\theta}$.
We say that a point $Q\in X$ is generic with respect to $(E,\theta)$,
if the projection $\Sigma_{E,\theta}\lrarr X$ is unramified over $Q$.

Suppose that $P$ is generic with respect to $(E,\theta)$.
We take a holomorphic coordinate
$(z_1,\ldots,z_n)$ around $P$.
We have the description
$\theta=\sum_{i=1}^n f_i\,dz_i$.
Let $f_i^{\dagger}$ denote the adjoint of $f_i$
with respect to $h$.
Because $[f_i,f_j^{\dagger}]=0$
for any $i$ and $j$ according to Lemma \ref{lem;12.12.20.10},
we have an orthogonal decomposition
$E_{|P}=\bigoplus_{\vecalpha\in\cnum^n}
 E_{\vecalpha}$
such that $f_{i|P}=\bigoplus_{\vecalpha} \alpha_i\id_{E_{\vecalpha}}$
and $f_{i|P}^{\dagger}=\bigoplus_{\vecalpha}\alphabar_i\id_{E_{\vecalpha}}$,
where $\alpha_i$ denotes the $i$-th component of $\vecalpha$.
On a small neighbourhood $X_P\ni P$,
we have a unique holomorphic decomposition
$E_{|X_P}=\bigoplus_{\vecalpha\in\cnum^n}
 E'_{\vecalpha}$
such that
$f_i(E'_{\vecalpha})\subset
 E'_{\vecalpha}$
and $E'_{\vecalpha|P}=E_{\vecalpha}$.

Because  $P$ is supposed to be generic with respect to $(E,\theta)$,
there exists a holomorphic function
$\beta_{\vecalpha,i}$ such that
$f_{i|E'_{\vecalpha}}=
 \beta_{\vecalpha,i}\id_{E'_{\vecalpha}}$.
Because we have $\del_E\theta=0$,
we obtain that 
$\sum_{i=1}^n \beta_{\vecalpha,i}dz_i$
is closed.
Hence, there exist holomorphic functions $\gminib(\vecalpha)$
such that
$f=\bigoplus d\gminib(\vecalpha)\id_{E'_{\vecalpha}}$.
Hence, we are done in the case that $P$ is generic.
Note that  the decomposition
$E_{|X_P}=\bigoplus E'_{\vecalpha}$ is orthogonal,
and hence flat with respect to the Chern connection
$\delbar_E+\del_E$.

Let us consider the case that $P$ is not necessarily generic.
Let $X_P$ be a simply connected small neighbourhood of $P$.
There exists a point $Q\in X_P$ which is generic
with respect to $(E,\theta)$.
On a small neighbourhood $X_Q\ni Q$ in $X_P$,
we have holomorphic functions $\gminia_i$
and 
an orthogonal decomposition
$(E,\theta)_{|X_Q}=\bigoplus (E_i,d\gminia_i\id_{E_i})$.
Because the decomposition is flat with respect to
$\delbar_E+\del_E$,
we can uniquely extend it to a flat decomposition
$E_{|X_P}=\bigoplus E_i'$,
which is orthogonal with respect to $h$.
It is easy to deduce that 
$\theta(E_i')\subset E_i'$,
and that 
$\theta_{|E_i'}$ has a unique eigenvalue
$d\gminia'_i$.
\hfill\qed

\subsection{Semi-infinite variation of Hodge structure}
\label{subsection;12.12.23.100}

\subsubsection{Real structure and 
semi-infinite variation of Hodge structure}
\label{subsection;12.12.3.12}

Let $X$ be a complex manifold.
Let $(E,\delbar_E,\theta,h)$ be a harmonic bundle
on $X$.
The dual vector bundle $E^{\lor}$
is naturally equipped with a Higgs field $\theta^{\lor}$
and a hermitian metric $h_{E^{\lor}}$,
so that $(E^{\lor},\delbar_{E^{\lor}},\theta^{\lor},h^{\lor})$
is a harmonic bundle.
Let $C$ be a non-degenerate symmetric pairing of
$(E,\delbar_E)$,
such that the induced map
$\Phi_C:E\lrarr E^{\lor}$
gives an isomorphism of harmonic bundles,
i.e.,
it is an isometric with respect to
$h$ and $h^{\lor}$,
and $\theta$ is symmetric with respect to $C$.
We have
\[
 \delbar C(u,v)=C(\delbar_Eu,v)+C(u,\delbar_Ev),
\quad\quad
 C(\theta u,v)-C(u,\theta v)=0.
\]
Because $\Phi_C$ is an isometry,
we have
\[
 \del C(u,v)=C(\del_Eu,v)+C(u,\del_Ev),
\quad\quad
 C(\theta^{\dagger}u,v)-C(u,\theta^{\dagger}v)=0.
\]
Hence, we obtain
\begin{equation}
 \label{eq;12.12.13.3}
 \delbar C(u,v)=
 C\bigl(
 (\delbar_E+\lambda \theta^{\dagger})u,\,
 v\bigr)
+C\bigl(
 u,\,(\delbar_E-\lambda\theta^{\dagger})v
 \bigr)
=0
\end{equation}
\begin{equation}
 \label{eq;12.12.13.4}
 \del C(u,v)=
 C\bigl(
 (\del_E+\lambda^{-1} \theta)u,\,
 v\bigr)
+C\bigl(
 u,\,(\del_E-\lambda^{-1}\theta)v
 \bigr)
=0
\end{equation}
Let $j:\nbigx\lrarr\nbigx$ be given by
$j(\lambda,x)=(-\lambda,x)$.
According to (\ref{eq;12.12.13.3}) 
and (\ref{eq;12.12.13.4}),
 $C$ naturally gives a $\nabla^{\lambda}$-flat pairing
$\nbige\otimes j_{\ast}\nbige
\lrarr\nbigo_{\nbigx^{\ast}}$.

\vspace{.1in}

Let $n$ be an integer.
Suppose moreover that $C$ is homogeneous of degree 
$\sqrt{-1}nm$ with respect to $\gminiv$,
i.e.,
$L^E_{\gminiv}C=\sqrt{-1}nmC$.
We have the induced pairing
$C:\nbige\otimes j^{\ast}\nbige\lrarr\nbigo_{\nbigx}$.
We set $P:=\lambda^{-n}C$.
The following lemma is clear by construction.
\begin{lem}
We have $L^{\nbige}_{\gminivtilde}P=0$,
and $\nablatilde P=0$.
\hfill\qed
\end{lem}

The pairing $P$ and 
the polarization $S_h$
induce a real structure
of the local system on $\cnum^{\ast}\times X$
corresponding to $(\nbige,\nablatilde)$.
(See \S\ref{subsection;12.10.3.100} below.)
Thus, it gives a variation of pure polarized TERP-structure.
Because we are given  an ``Euler field'' $\gminiv$
such that 
$\gminiv\nbige\subset\nbige$,
it is graded $\frac{\infty}{2}$-VHS
with real structure
which is pure and polarized.
(See \cite{Iritani}. The formulation here is slightly changed.)
If the base space is a punctured disc,
it is related with a Sabbah orbit or Hertling-Sevenheck orbit,
(``nilpotent orbit'' in the sense of Hertling and Sevenheck)
studied in \cite{Hertling-Sevenheck}.

\subsubsection{Appendix}
\label{subsection;12.10.3.100}

\paragraph{Symmetric pairing}

Let $V$ be a finite dimensional $\cnum$-vector space.
The natural pairing
$V^{\lor}\times V\lrarr \cnum$
is denoted by $\langle\cdot,\cdot\rangle$.
Let $C:V\times V\lrarr\cnum$ be a bilinear map.
It is equivalent to a linear map
$\Phi_C:V\lrarr V^{\lor}$
given by
$\langle \Phi_C(u),v \rangle=C(u,v)$.
It is symmetric if and only if
the dual $\Phi_C^{\lor}:V\lrarr V^{\lor}$
is equal to $\Phi_C$.
If $C$ is non-degenerate,
a pairing $C^{\lor}:V^{\lor}\times V^{\lor}\lrarr\cnum$
is induced by
$C^{\lor}(f,g)=C\bigl(
 \Phi_C^{-1}(f),\Phi_C^{-1}(g)
 \bigr)$.
Note that we have $\Phi_{C^{\lor}}=\Phi_C^{-1}$,
if $C$ is non-degenerate and symmetric.
Indeed,
\[
 \langle
 \Phi_{C^{\lor}}(f),\,g
 \rangle
=C^{\lor}(f,g)
=C\bigl(
 \Phi_C^{-1}(f),\Phi_C^{-1}(g)
 \bigr)
=\langle
 f,\Phi_C^{-1}(g)
 \rangle
=\langle
 (\Phi_C^{-1})^{\lor}(f),\,g
 \rangle.
\]

\paragraph{Conjugate}

Let $\Vbar$ be the vector space
such that $V=\Vbar$ as a $\seisuu$-module,
and that the scalar multiplication of $\alpha\in\cnum$
is given by $\alpha\bullet u=\alphabar\,u$.
Let $\overline{v}$ denote
the element of $\Vbar$ corresponding to
an element $v\in V$.
Then, $\alpha\,\vbar=\overline{(\alphabar v)}$.
We may regard $\overline{(\cdot)}$
as a $\cnum$-linear isomorphism
$V\simeq \Vbar$.

We have a natural anti-$\cnum$-linear map
$\Psi:\Hom_{\cnum}(V,\cnum)
\lrarr
\Hom_{\cnum}(\Vbar,\cnum)$
given by 
$\Psi(f)(\vbar)=\overline{f(v)}$.
It naturally gives a $\cnum$-isomorphism
$\overline{\Hom_{\cnum}(V,\cnum)}
\simeq
 \Hom_{\cnum}(\Vbar,\cnum)$,
expressed as
$\overline{(V^{\lor})}\lrarr \Vbar^{\lor}$.

We have the natural pairing
$\Vbar^{\lor}\times \Vbar\lrarr\cnum$.
Under the identifications,
we have
$\langle \fbar,\vbar\rangle
=\overline{\langle f,v\rangle}$.

\paragraph{Hermitian pairing}

Let $h:V\times \Vbar\lrarr\cnum$ 
be a bilinear map,
which is equivalent to
a linear map
$\Phi_h:V\lrarr \Vbar^{\lor}$
given by
$\langle \Phi_h(u),\vbar\rangle
=h(u,v)$.
It is hermitian,
if and only if
$(\overline{\Phi_h})^{\lor}
=\Phi_h$.
Indeed, we have
\[
h(u,v)
=\langle\Phi_h(u),\vbar\rangle
=\langle u,\Phi_h^{\lor}(\vbar)\rangle
=\overline{
 \langle \ubar,\Phibar_h^{\lor}(v)
 \rangle}.
\]
Because 
$\langle \ubar,\Phi_h(v)\rangle
=h(v,u)$,
we are done.

If $h$ is non-degenerate and hermitian,
we have the induced hermitian pairing
$h^{\lor}$ and $\overline{h}^{\lor}$
on $V^{\lor}$ and $\Vbar^{\lor}$,
respectively.
We have
$\Phi_{\overline{h}^{\lor}}=\Phi_h^{-1}$.
Indeed,
\[
 \langle
 \Phi_{\overline{h}^{\lor}}(\fbar),g
 \rangle
=\overline{h}^{\lor}(\fbar,\gbar)
=h\bigl(
 \Phi_h^{-1}(\fbar),\,\Phi_h^{-1}(\gbar)
 \bigr)
=\langle
 \fbar,\overline{\Phi_h^{-1}(\gbar)}
 \rangle
=\langle
 \fbar,\,\Phibar^{-1}_h(g)
 \rangle
=\langle
 \Phi_h^{-1}(\fbar),g
 \rangle
\]
We also have
$\Phibar_{h^{\lor}}=
 \Phi_{\overline{h}^{\lor}}$.

\paragraph{Compatibility}

Let $C$ be a non-degenerate symmetric pairing of $V$.
Let $h$ be a non-degenerate hermitian pairing of $V$.
We say that $C$ and $h$ are compatible,
if 
$h^{\lor}\bigl(\Phi_C(u),\Phi_C(v)\bigr)
=h(u,v)$,
i.e.,
$\Phi_C$ gives an isometry
$(V,h)\simeq (V^{\lor},h^{\lor})$.
In that case,
the following diagram is commutative:
\[
 \begin{CD}
 V @>{\Phi_C}>> V^{\lor} \\
 @V{\Phi_h}VV @VV{\Phi_{h^{\lor}}}V \\
 \Vbar^{\lor} @>{\Phibar_{C^{\lor}}}>> \Vbar
 \end{CD}
\]
By the compatibility,
we have
$h^{\lor}(\Phi_C(u),f)=h(u,\Phi_{C^{\lor}}(f))$.
The left hand side is equal to
$\langle
 \Phi_{h^{\lor}}\Phi_C(u),\,\fbar
 \rangle$,
and the right hand side is equal to
$\langle
 \Phi_h(u),\,\overline{\Phi_{C^{\lor}}(f)}
 \rangle
=\langle
 \Phibar_{C^{\lor}}\Phi_h(u),\fbar
 \rangle$.
Hence, we have
$\Phibar_{C^{\lor}}\circ\Phi_h
=\Phi_{h^{\lor}}\circ\Phi_C=:\kappa$,
which is a linear map
$V\lrarr \Vbar$.
We have
$\kappabar\circ\kappa=\id_V$.
Indeed,
$\kappabar\circ\kappa
=\Phibar_{h^{\lor}}\circ\Phibar_C
 \circ\Phibar_{C^{\lor}}
 \circ\Phi_h
=\Phibar_{h^{\lor}}\circ\Phi_h=1$.
Hence, $\kappa$ gives a real structure
of the vector space $V$.
In this case,
we have the relation
$h(u,v)=C(u,\overline{\kappa(v)})$.

\paragraph{Induced vector bundle
 with real structure on $\cnum^{\ast}$}

Let $V$ be a $\cnum$-vector space
with a non-degenerate hermitian pairing $h$
and a non-degenerate symmetric pairing $C$.
Assume that they are compatible.

Let $\sigma$, $\gamma$ and $j$
be the involutions of $\cnum_{\lambda}^{\ast}$,
given by
$\sigma(\lambda)=-\lambdabar^{-1}$,
$\gamma(\lambda)=\lambdabar^{-1}$
and $j(\lambda)=-\lambda$,
respectively.

We set $\nbigv:=V\otimes\nbigo_{\cnum^{\ast}_{\lambda}}$.
Let $S_h$ be the pairing
$\nbigv\otimes\sigma^{\ast}\nbigv\lrarr\nbigo_{\cnum^{\ast}}$
given by
$S_h(u,\sigma^{\ast}v)=h(u,\sigma^{\ast}v)$.
Let $S_C$ be the pairing
$\nbigv\otimes j^{\ast}\nbigv\lrarr\nbigo_{\cnum^{\ast}}$
given by
$S_C(u,j^{\ast}v)=C(u,j^{\ast}v)$.
The pairing $S_h$ satisfies the symmetry
$\overline{\sigma^{\ast}S_h(\sigma^{\ast}u,v)}
=S_h(v,\sigma^{\ast}u)$.
The pairing $S_C$ satisfies the symmetry
$j^{\ast}S_C(j^{\ast}u,v)=S_C(v,j^{\ast}u)$.
Let $\Phi_{S_h}:\nbigv\lrarr\sigma^{\ast}\nbigv^{\lor}$
and $\Phi_{S_C}:\nbigv\lrarr j^{\ast}\nbigv^{\lor}$
denote the induced morphisms.
We have the following commutativity by the compatibility
of $h$ and $C$:
\[
\begin{CD}
 \nbigv@>{\Phi_{S_C}}>> j^{\ast}\nbigv^{\lor}\\
 @V{\Phi_{S_h}}VV @VV{j^{\ast}\Phi_{S_h^{\lor}}}V \\
 \sigma^{\ast}\nbigv^{\lor}
 @>{\sigma^{\ast}\Phi_{S_C^{\lor}}}>>
 \gamma^{\ast}\nbigv
\end{CD}
\]
We set $\kappa:=j^{\ast}\Phi_{S_h^{\lor}}\circ\Phi_{S_C}$
which is a morphism
$\nbigv\lrarr\gamma^{\ast}\nbigv$
such that
$\gamma^{\ast}\kappa\circ\kappa=\id$.

If we set $P=\lambda^{-n}S_C$,
then we have the following $(-1)^n$-symmetry:
\[
 j^{\ast}(P)(j^{\ast}u,v)
=(-1)^nP(v,j^{\ast}v)
\]
We have $P^{\lor}=(-1)^n\lambda^n S_C^{\lor}$.
We have $\Phi_P=\lambda^{-n}\Phi_{S_C}$
and $\sigma^{\ast}\Phi_{P^{\lor}}=\sigma^{\ast}
 \bigl((-1)^{n}\lambda^n\Phi_{S_C^{\lor}}\bigr)
=\lambda^{-n}\sigma^{\ast}\Phi_{S_C^{\lor}}$.
Hence, the following diagram is commutative:
\[
\begin{CD}
 \nbigv@>{\Phi_{P}}>> j^{\ast}\nbigv^{\lor}\\
 @V{\Phi_{S_h}}VV @VV{j^{\ast}\Phi_{S_h^{\lor}}}V \\
 \sigma^{\ast}\nbigv^{\lor}
 @>{\sigma^{\ast}\Phi_{P^{\lor}}}>>
 \gamma^{\ast}\nbigv
\end{CD}
\]
We set $\kappa_n:=j^{\ast}\Phi_{S_h^{\lor}}\circ\Phi_{P}
=\lambda^{-n}\kappa$.
We have $\gamma^{\ast}\kappa_n\circ\kappa_n=\id$.
Hence, it also gives a real structure.

\vspace{.1in}
Assume that $\nbigv$ is equipped with a connection $\nabla$
such that
the pairings $S_h$ and $P$ are $\nabla$-flat.
Then, $\kappa$ is also flat.
As in \cite{Hertling},
it determines a real structure of the corresponding local system.
See \S2 of \cite{mochi8} for a review.

\subsection{The prolongation of the action}
\label{subsection;12.12.14.31}

Let $X$ be a complex manifold
with a simply normal crossing hypersurface $D$.
Let $D=\bigcup_{i\in\Lambda}D_i$
denote the irreducible decomposition.
Let $\gminiv$ be a holomorphic vector field on $X\setminus D$
such that $\gminiv^{1,0}\in\Theta_X(\log D)$.
Let $\harmonicbundle$ be a good wild harmonic bundle
on $X\setminus D$
which is homogeneous of degree $m$
with respect to $\gminiv$.

We set $\nbigx:=\cnum_{\lambda}\times X$
and $\nbigd:=\cnum_{\lambda}\times D$.
As explained in \S\ref{subsection;12.12.13.10},
we have the holomorphic vector bundle
$\nbige$ on $\nbigx$
with a meromorphic flat connection $\nablatilde$.
According to \cite{mochi2} and \cite{Mochizuki-wild}
(see \cite{mochi8} for a review),
$\nbige$ is naturally extended to a filtered vector bundle
$\nbigq_{\ast}\nbige
=\bigl(
 \nbigq_{\veca}\nbige\,\big|\,
 \veca\in\real^{\Lambda}
 \bigr)$ 
on $(\nbigx,\nbigd)$.
Let $\nbigq\nbige$ denote
the locally free $\nbigo_{\nbigx}(\ast \nbigd)$-module
$\bigcup_{\veca} \nbigq_{\veca}\nbige$.
According to \cite{Hertling-Sevenheck2} in the tame case
and \cite{mochi8} in the wild case,
$\nablatilde$ gives a meromorphic flat connection of
$\nbigq\nbige$:
\[
 \nablatilde:
 \nbigq\nbige\lrarr
 \nbigq\nbige\otimes
 \Omega^1_{\nbigx}(\log \nbigx^0)
 \otimes\nbigo_{\nbigx}(\nbigx^0).
\]
Here, $\nbigx^0:=\{0\}\times X$.

We have a $\cnum$-endomorphism 
$L^{\nbigq\nbige}_{\gminivtilde}$ of
$\nbigq\nbige$ given by 
$L^{\nbigq\nbige}_{\gminivtilde}:=\nablatilde_{\gminivtilde}$.
It is an extension of $L^{\nbige}_{\gminivtilde}$ on $\nbige$.
It satisfies
$L^{\nbigq\nbige}_{\gminivtilde}(fs)
=\gminivtilde(f)\,s+f\,L^{\nbigq\nbige}_{\gminivtilde}(s)$
for local sections
$f\in\nbigo_{\nbigx}(\ast\nbigd)$
and $s\in \nbigq\nbige$.

\begin{prop}
For any $\veca\in\real^{\Lambda}$,
we have 
$L^{\nbigq\nbige}_{\gminivtilde}
\bigl(
\nbigq_{\veca}\nbige
\bigr)
\subset
\nbigq_{\veca}\nbige$.
\end{prop}
\pf
We have only to consider the case
$X=\bigl\{(z_1,\ldots,z_n)\,\big|\,|z_i|<1\bigr\}$
and 
$D=\bigcup_{i=1}^{\ell} D_i$,
where $D_i=\{z_i=0\}$.
Let us consider the case that
the harmonic bundle $\harmonicbundle$
is unramifiedly good wild.
Let $O$ denote the origin $(0,\ldots,0)$.
We have a good set
$\nbigi\subset \nbigo_X(\ast D)/\nbigo_X$
and a formal decomposition
\[
 \bigl(\nbigq_{\ast}\nbige,\nablatilde\bigr)
 _{|\widehat{\cnum_{\lambda}\times\{O\}}}
=\bigoplus_{\gminia\in\nbigi}
 \bigl(
 \nbigq_{\ast}\nbigehat_{\gminia},\nablahat_{\gminia}
 \bigr)
\]
such that 
$\nablahat^{\reg}_{\gminia}:=
 \nablahat_{\gminia}
-d(\gminiatilde/\lambda)\id_{\nbigq\nbigehat_{\gminia}}$
are logarithmic along $\nbigd$,
where we choose lifts $\gminiatilde\in\nbigo_X(\ast D)$ 
of $\gminia$.
Around any point
$(\lambda_0,O)$ with $\lambda_0\neq 0$,
for any $\veca\in\real^{\ell}$,
we have 
$\nablahat^{\reg}_{\gminia}
 \nbigq_{\veca}\nbigehat_{\gminia}
\subset
 \nbigq_{\veca}\nbigehat_{\gminia}
 \otimes\Omega^1_{\nbigx}(\log\nbigd)$,
and hence
$\nablahat^{\reg}_{\gminia,\gminivtilde}
 \bigl(
 \nbigq_{\veca}\nbigehat_{\gminia}
 \bigr)
\subset
 \nbigq_{\veca}\nbigehat_{\gminia}$.
Because $[L^{E}_{\gminiv},\theta]=\sqrt{-1}m\theta$,
we obtain that
$\gminiv(\gminiatilde)-\sqrt{-1}m\gminiatilde$
is holomorphic.
We obtain that
$d\bigl(\gminiatilde/\lambda)(\gminivtilde\bigr)$
is holomorphic around $(\lambda_0,O)$.
Hence, we obtain that
$\nablahat_{\gminia,\gminivtilde}
 \bigl(
 \nbigq_{\veca}\nbigehat_{\gminia}
 \bigr)
\subset
 \nbigq_{\veca}\nbigehat_{\gminia}$
around $(\lambda_0,O)$.
We can conclude that
$L^{\nbigq\nbige}_{\gminivtilde}
 \bigl(\nbigq_{\veca}\nbige\bigr)
\subset
 \nbigq_{\veca}\nbige$
outside $\nbigd^0$.
By using the Hartogs property,
we obtain
$L^{\nbigq\nbige}_{\gminivtilde}
 \bigl(
 \nbigq_{\veca}\nbige\bigr) 
\subset
 \nbigq_{\veca}\nbige$
on $\nbigx$.

Let us consider the case that
$\harmonicbundle$ is not necessarily unramified.
We set $X'=\{(\zeta_1,\ldots,\zeta_n)\,\big|\,|\zeta_i|<1\}$
and $D'=\bigcup_{i=1}^{\ell}\{\zeta_i=0\}$.
Let $\varphi:(X',D')\lrarr (X,D)$ be the ramified covering
given by
$\varphi(\zeta_1,\ldots,\zeta_n)
=(\zeta_1^e,\ldots,\zeta_{\ell}^e,\zeta_{\ell+1},\ldots,\zeta_{n})$
such that 
$(E',\delbar_{E'},\theta',h'):=\varphi^{-1}\harmonicbundle$
is unramified.
Let $\gminiv'$ be the holomorphic vector field
of $X'\setminus D'$ which is the lift of $\gminiv$.
The harmonic bundle
$(E',\delbar_{E'},\theta',h')$
is homogeneous of degree $\sqrt{-1}m$
with respect to $\gminiv'$.
We have
$L^{\nbigq\nbige'}_{\gminivtilde'}
 \bigl(
 \nbigq_{\veca'}\nbige'
 \bigr)
\subset
\nbigq_{\veca'}\nbige'$
for any $\veca'\in\real^{\ell}$.
Note that a section $f$ of $\nbigq\nbige$
is contained in
$\nbigq_{\veca}\nbige$
if and only if
$\varphi^{\ast}f\in\nbigq_{e\veca}\nbige'$.
Then, we obtain
$L^{\nbigq\nbige}_{\gminivtilde}
 \bigl(
 \nbigq_{\veca}\nbige
 \bigr)
\subset\nbigq_{\veca}\nbige$.
\hfill\qed

\subsection{tr-TLE structure}

\subsubsection{Limit mixed twistor structure}

Let $X:=\bigl\{(z_1,\ldots,z_n)\,\big|\,|z_i|<1\bigr\}$
and $D:=\bigcup_{i=1}^{\ell}\{z_i=0\}$.
Let $\gminiv$ be a holomorphic vector field of $X\setminus D$
such that $\gminiv^{1,0}\in\Theta_X(\log D)$.
Let $\harmonicbundle$ be a {\em tame} harmonic bundle
on $X\setminus D$
which is homogeneous of degree $\sqrt{-1}m$
with respect to $\gminiv$.

For any $\veca\in\real^{\ell}$,
we have the limit mixed twistor structure
$\bigl(S^{\can}_{\veca}(E),\vecN\bigr)$.
(See \cite{mochi2}. See \cite{mochi8} for a review.)
Let us briefly recall the construction.
We have the following coherent $\nbigo$-module
on $\cnum_{\lambda}\times \bigcap_{i=1}^{\ell}\{z_i=0\}$:
\[
 \Gr^{\nbigq}_{\veca}(\nbige)
:=\frac{\nbigq_{\veca}(\nbige)}{\sum_{\vecb\lneq\veca}\nbigq_{\vecb}\nbige}
\]
Here, $\vecb\lneq\veca$ for $\vecb=(b_i),\veca=(a_i)\in\real^{\ell}$
means $b_i\leq a_i$ for any $i$ and $\vecb\neq\veca$.
According to \cite{mochi2},
$\Gr^{\nbigq}_{\veca}(\nbige)$ is locally free.
Hence, we have the vector bundle
$\nbigg_{\veca}(E):=
 \Gr^{\nbigq}_{\veca}(\nbige)
 _{|\cnum_{\lambda}\times\{O\}}$
on $\cnum_{\lambda}$.
They are equipped with the morphisms
$\nbign_i:
 \nbigg_{\veca}(E)
\lrarr
 \lambda^{-1}
 \nbigg_{\veca}(E)$
obtained as the nilpotent part of the residues
$\Res_{z_i}(\nablatilde)$.
We also obtain the vector bundle
$\nbigg^{\dagger}_{-\veca}(E)$
on $\cnum_{\mu}$
with the morphisms
$\nbign_i^{\dagger}:
\nbigg^{\dagger}_{-\veca}(E)
\lrarr
\mu^{-1}
\nbigg^{\dagger}_{-\veca}(E)$
from $(\nbigq_{\ast}\nbige^{\dagger},\nablatilde^{\dagger})$.
We have an isomorphism
$\Phi_{\veca}:
 \nbigg_{\veca}(E)_{|\cnum_{\lambda}^{\ast}}
\simeq
 \nbigg^{\dagger}_{-\veca}(E)_{|\cnum_{\mu}^{\ast}}$,
such that
$\Phi_{\veca}\circ\nbign_i
=-\nbign_i^{\dagger}\circ\Phi_{\veca}$.
(See \cite{mochi8} for a review of the construction.)
Thus, we obtain a vector bundle
$S^{\can}_{\veca}(E)$ on $\proj^1$
with morphisms
$\nbign_i:S^{\can}_{\veca}(E)
\lrarr
 S^{\can}_{\veca}(E)\otimes\nbigo_{\proj^1}(2)$.
Let $W$ denote the weight filtration of
$\sum_{i=1}^{\ell} \nbign_i$.
According to \cite{mochi2},
$(S^{\can}_{\veca}(E),W)$ is a mixed twistor structure,
i.e.,
$\Gr^W_nS^{\can}_{\veca}(E)$
are isomorphic to a direct sum of $\nbigo_{\proj^1}(n)$.

Because the harmonic bundle is homogeneous,
the associated polarized variation of
pure twistor structure is integrable
(\S\ref{subsection;12.12.13.10}).
As explained in \cite{Hertling-Sevenheck2}
(see \cite{mochi8} for a review),
$S^{\can}_{\veca}(E)$
is equipped with the induced meromorphic connection
$\lefttop{\veca}\nabla$ in the case of
polarized integrable variation of twistor structure.
Let $\gminiu:=
 \sqrt{-1}(\lambda\del_{\lambda}-\lambdabar\del_{\lambdabar})$.
We also have an action
$L^{S^{\can}}_{m\gminiu}$
on $S^{\can}_{\veca}(E)$
induced by $L^{\nbigq\nbige}_{\gminivtilde}$
(\S\ref{subsection;12.12.14.31}).
It is equal to
$\lefttop{\veca}\nabla_{m\gminiu}+\sum f_i(a_i+\nbign_i)$,
where $\gminiv_{|O}=\sum f_i\,(z_i\del_{z_i})_{|O}$.
The morphisms $\nbign_i$ are compatible with
$\lefttop{\veca}\nabla$ and $L^{S^{\can}}_{m\gminiu}$.

\subsubsection{Mixed Hodge structure}
\label{subsection;12.12.13.20}

Note that $\gminiu$ is the fundamental vector field
of the $S^1$-action on $\proj^1$
given by $(t,\lambda)\longmapsto t\lambda$.
For our purpose in \S\ref{section;12.12.14.200},
it is reasonable to suppose the following.
\begin{description}
\item[Assumption]
For any $\veca\in\real^{\ell}$,
the $\gminiu$-action on $S^{\can}_{\veca}(E)$
is lifted to an $S^1$-action.
\end{description}
It is naturally extended to a $\cnum^{\ast}$-action.
(See Lemma \ref{lem;12.10.3.20} below.)
The connection $\lefttop{\veca}\nabla$
and the morphisms $\nbign_i$
are $\cnum^{\ast}$-invariant.
In particular, the weight filtration is also
$\cnum^{\ast}$-equivariant.
Hence,
$(S^{\can}_{\veca}(E),W)$
comes from a complex mixed Hodge structure
$(H,F,G,W)$
with morphisms $N_i$,
by the Rees construction.
(\cite{s3} . See also \cite{mochi2}.)

Recall that, by using the canonical decomposition
of mixed Hodge structure by Deligne,
we obtain a decreasing filtration $F^{\op}$
such that
(i) it is $0$-opposite to $F$,
(ii) $N_i (F^{\op})^j\subset (F^{\op})^{j+1}$
for any $i=1,\ldots,\ell$ and $j\in\seisuu$.
Let us briefly review it.
We have the Deligne decomposition
$H=\bigoplus_{p,q\in\seisuu} I^{p,q}$
which satisfies
$F^p=\bigoplus_{r\geq p} I^{r,s}$
and $W_k=\bigoplus_{r+s\leq k}I^{r,s}$.
Moreover, 
for any morphism of mixed Hodge structures
$f:(H,F,G,W)\lrarr (H',F,G,W)$,
we have 
$f\bigl(I^{r,s}(H)\bigr)\subset I^{r,s}(H')$.
(See \cite{Peters-Steenbrink}, for example.)
We set
$(F^{\op})^q:=\bigoplus_{q+r\leq 0}I^{r,s}$.
Then, the filtrations $F$ and $F^{\op}$
are $0$-opposed,
i.e., $F^p\cap (F^{\op})^q=0$ for $p+q>0$
and $H=\bigoplus_{p+q=0} F^p\cap F^q$.
Because $N_i$ $(i=1,\ldots,\ell)$ give morphisms
 from $(H,F,G,W)$ to its $(-1)$-Tate twist,
we obtain that
$N_i\bigl((F^{\op})^q\bigr)\subset (F^{\op})^{q+1}$
by using the functoriality.

Applying the Rees construction
to two filtrations $F$ and $F^{\op}$,
we obtain a vector bundle
$\nbigr_{\veca}$ on $\proj^1$,
which is pure twistor structure of weight $0$.
It satisfies
$\nbigr_{\veca|\cnum_{\lambda}}
=S^{\can}_{\veca|\cnum_{\lambda}}$.
We have induced morphisms
\[
\nbign_i:
 \nbigr_{\veca}
\lrarr
 \nbigr_{\veca}\otimes\nbigo_{\proj^1}(\{0\}-\{\infty\}),
\quad
 (i=1,\ldots,\ell).
\]
The bundle $\nbigr_{\veca}$ is equipped with a $\cnum^{\ast}$-action,
which induces an action $L^{\nbigr_{\veca}}_{m\gminiu}$
of $m\gminiu$.
Because we have
$\lefttop{\veca}\nabla_{m\gminiu}
=L^{\nbigr_{\veca}}_{m\gminiu}-\sum f_i(a_i+\nbign_i)$
on $\cnum_{\lambda}$,
we have the induced connections
\[
\lefttop{\veca}\nabla:
 \nbigr_{\veca}
\lrarr
 \nbigr_{\veca}
 \otimes\Omega^1_{\proj^1}(2\{0\}+\{\infty\}).
\]
The residue of $\lefttop{\veca}\nabla$ at $\infty$ is semisimple.

\begin{rem}
\label{rem;13.1.2.1}
Hertling and Sevenheck {\rm\cite{Hertling-Sevenheck}}
obtained a mixed Hodge structure 
from their Sabbah orbit based on the above theorem
that $(S^{\can}_{\veca}(E),\vecN)$ are mixed twistor structure
in {\rm\cite{mochi2}}.
It seems related with the above mixed Hodge structure
in the case that $n=\ell=m=1$ and
$\gminiv=\sqrt{-1}(z\del_z-\zbar\del_{\zbar})$.
But, the construction in this paper is slightly different from theirs,
at least apparently.
Namely,
we construct a mixed Hodge structure
on the $L_{\gminiu}$-equivariant sections,
and they constructed a mixed Hodge structure
on the space of multi-valued sections
which are flat with respect to $\lefttop{\veca}\nabla$.
Their construction seems to have its origin
in Singularity theory.
Our construction might be natural
from the view point of ``Simpson's Meta theorem''
{\rm\cite{s3}}.
We postpone to clarify the precise relation.
\hfill\qed
\end{rem}

\subsubsection{tr-TLE structure}
\label{subsection;12.12.14.40}

We set $\gbigx:=\proj^1\times X$
and $\gbigd:=\proj^1\times D$.
We set $\gbigx^{\infty}:=\{\infty\}\times X$.
From $(\nbigq_0\nbige,\nablatilde)$,
we obtain an $\nbigo_{\gbigx}
\bigl(\ast \gbigx^{\infty}\bigr)$-module
$\nbigqtilde_0\nbige$
such that $\nablatilde$ is regular singular
along $\gbigx^{\infty}$.
We continue to suppose the assumption
in \S\ref{subsection;12.12.13.20}.
We regard $\proj^1$
as $\cnum_{\lambda}\cup\cnum_{\mu}$
by $\lambda\mu=1$.
We have the natural isomorphism
$\nbigqtilde_{0}\nbige_{|\cnum_{\mu}^{\ast}\times \{O\}}
\simeq\bigoplus_{\veca\in\openclosed{-1}{0}^{\ell}}
 \nbigr_{\veca|\cnum_{\mu}^{\ast}}$
which is the simultaneous
generalized eigen decomposition
of $\Res_{z_i}(\nablatilde)$ $(i=1,\ldots,\ell)$.
We naturally regard
 $\bigoplus_{\veca\in\openclosed{-1}{0}^{\ell}}
 \nbigr_{\veca|\cnum_{\mu}}$
as an $\nbigo_{\cnum_{\mu}}$-submodule
$\nbigqtilde_0\nbige_{|\cnum_{\mu}\times\{O\}}$.
It determines
a locally free $\nbigo_{\proj^1}$-submodule
$\nbigr\subset
 \nbigqtilde_0\nbige_{|\proj^1\times\{O\}}$.

\begin{lem}
\label{lem;12.12.20.20}
$\nbigr$ is a pure twistor structure of weight $0$.
\end{lem}
\pf
For any $-1<b\leq 0$,
let $\vecb_i\in\real^{\ell}$ denote the element
whose $j$-th entry is $b$ if $j=i$, or $0$ otherwise.
We have the naturally defined map
$\nbigq_{\vecb_i}\nbige_{|\cnum_{\lambda}\times\{O\}}
\lrarr
 \nbigq_{0}\nbige_{|\cnum_{\lambda}\times \{O\}}$.
The image is denoted by
$\lefttop{i}F_b(\nbigq_0\nbige_{|\cnum_{\lambda}\times\{O\}})$.
Thus, we obtain the filtrations
$\lefttop{i}F$ $(i=1,\ldots,\ell)$.
We set $\lefttop{\ellsitabar}F_{\veca}:=
 \bigcap_{i=1}^{\ell}\bigl(\lefttop{i}F_{a_i}\bigr)$.
According to \cite{mochi2},
we have an isomorphism
\[
 \nbigg_{\veca}(E)
\simeq
 \frac{\lefttop{\ellsitabar}F_{\veca}
 (\nbigq_0\nbige_{|\cnum_{\lambda}\times\{O\}})}
 {\sum_{\vecb\lneq\veca}
 \lefttop{\ellsitabar}F_{\vecb}(\nbigq_0\nbige_{|\cnum_{\lambda}\times\{O\}})},
\]
and we have a splitting
$\lefttop{\ellsitabar}F_{\veca}
\simeq
 \bigoplus_{\vecb\leq\veca}
 \nbigg_{\veca}(E)$.
The filtrations $\lefttop{i}F$ naturally induces
filtrations $\lefttop{i}F$ of $\nbigr$,
and we have an isomorphism
\[
 \nbigr_{\veca}
\simeq
 \frac{\lefttop{\ellsitabar}F_{\veca}\nbigr}
 {\sum_{\vecb\lneq\veca}\lefttop{\ellsitabar}F_{\vecb}\nbigr}.
\]
In other words, for some $N>0$,
we can order $\veca_1,\ldots,\veca_N$
and we can take a filtration $\nbigf_j$ $(j=1,\ldots,N)$ of $\nbigr$
such that $\nbigf_j/\nbigf_{j-1}\simeq \nbigr_{\veca_j}$.
Because $\nbigr_{\veca}$ are pure twistor structure of weight $0$,
by construction,
we obtain that $\nbigr$ is also a pure twistor structure of weight $0$.
\hfill\qed

\begin{prop}
\label{prop;12.12.14.41}
We have a locally free $\nbigo_{\gbigx}$-submodule
$\gbigq_0\nbige\subset\nbigqtilde_0\nbige$
such that
(i) $\gbigq_0\nbige_{|\proj^1\times\{O\}}=\nbigr$,
(ii) $\nablatilde$ is logarithmic along $\gbigx^{\infty}$
 with respect to $\gbigq_0\nbige$.
\end{prop}
\pf
It follows from the correspondence
between regular filtered flat bundles
and filtered local systems
(\cite{Simpson90}, \cite{mochi5}).
It is also obtained by successive use of the elementary transforms.
We shall give an outline of the latter argument.
We have a meromorphic flat bundle
$(\nbigqtilde\nbige,\nablatilde)$
such that
(i) $\nbigqtilde\nbige_{|\cnum_{\lambda}\times X}
=\nbigq\nbige$,
(ii) $\nablatilde$ is regular singular along $\gbigx^{\infty}$.
We have 
$\gbigx=(\cnum_{\mu}\times X)\cup(\cnum_{\lambda}\times X)$
by $\mu=\lambda^{-1}$.
We have only to construct a submodule
of $\nbigv:=\nbigqtilde\nbige_{|\cnum_{\mu}\times X}$
with the desired property.
Let $M$ be the Deligne lattice of  $(\nbigv,\nablatilde)$.
We have
$M_{|\cnum^{\ast}_{\mu}\times X}
=\nbigq_0\nbige_{|\cnum^{\ast}_{\mu}\times X}$
and 
$M_{|\cnum_{\mu}\times\{O\}}
 \otimes
 \nbigo_{\cnum_{\mu}}(\ast 0)
=\nbigqtilde_0\nbige_{|\cnum_{\mu}\times\{O\}}
=:\nbigu$.
The induced connection of $\nbigu$
is denoted by $\nabla_{\nbigu}$.
It is equal to
$\bigoplus(\lefttop{\veca}\nabla)$.

We recall a general procedure called elementary transform.
Let $Y:=\{(w_1,\ldots,w_n)\,|\,|w_i|<1\}$,
$D_Y:=\bigcup_{i=1}^{\ell}\{w_i=0\}$
and $Y_1:=\{w_1=0\}$.
Let $V$ be a holomorphic vector bundle on $Y$
with a logarithmic flat connection $\nabla$.
The holomorphic vector bundle $V_{1}:=V/w_1V$ on $Y_1$
is equipped with the endomorphism
$\Res_{w_1}(\nabla)$ obtained as the residue with respect to
$w_1$.
It is also equipped with an induced logarithmic flat connection
$\nabla_1$,
i.e.,
for a holomorphic local section $s$ of $V_{1}$,
take a holomorphic local section $\stilde$ of $V$
such that $\stilde_{|Y_1}=s$,
and $\nabla_1s=\sum_{j=2}^n \stilde_{j|Y_1}\,dw_j$,
where $\nabla \stilde=\sum_{j=1}^n\stilde_jdw_j$.
Note that $V\otimes_{\nbigo_Y}\nbigo_Y(Y_1)$
is equipped with the induced logarithmic connection.
Let $H\subset V_{1}$ be a subbundle,
which is invariant with respect to
$\Res_{w_1}(\nabla)$
and $\nabla_1$.
We obtain a holomorphic vector bundles
$V_{Hi}$ $(i=1,2)$ by the following
exact sequences of $\nbigo_Y$-modules:
\begin{equation}
\label{eq;12.12.24.1}
 0\lrarr V_{H1}\lrarr V\lrarr V_1/H\lrarr 0
\end{equation}
\begin{equation}
\label{eq;12.12.24.2}
 0\lrarr V_{H2}\lrarr V\otimes\nbigo_Y(Y_1)
 \lrarr (V_1/H)\otimes\nbigo_Y(Y_1)\lrarr 0
\end{equation}
It is easy to see that
$V_{Hi}$ have induced logarithmic flat connections
$\nabla$.
They are called the elementary transforms of $V$
along $H$
of type (\ref{eq;12.12.24.1})
and (\ref{eq;12.12.24.2}), respectively.
The kernel of the surjective morphisms
$V_{H1}/w_1V_{H1}\lrarr H$
and 
$V_{H2}/w_1V_{H2}\lrarr H\otimes_{\nbigo_Y}\nbigo_Y(Y_1)$
are denoted by $K_1$ and $K_2$.
Then, $V$ is the elementary transform of
$V_{H1}$ 
(resp. $V_{H2}$)
along $K_1$ (resp. $K_2$)
of type (\ref{eq;12.12.24.2})
(resp. (\ref{eq;12.12.24.1})).

Let us return to the proof of Proposition \ref{prop;12.12.14.41}.
We have locally free $\nbigo_{\cnum_{\mu}}$-submodules
$\nbigr_i\subset \nbigu$ $(i=0,\ldots,N)$
such that 
(i) $\nbigr_0=\nbigr$,
(ii) $\nbigr_N=M_{|\cnum_{\mu}\times\{O\}}$,
(iii) $\nbigr_i$ is the elementary transform of $\nbigr_{i-1}$
along 
$\Res(\nabla_{\nbigu})$-invariant subspaces
$H_i\subset\nbigr_{i-1}/\mu\nbigr_{i-1}$.
More precisely,
$H_i$ are direct sums of
some generalized eigen spaces of
$\Res(\nabla_{\nbigu})$.
We have the endomorphisms $F_i$ of $\nbigr$
induced by $\Res_i(\nablatilde)$.
By an induction,
$\nbigr_i$ are equipped with
induced endomorphisms $F_i$,
and $H_i$ are $F_i$-invariant.

We have either one of the morphisms
$\nbigr_i/\mu\nbigr_i\lrarr H_i$
in the case (\ref{eq;12.12.24.1})
or 
$\nbigr_i/\mu\nbigr_i\lrarr H_i
 \otimes_{\nbigo_{\cnum_{\mu}}}\nbigo_{\cnum_{\mu}}(0)$
in the case (\ref{eq;12.12.24.2}).
The kernel is denoted by $K_i$.
They are invariant with respect to $F_i$
and $\Res(\nabla_{\nbigu})$,
and we obtain $\nbigr_{i-1}$ from $\nbigr_i$
by applying an elementary transform along $K_i$
of type (\ref{eq;12.12.24.2}) or (\ref{eq;12.12.24.1}),
respectively.

We set $M_N:=M$.
We have a logarithmic connection
$\nablatilde'$ on $M'_N:=M_N/\mu M_N$.
Its residues at $(0,O)\in\cnum_{\mu}\times X$
along $z_i$ are given by $F_{i}$.
Note that $M'_N$ is the Deligne lattice.
We have the subspace
$K_N\subset \nbigr_N\simeq M'_N/(z_1,\ldots,z_n)M_N'$.
Because it is invariant with respect to $F_i$,
we have $\nablatilde'$-invariant subbundle
$\Ktilde_N'\subset M'_N$.
It is $\Res_{\mu}(\nablatilde)$-invariant.
If $\nbigr_{N}$ is the elementary transform of
$\nbigr_{N-1}$ of type (\ref{eq;12.12.24.1})
(resp. (\ref{eq;12.12.24.2})),
let $M_{N-1}$ be the elementary transform of $M_N$
along $\Ktilde_N'$ of type (\ref{eq;12.12.24.2})
(resp. (\ref{eq;12.12.24.1})).
By a descending induction with a similar procedure,
we construct $M_{i}$ $(i=0,\ldots,N)$.
Then, $M_0$ has the desired property.
\hfill\qed

\vspace{.1in}
By Lemma \ref{lem;12.12.20.20},
there exists a neighbourhood $U\subset X$ of $O$
such that,
for each $P\in U$,
the restriction
$\gbigq_0\nbige_{|\proj_{\lambda}^1\times\{P\}}$
is pure twistor structure of weight $0$.
We state it as a proposition.
\begin{prop}
\label{prop;12.12.14.3}
There exists a neighbourhood $U\subset X$ of $O$
such that
$(\nbigq_0\nbige,\nablatilde)_{|\cnum_{\lambda}\times (U\setminus D)}$
is extended to 
a tr-TLE structure
$(\gbigq_0\nbige,\nablatilde)_{|\cnum_{\lambda}\times(U\setminus D)}$
in the sense of Hertling {\rm\cite{Hertling}}.
\hfill\qed
\end{prop}

\begin{lem}
The action $L^{\nbigq_0\nbige}_{\gminivtilde}$ on
$\nbigq\nbigetilde$
is extended to an action 
$L^{\gbigq_0\nbige}_{\gminivtilde}$
on $\gbigq_0\nbigetilde$.
\end{lem}
\pf
We define
$L^{\gbigq_0\nbige}_{\gminivtilde}:=
 \nablatilde_{\gminivtilde}$.
Because $\nablatilde$ is logarithmic
along $\gbigx^{\infty}$,
we obtain that 
$L^{\gbigq_0\nbige}_{\gminivtilde}
 \gbigq_0\nbige
\subset
 \gbigq_0\nbige$.
\hfill\qed

\vspace{.1in}

Let $\beta\in\cnum$,
and let $\gminiv_1:=\beta\gminiv^{1,0}+\betabar\gminiv^{0,1}$.
We set 
$\gminivtilde_1:=
 \beta\gminivtilde^{1,0}+\betabar\gminivtilde^{0,1}$.
Let $P$ be any point of $U\setminus D$, 
and let $\gamma_t$ be an integral curve of 
$\gminiv_1$ starting from $P$
such that $\gamma_0=P$.
The holomorphic vector field induces
an isomorphism
$\proj^1_{\lambda}\times\{P\}\simeq
 \proj^1_{\lambda}\times\{\gamma_t\}$.
We have an induced isomorphism
$(\gbigq\nbige,\nablatilde)_{|\proj^1_{\lambda}\times\{P\}}
\simeq
 (\gbigq\nbige,\nablatilde)_{|\proj^1_{\lambda}\times\{\gamma_t\}}$.
(See the proof of Lemma \ref{lem;12.10.3.20} below.)
Hence, 
$\gbigq\nbige_{|\proj^1_{\lambda}\times\{\gamma_t\}}$
is also a pure twistor structure of weight $0$.
\begin{cor}
\label{cor;12.12.14.42}
In Proposition {\rm\ref{prop;12.12.14.3}},
$U$ can be an open subset in $X$
which is preserved by the flow by
$\beta\gminiv^{1,0}+\betabar\gminiv^{0,1}$
for any $\beta\in\cnum$.
\hfill\qed
\end{cor}

\subsubsection{$S^1$-action and $\cnum^{\ast}$-action (Appendix)}

We will be mainly interested in the case that
the vector field $\gminiv$ is obtained as the fundamental
vector field of an $S^1$-action.
So, as an appendix,
we give a remark
that an $S^1$-action on a holomorphic vector bundle
is extended to a $\cnum^{\ast}$-action.

Let $Y$ be a complex manifold
equipped with a holomorphic $\cnum^{\ast}$-action $\rho$.
The restriction of $\rho$ to $S^1$ is denoted by $\rho_1$.
Let $E$ be a holomorphic vector bundle on $Y$
with an $S^1$-action $\rho_{E1}$ such that
(i) $(E,\rho_{E1})$ is an $S^1$-equivariant vector bundle over $(Y,\rho_1)$,
(ii) the fundamental vector field of the action $\rho_{E1}$ is holomorphic.
In this subsection,
we shall identify $TY$ and $T^{1,0}Y$.

\begin{lem}
\label{lem;12.10.3.20}
The $S^1$-action $\rho_{E1}$ is uniquely extended to
a holomorphic $\cnum^{\ast}$-action $\rho_E$ on $E$
such that $(E,\rho_E)$ is a $\cnum^{\ast}$-vector bundle
on $(Y,\rho)$.
\end{lem}
\pf
Let $\proj_E$ denote the projective completion of $E$.
The $S^1$-action $\rho_{E1}$ on $E$ is uniquely extended 
to an $S^1$-action of $\rho_{E2}$ on $\proj_E$,
whose fundamental vector field $V$ is also holomorphic.
We obtain the holomorphic vector field $\sqrt{-1}V$,
and we have $[V,\sqrt{-1}V]=0$.

Let $\nu$ be the $\real$-action on $Y$
given by $\nu(t,Q)=\rho(e^{-t},Q)$.
Let $V_1$ denote the fundamental vector field of $\nu$.
Let $\pi:\proj_E\lrarr Y$ denote the projection.
Let $T_P\pi:T_P\proj_E\lrarr T_{\pi(P)}Y$
denote its derivative.
Then, $T_{P}\pi(\sqrt{-1}V_{|P})=V_{1|Q}$
for any $P\in\pi^{-1}(Q)$.

Assume that $V_{1|Q}\neq 0$.
By using the properness of $\pi$,
there exist $\epsilon>0$
and $C^{\infty}$-map
$F:(-\epsilon,\epsilon)\times \pi^{-1}(Q)
\lrarr \proj_E$
such that
$F(t,P)$ is an integral curve $\gamma(t)$
of $\sqrt{-1}V$ such that $\gamma(0)=P$.
For each $t\in(-\epsilon,\epsilon)$,
$F(t,\bullet)$ gives a diffeomorphism
of $\pi^{-1}(Q)\simeq \pi^{-1}(\nu(t,Q))$.
Then, it is easy to prove that
an integral curve of $\sqrt{-1}V$ through any point of $\proj_E$
can be defined over $\real$.
Hence, we have an $\real$-action $\nu_{\proj_E}$
on $\proj_E$ whose fundamental vector field is 
$\sqrt{-1}V$.
Because it preserves $\proj_E\setminus E$,
we obtain an $\real$-action $\nu_E$
on $E$.
It commutes with $\rho_{E1}$,
and satisfies $\phi\circ\nu_E(t,P)=\nu(t,\pi_1(P))$,
where $\pi_1:E\lrarr Y$.
Hence, we obtain a holomorphic $\cnum^{\ast}$-action
$\rho_E$ such that 
$\pi_1$ is equivariant with respect to
$\rho_E$ and $\rho$.
Let $\alpha_1,\alpha_2\in\cnum$ 
and $v_1,v_2\in E_{|Q}$ for $Q\in Y$.
We have
$\rho_E(t,\alpha_1 v_1+\alpha_2 v_2)
=\alpha_1\rho_E(t,v_1)+\alpha_2\rho_E(t,v_2)$
for $t\in S^1$.
Hence, it holds for $t\in \cnum^{\ast}$,
i.e., $(E,\rho_E)$ is a $\cnum^{\ast}$-equivariant
bundle over $(Y,\rho)$.
\hfill\qed

\section{$\frac{\infty}{2}$-VHS associated to
 Toda like harmonic bundle}
\label{section;12.12.14.200}
\subsection{Explicit description of 
the associated meromorphic flat  bundle}
\label{subsection;12.12.16.11}

\subsubsection{Refinement for meromorphic prolongment}
\label{subsection;12.12.20.1}

We set $X:=\proj^1_q$, $D_0:=\{0\}$, $D_{\infty}:=\{\infty\}$
and $D:=D_0\cup D_{\infty}$.
We set $\gminiv:=\sqrt{-1}r(q\del_q-\qbar\del_{\qbar})$,
which is the fundamental vector field of 
the $S^1$-action on $\proj^1$
given by $(t,q)\longmapsto t^rq$.
We set $\nbigx:=\cnum_{\lambda}\times X$
and $\gbigx:=\proj^1\times X$.
We use the notation $\nbigd$, $\gbigd$, etc.,
in similar meanings.
We set $\gbigx^{\lambda}:=\{\lambda\}\times X\subset\gbigx$
for $\lambda\in\proj^1$.

Because $t^{\ast}\theta=t^m\theta$ and $t^{\ast}h=h$,
the harmonic bundle
$(E,\theta,h):=(E_r,\theta_{r,m},h_{\veca})$ on $X\setminus D$
is homogeneous of degree $m$
with respect to $\gminiv$.
As explained in \S\ref{subsection;12.12.13.10} and
\S\ref{subsection;12.12.14.31},
we have the locally free $\nbigo_{\nbigx}(\ast\nbigd)$-module
$\nbigq\nbige$
with the meromorphic flat connection $\nablatilde$.
We have the $\nbigo_{\nbigx}(\ast \nbigd_{\infty})$-submodule
$\nbigq_{0}\nbige\subset\nbigq\nbige$.
The connection $\nablatilde$ is logarithmic
along $\nbigd_0$.

\begin{lem}
\label{lem;12.10.3.30}
$\nbigq\nbige$
is uniquely extended to
a locally free $\nbigo_{\gbigx}
 \bigl(\ast(\gbigd\cup\gbigx^{\infty})\bigr)$-module
$\nbigqtilde\nbige$
such that $\nablatilde$ gives a meromorphic
flat connection of $\nbigqtilde\nbige$
which is regular singular at any point of
$\gbigx^{\infty}\setminus\gbigd_{\infty}$.
\end{lem}
\pf
Note that the harmonic bundle is not tame at $\infty$.
We have only to consider the issue
locally around $(\lambda,q)=(\infty,\infty)$.
We use the Riemann-Hilbert correspondence
between meromorphic flat bundles 
and local systems with Stokes structure,
studied in \S4 of \cite{Mochizuki-wild}.
See \cite{Mochizuki-Stokes-Review}
for a review.

By using the descent and pull back,
we may assume $m=r$.
Then, the eigenvalues of
$\theta_{r,r}$ are of the form
$\alpha\,dq$ $(\alpha\in\cnum)$.
Hence, the irregular values of
$(\nbigq\nbige,\nablatilde)$
along $q=\infty$ are of the form
$\alpha \lambda^{-1} q$ 
$(\alpha\in\cnum)$.

It is clear that $\nbigq\nbige$
is extended
along $\gbigx^{\infty}\setminus\{(\infty,\infty)\}$
in a regular singular way.
We take a blow up
$\varphi:Z\lrarr \gbigx$
at $(\infty,\infty)$.
Then, the pole and the zero of the function
$\varphi^{-1}(\lambda^{-1}q)$
are separated.
The pole of $\varphi^{-1}(\lambda^{-1}q)$
is the strict transform of 
$\proj_{\lambda}^1\times\{\infty\}$.

We may restrict our interest to a neighbourhood of 
$\varphi^{-1}(\infty,\infty)$.
Let $D_1'$ be the exceptional divisor.
Let $\pi:\Ztilde\lrarr Z$ be the real blow up
along 
$\varphi^{-1}\bigl(\proj_{\lambda}^1\times\{\infty\}\bigr)
 \cup\varphi^{-1}(\gbigx^{\infty})$.
Let $\nbigl$ be the local system corresponding to
$(\nbigq\nbige,\nablatilde)$.
It is equipped with the Stokes structure
outside $\pi^{-1}(D_1')$.
(See \cite{Mochizuki-Stokes-Review}
for a review of the Stokes structure.
We shall use the notation there.)

Let $D_2'$ be the strict transform of
$\proj_{\lambda}^1\times\{\infty\}$.
Let $P$ be the intersection of $D_1'$ and $D_2'$.
Let $Q\in\pi^{-1}(P)$.
For a local coordinate $(u,v)$ around $P$
such that $q=uv$ and $\mu=v$,
the irregular values are of the form
$\alpha\,u^{-1}$ $(\alpha\in\cnum)$.
We can take a point 
$Q'\in\pi^{-1}(D'_2\setminus P)$ such that
(i) $Q'$ is sufficiently close to $Q$,
(ii) the orders $\leq_{Q'}$ and $\leq_{Q}$
are the same on the set of the irregular values of
$(\nbigq\nbige,\nablatilde)$.
We have a natural isomorphism
$\nbigl_{Q}\simeq\nbigl_{Q'}$.
Hence, we obtain a filtration of $\nbigl_Q$
induced by the Stokes filtration of $\nbigl_{Q'}$.
Thus, we obtain the Stokes structure of $\nbigl$.
By using the Riemann-Hilbert correspondence,
we obtain a good meromorphic flat bundle on $Z$.
Its push-forward to $\gbigx$
gives the desired extension 
of $(\nbigq\nbige,\nablatilde)$.
\hfill\qed

\begin{rem}
The first claim of the previous lemma holds
in a more general situation,
which will be argued elsewhere.
\hfill\qed
\end{rem}

We have a locally free
$\nbigo_{\nbigx}(\ast\nbigd_{\infty})$-module
$\nbigq_0\nbige\subset\nbigq\nbige$.
\begin{lem}
It is extended to a locally free
$\nbigo_{\gbigx}\bigl(
 \ast(\gbigd_{\infty}\cup\gbigx^{\infty})\bigr)$-module
$\nbigqtilde_0\nbige$.
\end{lem}
\pf
We have only to check the claim around
$(\lambda,q)=(\infty,0)$.
We have the Deligne lattice $M$
of $(\nbigqtilde\nbige,\nablatilde)$ around $(\infty,0)$.
Then, it is easy to see that
$\nbigqtilde_0\nbige
=M(\ast \gbigx^{\infty})$
around $(\infty,0)$.
\hfill\qed

\vspace{.1in}

Let $\gminiv_m:=m^{-1}\gminiv$.
The harmonic bundle
$(E,\theta,h)$ is homogeneous of degree $1$
with respect to $\gminiv_m$.
Even if we replace $\gminiv$ with $\gminiv_m$,
we clearly obtain the same meromorphic flat connection
$\nablatilde$.

For a positive integer $\ell$, 
let $\varphi_{\ell}:\proj^1\lrarr\proj^1$
be given by $\varphi(q)=q^{\ell}$.
It gives $\varphi_{\ell\ast}(\gminiv_{\ell m})=\gminiv_m$.
The harmonic bundle
$\varphi_{\ell}^{-1}(E,\theta,h)$
is homogeneous of degree $\ell m$
with respect to $\gminiv$.
Then, it is easy to observe that
the meromorphic flat bundle associated to
$\varphi_{\ell}^{-1}(E,\theta,h)$
is naturally isomorphic to
$\varphi_{\ell}^{\ast}(\nbigq\nbige,\nablatilde)$,
where the induced map
$\nbigx\lrarr\nbigx$, $(\lambda,q)\longmapsto (\lambda,q^{\ell})$
is also denoted by $\varphi_{\ell}$.

\subsubsection{Refinement for tr-TLE structure}

\begin{prop}
\label{prop;12.10.6.20}
We have a locally free 
$\nbigo_{\gbigx}(\ast\gbigd_{\infty})$-submodule
$\gbigq_0\nbige\subset
 \nbigqtilde\nbige$
such that 
$(\gbigq_0\nbige,\nablatilde)$
is tr-TLE structure in the sense of Hertling {\rm\cite{Hertling}},
i.e.,
(i) $\nablatilde$ is logarithmic with respect to $\gbigq_0\nbige$
along $\gbigx^{\infty}$,
(ii)  $\gbigq_0\nbige_{|\proj^1\times\{Q\}}$
is a pure twistor structure of weight $0$
for each $Q\in X\setminus D_{\infty}$.
\end{prop}
\pf
We have already studied the issue 
locally around $q=0$ in Proposition \ref{prop;12.12.14.41},
by using the Deligne lattice.
For the construction of $\gbigq_0\nbige$
in this proposition,
we have only to replace the Deligne lattice
with $M\otimes\nbigo_{\gbigx}(\ast\gbigd_{\infty})$,
where $M$ is the Deligne-Malgrange lattice of
$\nbigqtilde\nbige$.
(See \cite{malgrange} for Deligne-Malgrange lattice,
or canonical lattice.
 See \cite{mochi9} for a review.)
As explained in Proposition \ref{prop;12.12.14.3}
and Corollary \ref{cor;12.12.14.42},
it is a tr-TLE structure.
\hfill\qed

\vspace{.1in}
We set
$\gbigq\nbige:=
 \gbigq\nbige\otimes
 \nbigo_{\gbigx}(\ast\gbigd_0)$.
Note that we have a natural isomorphism
$\gbigq\nbige_{|\gbigx^0}
\simeq
 \nbigv_{r}$
by construction.
We shall identify them.
Let $p:\gbigx\lrarr X$ denote the projection.
\begin{prop}
\label{prop;12.12.19.1}
We have a unique holomorphic isomorphism
$\Phi:
 p^{\ast}\nbigv_{r}
\simeq
 \gbigq\nbige$
such that
$\Phi_{|\lambda=0}$ is the identity.
\end{prop}
\pf
We have a natural isomorphism
$\gbigq\nbige\simeq
 p^{\ast}
 p_{\ast}\gbigq\nbige$.
Because
$\gbigq\nbige_{|\{0\}\times X}
=
 \nbigv_{r}$,
we have a natural morphism
$p_{\ast}\gbigq\nbige\lrarr
 \nbigv_{r}$.
Because its restriction to $X^{\ast}$
is an isomorphism,
it is an isomorphism on $X$.
Hence, we obtain a natural isomorphism
$p^{\ast}\nbigv_{r}
\simeq
 \gbigq\nbige$.
\hfill\qed

\vspace{.1in}

Under the isomorphism $\Phi$,
we have
$\nablatilde=
 p^{-1}\nabla_1+\lambda^{-1}\theta
 +\nablatilde_{\lambda}d\lambda$
for some connection $\nabla_1$ of 
$\nbigv_{r}$.
Recall the decomposition
$\nbigv_{r}=\bigoplus_{i=1}^r \nbigv^{(i)}$,
where $\nbigv^{(i)}=\nbigo_{X}(\ast D)e_i$.

\begin{prop}
\label{prop;12.12.4.11}
$\nabla_1$ is isomorphic to 
$\nabla_0-\bigoplus a_i\id_{\nbigv^{(i)}}dq/q$,
where $\nabla_0e_i=0$ for $i=1,\ldots,r$.
\end{prop}
\pf
We have only to consider the case $m=1$
by the remark given in the last of \S\ref{subsection;12.12.20.1}.
Let $\sigma_0$ be an automorphism of $\gbigx$
given by $\sigma_0(\lambda,q)=(\tau\lambda, q)$,
where $\tau$ is the primitive $r$-the root.
We have the $\cnum^{\ast}$-action
on $\nbigx$
given by $t(\lambda,q)=(t\lambda,t^{r}q)$.
The automorphism $\sigma$ of $\nbigv_r$
induces an automorphism $\sigma$ of $p^{\ast}\nbigv_r$
over $\sigma_0$.
The $\cnum^{\ast}$-action on $\nbigv_r$
induces the action on $p^{\ast}\nbigv_r$.
Then, $\lambda^{-1}\theta$
and $\nablatilde$ are $\sigma$-equivariant
and $\cnum^{\ast}$-equivariant.
Hence, $\nabla_1$ is also $\sigma$-equivariant 
and $\cnum^{\ast}$-equivariant.
Let $\nbigv_{r,i}:=\nbigo_{\proj^1}(\ast D)e_i$.
Then, $\nabla_1$ preserves the decomposition
$\nbigv_{r}=\bigoplus \nbigv_{r,i}$,
and it is expressed 
as $\nabla_0-\bigoplus b_i\id_{\nbigv_{r,i}}dq/q$,
where $b_i\in\cnum$.
The numbers $b_i$ must be equal to
the parabolic weights.
(Recall the dependence of the parabolic weights and 
the eigenvalues of the residues $\Res(\nabla^{\lambda})$
on $\lambda$
in \cite{Hertling-Sevenheck}, \cite{mochi2}, \cite{Simpson90}.)
\hfill\qed

\subsubsection{Explicit description}
\label{subsection;12.12.4.100}

Let us give an explicit description of $\nablatilde$
under the isomorphism
$\gbigq\nbige
\simeq
 p^{\ast}\nbigv_r$
in Proposition \ref{prop;12.12.19.1}.
For given complex numbers
$\gamma_1,\ldots,\gamma_r$,
let $\diag[\gamma_1,\ldots,\gamma_r]$
denote the diagonal matrix
whose $(j,j)$-entries are $\gamma_j$.
By the construction of $\cnum^{\ast}$-action
and $\nablatilde$,
we have
$\nablatilde_{m\lambda\del_{\lambda}+rq\del_q}
 \vece
=\vece\,\diag[m,2m,\ldots,rm]$.
According to Proposition \ref{prop;12.12.4.11},
we have
\begin{equation}
 \label{eq;12.12.4.12}
 q\nablatilde_q\vece
=\vece\Bigl(
-\diag[a_1,\ldots,a_r]
+\frac{m}{\lambda} \nbigk(r,m)
 \Bigr).
\end{equation}
Hence, we obtain
\begin{equation}
 \label{eq;12.12.4.13}
 \lambda\nablatilde_{\lambda}\vece
=\vece\Bigl(
 \diag[1,2,\ldots,r]
+\frac{r}{m}
 \diag[a_1,\ldots,a_r]
-\frac{r}{\lambda}
 \nbigk(r,m)
 \Bigr).
\end{equation}
The formulas (\ref{eq;12.12.4.12}) and (\ref{eq;12.12.4.13})
describe $\nablatilde$.

In particular, if $m=r$,
we have the following:
\begin{equation}
 \label{eq;12.12.4.20}
  q\nablatilde_q\vece
=\vece\Bigl(
-\diag[a_1,\ldots,a_r]
+\frac{r}{\lambda} \nbigk(r,r)
 \Bigr)
\end{equation}
\begin{equation}
\label{eq;12.12.4.21}
 \lambda\nablatilde_{\lambda}\vece
=\vece\Bigl(
 \diag[1,2,\ldots,r]
+\diag[a_1,\ldots,a_r]
-\frac{r}{\lambda}
 \nbigk(r,r)
 \Bigr)
\end{equation}
Let $C_r$ be the matrix whose $(i,j)$-th entry is
$1$ if $i\equiv j+1$ modulo $r$,
or $0$ otherwise.
Then, by specializing (\ref{eq;12.12.4.21}) 
to $q=1$,
we obtain
\[
  \lambda\nablatilde_{\lambda}\vece_{|q=1}
=\vece_{|q=1}\Bigl(
 \diag[1,2,\ldots,r]
+\diag[a_1,\ldots,a_r]
-\frac{r}{\lambda}C_r
 \Bigr).
\]
By specializing (\ref{eq;12.12.4.20}) to $\lambda=1$,
we obtain
\[
  q\nablatilde_q\vece_{|\lambda=1}
=\vece_{|\lambda=1}
\Bigl(
-\diag[a_1,\ldots,a_r]
+r\nbigk(r,r)
 \Bigr). 
\]
We set $v_i:=q^{-i}e_i$.
Then, we obtain
\[
 q\nablatilde_q\vecv
=\vecv\Bigl(
-\diag[1,\ldots,r]-\diag[a_1,\ldots,a_r]
+rC_rq
 \Bigr)
\]
If we set $\zeta=q^{-1}$,
we obtain the following:
\[
 \zeta\nablatilde_{\zeta}\vecv
=\vecv\Bigl(
 \diag[1,\ldots,r]+\diag[a_1,\ldots,a_r]
-\frac{r}{\zeta}C_r
 \Bigr)
\]
In particular,
we obtain the following:
\begin{cor}
\label{cor;12.12.4.22}
We have an isomorphism
$(\nbigqtilde\nbige,\nablatilde)_{|\gbigx^1}$
and
$(\nbigqtilde\nbige,\nablatilde)_{|\proj_{\lambda}^1\times\{1\}}$.
\hfill\qed
\end{cor}

\begin{rem}
We can also obtain the isomorphism in
Corollary {\rm\ref{cor;12.12.4.22}}
by using the homogeneity
without explicit computations.
\hfill\qed
\end{rem}

\subsection{Integral structure of Toda-like harmonic bundle}
\label{subsection;12.12.16.12}

\subsubsection{Special case}

Let $(E_r,\theta_{r,1},h)$ be the Toda-like harmonic bundle,
corresponding to
$\veca=(a_i)\in\gbigr_{r,1}(a)$
as in Theorem \ref{thm;12.9.24.200}.
We study when
$(\nbigqtilde\nbige,\nablatilde)_{|\proj^1_{\lambda}\times\{1\}}$
and
$(\nbigqtilde\nbige,\nablatilde)_{|\gbigx^1}$
have $\seisuu$-structures.

\begin{thm}
\label{thm;12.12.4.30}
$(\nbigqtilde\nbige,\nablatilde)_{|\gbigx^1}$
has a $\seisuu$-structure
(resp. $\rnum$-structure)
if and only if
$P_{\veca}(T)\in \seisuu[T]$
(resp. $\rnum$-structure).
(See Theorem {\rm\ref{thm;12.12.22.300}} below 
for general $m$.)
\end{thm}
\pf
It follows from Theorem \ref{thm;12.10.6.40},
Corollary \ref{cor;12.12.4.33}
and the formula (\ref{eq;12.12.4.12}).
\hfill\qed

\vspace{.1in}
If $(\nbigqtilde\nbige,\nablatilde)_{|\gbigx^1}$
has a $\seisuu$-structure,
then we have $\veca\in\rnum^r$.
Indeed, because we have
$\prod_{i=1}^{r}(T-e^{2\pi\sqrt{-1}a_i})\in\seisuu[T]$,
the number $e^{2\pi\sqrt{-1}a_i}$ and its conjugates over $\rnum$
are algebraic integers such that the absolute values are $1$.
It is a classical theorem of Kronecker that
$e^{2\pi\sqrt{-1}a_i}$ are roots of unity,
i.e., $a_i\in\rnum$.
Hence, when we are interested in $\seisuu$-structure,
it is natural to impose the condition $\veca\in\rnum^r$. 
In that case, 
we have $P_{\veca}(T)\in\seisuu[T]$
if and only if $P_{\veca}(T)\in\rnum[T]$,
and therefore
$(\nbigqtilde\nbige,\nablatilde)_{|\gbigx^1}$
has a $\seisuu$-structure
if and only if it has $\rnum$-structure.

\begin{rem}
The condition $P_{\veca}(T)\in \seisuu[T]$
implies an additional symmetry for
the sequence $(a_1,\ldots,a_n)$.
Hence, the harmonic bundle
is equipped with a symmetric pairing.
(See {\rm\S\ref{subsection;12.12.21.101}}.)
It induces a real structure by the procedure in
{\rm\S\ref{subsection;12.12.23.100}}.
By Proposition {\rm\ref{prop;12.12.4.101}},
we can adjust that a $\seisuu$-structure
is compatible with the $\real$-structure.

\hfill\qed
\end{rem}

\begin{rem}
If the condition in Theorem {\rm\ref{thm;12.12.4.30}}
is satisfied,
the meromorphic flat bundle
$(\nbigqtilde\nbige,\nablatilde)$ has a $\seisuu$-structure,
although we do not discuss a $\seisuu$-structure
of a meromorphic flat bundle on higher dimensional varieties
in this paper.
In particular,
$(\nbigqtilde\nbige,\nablatilde)_{|\proj^1_{\lambda}\times\{1\}}$
has a $\seisuu$-structure,
which also follows from Theorem 
{\rm\ref{thm;12.12.4.130}} below.
\hfill\qed
\end{rem}

\subsubsection{Criterion}

We shall freely use the notation in \S\ref{subsection;12.12.4.31}.
Let $h$ be the Toda-like harmonic metric of 
$(E,\theta):=(E_r,\theta_{r,r})$
corresponding to
$\veca\in\gbigr_{r,r}\cap\rnum^r$.
We have the associated meromorphic flat bundle
$(\nbigqtilde\nbige,\nablatilde)$ on $\gbigx$.
Note that
$(\nbigqtilde\nbige,\nablatilde)_{|\gbigx^1}$
has a $\seisuu$-structure
(resp. $\rnum$-structure)
if and only if
$(\nbigqtilde\nbige,\nablatilde)_{|\proj^1_{\lambda}\times\{1\}}$
has a $\seisuu$-structure
(resp. $\rnum$-structure)
by Corollary \ref{cor;12.12.4.22}.
We set $\vecdelta_d=(\overbrace{1,\ldots,1}^d)$
for a positive integer $d$.

\begin{thm}
\label{thm;12.12.4.130}
The following conditions are equivalent.
\begin{description}
\item[(P1)]
 $(\nbigqtilde\nbige,\nablatilde)_{|\proj^1_{\lambda}\times\{1\}}$
 has a $\seisuu$-structure.
\item[(P2)]
 $(\nbigqtilde\nbige,\nablatilde)_{|\proj^1_{\lambda}\times\{1\}}$
 has a $\rnum$-structure.
\item[(P3)]
 There exists $\gamma\in\cnum^{\ast}$ such that 
 (i) $\gamma^r\in\rnum$,
 (ii) $\prod_{i=1}^r(T-\gamma e^{2\pi\sqrt{-1}a_i/r})\in\rnum[T]$.
\end{description}
\end{thm}
\pf
We shall use the description of 
$(\nbigqtilde\nbige,\nablatilde)$
given in \S\ref{subsection;12.12.4.100}.
We have only to study the existence 
of $\rnum$-structure or $\seisuu$-structure
of $(\nbigqtilde\nbige,\nablatilde)_{|\gbigx^1}$.

Let us consider the case
$\nbigs(\veca,r)=\{r\}$.
By Proposition \ref{prop;12.12.4.110}
and the Kobayashi-Hitchin correspondence,
the filtered flat bundle
$(\nbigq_{\ast}\nbige,\nablatilde)_{|\gbigx^1}$
is stable.
Because the eigenvalues of the residue of  $\theta$ are $0$
at both $\{0\}$ and $\{\infty\}$,
we obtain that
$(\nbigqtilde\nbige,\nablatilde)_{|\gbigx^1}$
is irreducible
as a meromorphic flat bundle.
(See \cite{sabbah3}. 
See also Lemma 19.4.3 of \cite{Mochizuki-wild}.
Note that $(\nbigqtilde\nbige,\nablatilde)_{|\gbigx^1}$
is deformed as in \S11 \cite{Mochizuki-wild}.)
Then, the equivalence of the three conditions
follow from Theorem \ref{thm;12.12.26.12}.

\vspace{.1in}
Let us consider the case
$\nbigs(\veca,r)\neq\{r\}$.
We have the decomposition
$(\nbigqtilde\nbige,\nablatilde)
=\bigoplus_s
 \bigl(\nbigqtilde\nbige^{(s)},\nablatilde^{(s)}\bigr)$
induced by (\ref{eq;12.12.4.120}).
As remarked in the argument for the case
$\nbigs(\veca,r)=\{r\}$,
each $(\nbigqtilde\nbige^{(s)},\nablatilde^{(s)})_{|\gbigx^1}$
is irreducible.
Because there is no common eigenvalue
of $\theta^{(s)}$ and $\theta^{(s')}$ $(s\neq s')$,
we have 
$\Irr(\nablatilde^{(s)})\cap\Irr(\nablatilde^{(s')})=\emptyset$
if $s\neq s'$.
Let us show the equivalence of
the conditions {\bf (P1)} and {\bf (P2)}.
Clearly {\bf (P1)} implies {\bf (P2)}.
Suppose {\bf (P2)}.
By applying Proposition \ref{prop;12.12.4.131}
to the irreducible decomposition
$(\nbigqtilde\nbige,\nablatilde)_{|\gbigx^1}
=\bigoplus_s
 \bigl(\nbigqtilde\nbige^{(s)},\nablatilde^{(s)}\bigr)
 _{|\gbigx^1}$,
each $(\nbigqtilde\nbige^{(s)},\nablatilde^{(s)})_{|\gbigx^1}$
has a $\rnum$-structure.
By Theorem \ref{thm;12.12.26.12},
each $(\nbigqtilde\nbige^{(s)},\nablatilde^{(s)})_{|\gbigx^1}$
has a $\seisuu$-structure.
Hence, {\bf (P1)} holds.

We consider the following condition:
\begin{description}
\item[(P3')]
  There exists $\gamma\in\cnum^{\ast}$ such that 
(i) $\gamma^{r_0}\in\rnum$,
(ii) $\prod_{i=1}^{r_0}(T-\gamma e^{2\pi\sqrt{-1}a_i/r_0})
 \in\rnum[T]$.
\end{description}
If {\bf (P3')} holds,
by applying the result in the case $\nbigs(\veca,r)=\{r\}$
to each  $(\nbigqtilde\nbige^{(s)},\nablatilde^{(s)})_{|\gbigx^1}$,
we obtain that {\bf (P2)} holds.
Conversely, if {\bf (P2)} holds,
by Proposition \ref{prop;12.12.4.131},
each 
$\bigl(\nbigqtilde\nbige^{(s)},\nablatilde^{(s)}\bigr)
 _{|\gbigx^1}$
has a $\rnum$-structure.
By the result in the case 
$\nbigs(\veca,r)=\{r\}$,
we obtain that {\bf (P3')} holds.
Finally, the claim of the theorem follows from
the equivalence of 
the conditions {\bf (P3)} and {\bf (P3')}.
(See Lemma \ref{lem;12.12.12.20.2} below.)
\hfill\qed

\vspace{.1in}
For example, we obtain the following 
from Theorem \ref{thm;12.12.4.130}
and Lemma \ref{lem;12.12.27.10}
(or Proposition {\rm\ref{prop;12.12.28.1}}).
\begin{cor}
Suppose that $r$ is an odd prime.
Then,
$(\nbigqtilde\nbige,\nablatilde)_{|\proj^1_{\lambda}\times\{1\}}$
has a $\seisuu$-structure,
if and only if
there exists a half integer $\ell$ such that
$\prod_{i=1}^r(T-e^{2\pi\sqrt{-1}(a_i+\ell)/r})\in
 \seisuu[T]$.
\hfill\qed
\end{cor}

\begin{rem}
According to the correspondence
between Toda-like harmonic bundles
and solutions of Toda lattice with opposite sign
(see \S\ref{subsection;13.1.8.2}),
we obtain a criterion when 
the meromorphic flat bundle
associated to a solution of the Toda equation
has a $\seisuu$-structure.
\hfill\qed
\end{rem}

\paragraph{Variant}
Let $h$ be the Toda-like harmonic metric of $(E_r,\theta_{r,m})$
corresponding to
$\veca\in\gbigr_{r,m}\cap\rnum^r$.
The following proposition can be proved 
by the argument in the proof of Theorem \ref{thm;12.12.4.130}.

\begin{thm}
\label{thm;12.12.22.300}
The associated meromorphic flat bundle
 $(\nbigqtilde\nbige,\nablatilde)_{|\gbigx^1}$
 has a $\rnum$-structure,
if and only if
there exists
$\gamma\in\cnum^{\ast}$ such that 
(i) $\gamma^m\in \rnum$,
(ii) $\prod_{i=1}^r(T-\gamma e^{2\pi\sqrt{-1}a_i/m})
 \in\rnum[T]$.
In that case,
 $(\nbigqtilde\nbige,\nablatilde)_{|\gbigx^1}$
 has a $\seisuu$-structure, more strongly.
\hfill\qed
\end{thm}

\paragraph{Appendix}
Let $r_0,j_0$ be positive integers,
and set $r_1=r_0j_0$.
Let $\mu_{p}:=\bigl\{\kappa\in\cnum\,\big|\,
 \kappa^{p}=1
 \bigr\}$ for a positive integer $p$.
We have the homomorphism
$\Psi:\mu_{r_1}\lrarr\mu_{r_0}$
given by $\Psi(\kappa)=\kappa^{j_0}$.
Let $S_0$ be a subset $\mu_{r_0}$
with a function $f_0:S_0\lrarr\seisuu_{>0}$.
We put $S_1:=\Psi^{-1}(S_0)$
and $f_1:=f_0\circ\Psi$.
We take $\gamma_1\in\cnum^{\ast}$
and put $\gamma_0:=\gamma_1^{j_0}$.
For $i=0,1$, we set
\[
 P_i(T):=\prod_{b\in S_i}(T-\gamma_ib)^{f_i(b)}.
\]
Because $P_1(T)=P_0(T^{j_0})$,
the following lemma is clear.

\begin{lem}
\label{lem;12.12.12.20.2}
We have
$P_0(T)\in\rnum[T]$
if and only if 
$P_1(T)\in\rnum[T]$.
\hfill\qed
\end{lem}

\noindent
{\em Address\\
Research Institute for Mathematical Sciences,
Kyoto University,
Kyoto 606-8502, Japan,\\
takuro@kurims.kyoto-u.ac.jp
}

\end{document}